%% file: Elimination_cusp_dim4.tex
\documentclass[a4paper,reqno]{amsart}
	\pdfoutput=1

	\usepackage[utf8]{inputenc}
	\usepackage{xspace}
	\usepackage{hyperref}

	\usepackage{amssymb,amsmath,amsthm}
	\usepackage[alphabetic,abbrev]{amsrefs}

	\usepackage{graphicx}
	\usepackage[all]{xy}
	\usepackage{color}
	\usepackage{transparent}
	\usepackage{enumerate}
	\usepackage{amscd}
	\usepackage{subfigure}
	\usepackage{thmtools, thm-restate}

	\usepackage[nameinlink,capitalise,noabbrev]{cleveref}
		\Crefname{page}{page}{pages}

\input{customcommands}

\input{morecommands_SB}

\input{morecommands_KH}

\begin{document}

\title{Elimination~of~cusps~in~dimension~4 and~its~applications}

\author{Stefan Behrens}
\address{%
	Alfréd Rényi Institute of Mathematics, %
	Hungarian Academy of Sciences, %
	Reáltanoda utca 13--15, %
	H--1053 Budapest, %
	Hungary}
\email{sbehrens@renyi.hu}

\author{Kenta Hayano}
\address{%
Department of Mathematics, %
Graduate School of Science, %
Hokkaido University, %
Kita 10, Nishi 8, Kita-ku, %
Sapporo, Hokkaido 060-0810, Japan}
\email{k-hayano@math.sci.hokudai.ac.jp}


\begin{abstract}

Several new combinatorial descriptions of closed 4--manifolds have recently been introduced in the study of smooth maps from 4--manifolds to surfaces. 
These descriptions consist of simple closed curves in a closed, orientable surface and these curves appear as so called \emph{vanishing sets} of corresponding maps.
In the present paper we focus on homotopies canceling pairs of cusps so called \emph{cusp merges}. 
We first discuss the classification problem of such homotopies, showing that there is a one-to-one correspondence between the set of homotopy classes of cusp merges canceling a given pair of cusps and the set of homotopy classes of suitably decorated curves between the cusps. 
Using our classification, we further give a complete description of the behavior of vanishing sets under cusp merges in terms of mapping class groups of surfaces. 
As an application, we discuss the uniqueness of \emph{surface diagrams}, which are combinatorial descriptions of 4--manifolds due to Williams, and give new examples of surface diagrams related with Lefschetz fibrations and Heegaard diagrams. 
\end{abstract}

\maketitle


\input{intro}
\input{preliminaries}

\input{merging}
\input{parallel_transport_short}
\input{MCG_interpretation}
\input{applications}
\input{applications_moves}

\input{applications_LF_HD}

\appendix
\renewcommand*{\thesection}{\Roman{section}}
\input{appendix_normalization_new}

\input{appendix_symmetries}


\vspace{1em}

\noindent
{\bf Acknowledgments.}
During the work on this project the first author was supported by an IMPRS Scholarship of the Max Planck Institute for Mathematics in Bonn and the ERC grant LDTBud.
The second author was supported by JSPS Research Fellowships for Young Scientists~(24$\cdot$993) and Grant-in-Aid for Young scientists (B)~(26800027). 
The authors would like to thank the Max Planck Institute for hospitality and Carlos Moraga Ferrándiz for helpful comments on an earlier version of this manuscript.


\input{references}
\end{document}

%% file: customcommands.tex
\usepackage{amsthm}

	\Crefname{maintheorem}{Theorem}{Theorems}
\newtheorem{theorem}{Theorem}[section]
\newtheorem{proposition}[theorem]{Proposition}
\newtheorem{lemma}[theorem]{Lemma}

\newtheorem{fact}[theorem]{Fact}

\theoremstyle{definition}
\newtheorem{definition}[theorem]{Definition}

\newtheorem{example}[theorem]{Example}

\theoremstyle{remark}
\newtheorem{remark}[theorem]{Remark}

\newcommand{\eqand}{\quad\text{and}\quad}
\newcommand{\eps}{\varepsilon}

\newcommand{\Z}{\mathbb{Z}}

\newcommand{\R}{\mathbb{R}}
\newcommand{\C}{\mathbb{C}}

\newcommand{\ra}{\rightarrow}

\newcommand{\lra}{\longrightarrow}

\DeclareMathOperator{\coker}{coker}	
\DeclareMathOperator{\rk}{rk}	
\DeclareMathOperator{\im}{im}	

\DeclareMathOperator{\Diff}{Diff}	
\DeclareMathOperator{\id}{id}	
\DeclareMathOperator{\Crit}{Crit}	

\newcommand{\set}[1]{\left\{ #1 \right\}}
\newcommand{\Set}[2]{ \left\{ {#1} \,\middle|\, {#2} \right\} }
\newcommand{\scp}[1]{\left\langle { #1 } \right\rangle}

\newcommand{\inv}{^{-1}}

\newcommand{\del}{\partial}

\newcommand{\bs}{\boldsymbol}

\usepackage{xspace}

\newcommand{\wrt}{with respect to\xspace}
\newcommand{\rhs}{right hand side\xspace}

\newcommand{\nbhd}{neighborhood\xspace}
\newcommand{\nbhds}{neighborhoods\xspace}

%% file: morecommands_SB.tex
\newcommand{\Cinfty}{C^\infty}

\newcommand{\wf}{wrinkled fibration\xspace}

\newcommand{\wfs}{wrinkled fibrations\xspace}
\newcommand{\Wfs}{Wrinkled fibrations\xspace}
\newcommand{\swf}{simple wrinkled fibration\xspace}

\newcommand{\swfs}{simple wrinkled fibrations\xspace}

\newcommand{\mc}{\mathcal}
\newcommand{\mf}{\mathfrak}
\newcommand{\wt}[1]{\widetilde{#1}}

\newcommand{\FM}{\mathfrak{FM}} 
\newcommand{\CM}{\mathfrak{CM}} 
\newcommand{\CMel}{\mathfrak{CM}_0} 
\newcommand{\JC}{\Lambda} 
\newcommand{\JCfr}{\boldsymbol{\Lambda}} 

\newcommand{\lam}{\lambda}
\newcommand{\blam}{\boldsymbol{\lambda}}
\newcommand{\rev}{^{\mathrm{rev}}}

\newcommand{\scc}{simple closed curve\xspace}
\newcommand{\sccs}{simple closed curves\xspace}
\newcommand{\g}{\gamma}

\newcommand{\SD}{\mathfrak{S}}
	
	\let\S\Sigma
	\let\parasign\temp

\newcommand{\MCG}{\mathrm{Mod}}
\renewcommand{\Crit}{\mathcal{C}}

\renewcommand{\H}{\mathcal{H}}
\newcommand{\PT}{\varPi}

\newcommand{\PP}{\mathbb{P}}

%% file: morecommands_KH.tex
\def\Int{\operatorname{Int}}

\def\Re{\operatorname{Re}}
\def\GL{\operatorname{GL}}
\def\SO{\operatorname{SO}}
\def\O{\operatorname{O}}

\newcommand{\mb}[1]{\mathbb{#1}}
\newcommand{\Bf}[1]{\mathbf{#1}}
\newcommand{\Pa}{\partial}



%% file: intro.tex
\section{Introduction}
	\label{ch:introduction}
The study of maps from 4--manifolds to surfaces has received considerable attention in the last years.
Motivated by Perutz's generalization~\cites{Perutz1,Perutz2} of the Donaldson-Smith invariant \cite{DS} for Lefschetz fibrations to an invariant of broken Lefschetz fibrations introduced in \cite{ADK}, 
Lekili~\cite{Lekili} discussed homotopies between so called wrinkled fibrations, which are stable maps on a $4$--manifold with only indefinite folds and cusps.
Lekili's results were later improved by Williams~\cite{Williams1} and Gay-Kirby~\cites{GK_PNAS,GK_Morse2}. 
The results guarantee that for any two homotopic wrinkled fibrations there exists a homotopy such that all but finitely many maps in the homotopy are wrinkled fibrations and near the exceptional maps the homotopy has one of essentially six possible local models.
While it is easy to see that basic homotopies change the critical value sets as in \cref{fig_discriminant_homotopy}, for a map with one of the critical value configurations shown in \cref{fig_discriminant_homotopy} it is not always possible to realize the corresponding modification of critical values by a homotopy.
In order to guarantee the existence of such a homotopy one has to investigate what we call \emph{vanishing sets} in fibers of the map, which are reminiscent of ascending and descending manifolds in Morse theory (this problem is addressed in~\cite{Williams_crossings}, for example). 
Moreover, even if a homotopy can be found, it is not enough to understand the critical values but one should also keep track of how the vanishing sets behave throughout the homotopy.
On a related note, various existence results for maps with special properties have been obtained.
In particular, Williams's \emph{simple wrinkled fibrations}~\cite{Williams1} can be used to obtain combinatorial descriptions of closed 4--manifolds known as \emph{surface diagrams}.
These diagrams consist of collections of simple closed curves in surfaces arising as vanishing sets of \swfs. 
(Other combinatorial descriptions of 4--manifolds related to the idea of vanishing sets can be obtained from Baykur's \emph{simplified broken Lefschetz fibrations}~\cite{Baykur2} and Gay and Kirby's \emph{trisections}~\cite{GK_trisections}.)
Williams also discussed the uniqueness of simple wrinkled fibrations up to homotopy, introducing four basic homotopies which he called \emph{handleslides}, \emph{stabilizations}, \emph{multislides} and \emph{shifts}~\cite{Williams2}. 
Our main motivation was to understand how these homotopies affect the associated surface diagrams.
However, our results are applicable in more general contexts.
Given a \wf, the general strategy is to fix one reference fiber over each connected component of the regular value set and to collect the vanishing sets associated to neighboring critical values in those fibers.
If one wants to understand how the vanishing sets change in homotopies, the crucial problem is that two formerly disconnected regions of regular values can be joined.
As a result, after the homotopy there is one region that contains two previously unrelated fibers with their vanishing sets.
However, the two fibers can now be related by parallel transport so that the two collections of vanishing sets appear in a single fiber again.
But it turns out that the way they appear relative to each other usually depends on the homotopy.
This phenomenon can be caused by two kinds of homotopies known as \emph{$R_2$--moves} and \emph{cusp merges}. 
The former have been studied in~\cite{HayanoR2} and the latter are the main focus of the present paper.

It is interesting to note that cusp merge homotopies have already been studied 50 years ago in Levine's work \cite{Levine} on the elimination of cusps. 
It is well known that if the source manifold of a stable map to a surface is closed, then the number of cusps in this map is finite and has the same parity as the Euler characteristic of the source. 
Using cusp merges Levine showed that the Euler characteristic is the only obstruction for the elimination of cusps:~he proved that any stable map to a surface is homotopic to a map with at most one cusp. 
In some sense, part of the present paper is a natural extension of~\cite{Levine}, whence the title.
For the elimination of cusps it was enough to understand \emph{when} a pair of cusps can be eliminated. 
However, as explained in the previous paragraph, the recent developments of the topology of $4$--manifolds have led us to study \emph{how many ways} there are to eliminate a given pair of cusps in order to understand the behavior of vanishing sets.
An analogous situation occurred in Morse theory: in order to prove the h--~and~s--cobordism theorems it was enough to understand \emph{when} two critical points can be canceled, but Cerf's approach to the pseudo-isotopy problem made it necessary to study -- among many other things -- the ambiguities in the cancellation procedure (see~\cites{Cerf,HatcherWagoner}).

\subsection*{Statement of results and outline}
We first address the classification problem of cusp merges for a stable map on a $4$--manifold which cancel a given pair of indefinite cusps. 
For this purpose we will define a notion of \emph{elementary cusp merges} and characterize them up to homotopy by their so called \emph{framed joining curves}, which are certain (suitably decorated) curves that connect the cusps.
We will prove the following result in \cref{ch:merge homotopies} to which we also refer for  precise definitions.
The proof relies on some technical results that we outsource into \cref{ch:normalization,ch:symmetries}.
\begin{restatable}{maintheorem}{restateCMvsJC}\label{T:main theorem:cusp merges and joining curves}
Let~$f\colon X\ra B$ be a stable map and~$p,q\in X$ a pair of indefinite cusps.
Denote by $\CM_0(f;p,q)$ and $\JCfr(f;p,q)$ the spaces of elementary cusp merge homotopies and framed joining curves for~$p$ and~$q$.
Then there are canonical bijections between $\pi_0\big(\CM_0(f;p,q)\big)$ and $\pi_0\big(\JCfr(f;p,q)\big)$.
\end{restatable}
%

\begin{figure}[t!]
\centering
\subfigure[]{\includegraphics[height=16mm]{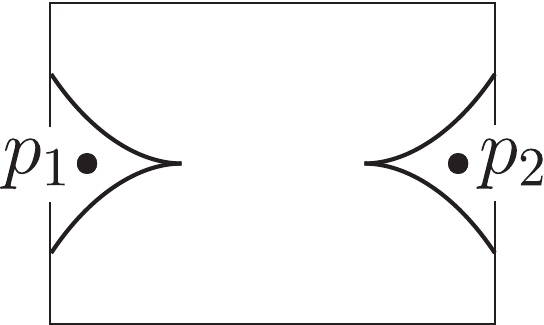}
\label{F:discriminant cusp merges1}}
\hspace{1em}
\subfigure[]{\includegraphics[height=16mm]{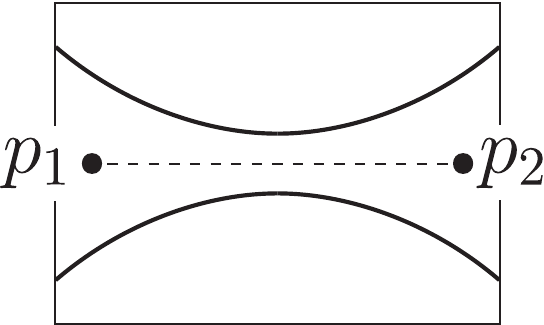}
\label{F:discriminant cusp merges2}}
\caption{Critical value sets in a cusp merge.}
\label{F:discriminant_cusp merges}
\end{figure}

%
Building on \cref{T:main theorem:cusp merges and joining curves} we go on to clarify how cusp merges between indefinite cusps affect vanishing sets. 
As explained, the problem comes from the fact that the two regions in the left and the right sides of \cref{F:discriminant cusp merges1} are joined as shown in \cref{F:discriminant cusp merges2} and the two fibers~$\S_1$ and~$\S_2$ over~$p_1$ and~$p_2$ can be identified by parallel transport along the dotted path in \cref{F:discriminant cusp merges2}.
As a consequence, the two collections of vanishing sets in~$\S_1$ and~$\S_2$ derived from the initial map can be considered in only one of these fibers after the homotopy.
We can thus understand the behavior of vanishing sets in cusp merges once we know which diffeomorphisms can appear as parallel transport. 
For the statement of our result, recall that for a non-separating \scc~$c\subset\S$ in a closed, orientable surface there is a so called \emph{surgery homomorphism}~$\Phi_c\colon\MCG(\S)(c)\ra\MCG(\S'_c)$ where~$\MCG(\S)(c)$ consists of all mapping classes that fix the isotopy class of~$c$ and~$\S'_c$ is the surface obtained by surgery on~$c$ (see~\cref{ch:mapping class groups} for more information).

\begin{restatable}{maintheorem}{restateCMvsMCG}\label{T:main theorem:cusp merges and MCG}
The subset of $\pi_0(\Diff(\S_1,\S_2))$ obtained from merging cusps along a fixed arc in~$B$ has a free and transitive action of either of the groups
	\begin{equation*}
	\mc K(c_i,d_i)=
	\ker\Phi_{c_i}\cap\ker\Phi_{d_i}\subset\MCG(\S_i),\quad i=1,2
	\end{equation*}
where~$c_i,d_i\subset\S_i$ are \sccs related to the vanishing sets of the cusps.
\end{restatable}
\noindent
We will also describe a set of generators for~$\mc K(c_i,d_i)$ in \cref{T:generating the image of Psi}.

\smallskip

In \cref{ch:applications} we discuss some applications of \cref{T:main theorem:cusp merges and MCG} to the theory of surface diagrams.
\cref{ch:surface diagram moves} is concerned with the uniqueness problem for surface diagrams.
More precisely, we use \cref{T:main theorem:cusp merges and MCG} to study how surface diagrams change in two types of homotopies introduced in~\cite{Williams2}, the so called \emph{multislides} and \emph{shifts}.
The analogous problem for the remaining homotopies (handleslides and stabilizations) was studied by the second author in~\cite{HayanoR2}.
Finally, in \cref{ch:Lefschetz and Heegaard} we discuss how to obtain surface diagrams from some constructions of \swfs involving cusp merges.
More precisely, we show how to construct surface diagrams for total spaces of Lefschetz fibrations from the knowledge of the (Lefschetz) vanishing cycles,
and for products of 3--manifolds with the circle from a given Heegaard diagram.
As concrete examples, we obtain new surface diagrams for~$S^4$ and~$S^1\times S^3$ (see \cref{fig_SD_S4,surface_diagram_sphere}).
These have the interesting property that they are \emph{not} related to the previously known diagrams given in~\cite{HayanoR2} by the moves discussed in \cref{ch:surface diagram moves} because the corresponding \swfs are not homotopic.

\smallskip

The observant reader will have noticed that we have neglected to mention \cref{ch:preliminaries,ch:vanishing sets chapter} so far.
These contain necessary background material for \cref{ch:merge homotopies,ch:mapping class group interpretation} that might also be of independent interest.
Although the main purpose of \cref{ch:preliminaries} is to introduce terminology and notation used in the subsequent sections, we go a bit further.
The recent developments in the study of stable maps on $4$--manifolds have heavily relied on folklore facts in singularity theory. 
However, the authors found it extremely difficult to find references including complete proofs of these results.  
For this reason, we take the chance to review these facts with outlines of proofs.
In the remaining \cref{ch:vanishing sets chapter} we introduce a notion of \emph{connections} for general smooth maps generalizing the usual concept for fiber bundles. 
We thereby provide a conceptual framework for discussing vanishing sets and parallel transport which has previously been done in an ad hoc fashion (if at all).

%% file: preliminaries.tex
\section{Generalities about Maps to Surfaces}
	\label{ch:preliminaries}
The purpose of this section is twofold. 
First and foremost, we will introduce necessary terminology and notation that will be used throughout and explain the necessary background for our results.
Since the proofs of some technical results require more sophisticated notions from singularity theory (such as stability of map-germs and their unfoldings), our review will be a little more extensive than one might expect.
Second, we take the chance to give a survey some ``generic'' properties of maps from 4--manifolds to surfaces and homotopies between them from the perspective of singularity theory.
Specifically, we want to address the following ``well known facts'' that have been used in a number of papers in the context of wrinkled fibrations and Morse 2--functions.
\begin{enumerate}[(1)]\label{pg:rough statements}
	\item 
		``Generic'' maps from 4--manifolds to surfaces have only folds and cusp singularities and if a critical value is covered by more than one critical point, then it is covered by two folds whose images meet transversely in the target.
	\item
		In ``generic'' homotopies of maps from 4--manifolds to surfaces all but finitely many maps have the above structure.
		While passing through the exceptional maps one of six phenomena occurs in which the critical image changes as shown in \cref{fig_discriminant_homotopy}.
\end{enumerate}
%
%
While the first statement indeed follows from rather standard (albeit not completely trivial) results in singularity theory, the second one is a little more complicated.
In fact, not much literature is available about global phenomena in families of maps.
What both statements have in common is that it is surprisingly difficult to find concrete references in the literature.
We are aware of only one source, namely~\cite{Chincaro}, but unfortunately it is very hard to find.
So as a third purpose, we aim to provide a more readily available reference.
To be clear, we do not claim any originality.
\cref{T:maps to surfaces} can be proved immediately using Mather's criterion of stability in \cite{MatherV}, while \cref{T:generic homotopy1,T:1-parameter versal unfoldings} are contained in~\cite{Chincaro}, although our proofs are slightly different.
We are aware that what we write will be cryptic to some and at the same time obvious to others.
However, we hope that some low dimensional topologists who are working with surface valued maps will find our exposition helpful, either as a reference or to obtain some idea about the inner workings of the ``black box'' known as singularity theory.

\subsection{Basic notation and terminology}
	\label{ch:singularity theory lingo}
Throughout we assume that all manifolds and maps are smooth.
As a precautionary measure we also assume that all manifolds are connected and oriented, and that sources of maps are closed unless otherwise noted.
We reserve the letters~$X$ and~$B$ for 4--manifolds and surfaces, respectively.
In more general situations we follow the tradition in singularity theory and consider maps between manifolds~$N$ and~$P$ of dimensions~$n$ and~$p$;
in this case we always assume that~$n\geq p$.
%

\subsubsection*{Maps and homotopies.}
Let~$\Cinfty(N,P)$ be the set of smooth maps from~$N$ to~$P$. 
Throughout the paper we endow this set with the \emph{Whitney $\Cinfty$~topology}%
	\footnote{Since we are assuming that $X$~is closed, the Whitney topology agrees with the topology of uniform convergence of all partial derivatives and is generally rather well-behaved.
	However, for non-compact~$X$ the story is more complicated.
	};
for its definition and basic properties we refer to~\cite{GG}*{II.\parasign3}. 
As usual in differential topology, objects of interest are studied up to diffeomorphism.
In this spirit, two maps $f\in\Cinfty(N,P)$ and $f'\in\Cinfty(N',P')$ are called \emph{(right-left) equivalent}\label{D:right-left eq} 
if there exist diffeomorphisms $\phi:N\to N'$ and $\psi:P\to P'$ which satisfy $f'\circ \phi = \psi\circ f$. 
%
%
%
%
%
%
%
By a \emph{homotopy} of maps we mean a smooth map $F\colon J\times N\ra P$ where $J$~is some fixed interval; we will mostly use~$J=[-1,1]$.
We frequently think of homotopies as 1--parameter families of maps~$F=(f_s)_{s\in J}$ where $f_s\colon N\ra P$ is defined by~$f_s(p)=F(s,p)$. 
We will mostly consider homotopies with fixed initial map which we simply denote by the lower case letter~$f$.
Another useful way to think of homotopies is to consider the map~$\wt{F}\colon J\times N \lra J\times P$ defined by $\wt{F}(s,p)=(s, f_s(p))$.
We call $\wt{F}$ the \emph{unfolding} associated to~$F$.
%
Of course, the objects~$F$, $(f_s)$, and~$\wt{F}$ all contain the same information but each perspective has its advantages.
A natural notion of \emph{equivalence} for homotopies $F\in \Cinfty(J\times N,P)$ and $F'\in \Cinfty(J\times N',P')$ is given by commutative diagrams of the form
	\begin{equation}\label{eq:equivalence of homotopies}
	\begin{CD}
	J\times N  @>{\wt{F}}>>  J\times P  @>{\mathrm{pr}_J}>>  J\\
	@VV{\Phi}V   @VV{\Psi}V   @VV{\theta}V\\
	J\times N'  @>{\wt{F}'}>> J\times P'  @>{\mathrm{pr}_J}>>  J
	\end{CD}
	\end{equation}
where $\Phi$, $\Psi$, and $\theta$ are diffeomorphisms and~$\mathrm{pr}_J$ denotes the projections onto~$J$.
We say that $F$~is \emph{constant} if~$\wt{F}=\id_J\times f$ or, equivalently, $f_s=f$ for all~$s\in J$.
We call $F$ \emph{trivial} if it is equivalent to the constant homotopy at the initial map of~$F$.
In this case, we can assume that~$\theta=\id_J$ so that we can consider~$\Phi$ and~$\Psi$ as families of diffeomorphisms~$(\phi_s)\in\Diff(N)$ and~$(\psi_s)\in\Diff(P)$.
Moreover, they satisfy $f_s=\psi_s\circ f \circ \phi_s\inv$ and we can assume that both families of diffeomorphisms start from the identity map;
we will refer to homotopies of this form as \emph{isotopies} of~$f$.%
	\footnote{The reader should be warned that ``isotopy'' has a slightly different meaning in~\cite{Lekili} and~\cite{Williams1}, there it refers to a homotopy that stays within a special class of maps.  Such homotopies are usually not related to isotopies of~$X$ and~$B$.}

\subsubsection*{Critical points.}
Recall that a \emph{critical point} of a smooth map $f\colon N\ra P$ is a point~$x\in N$ where the derivative~$df_x$ fails to have maximal rank.
The number~$p-\rk(df_x)$ is called the \emph{corank} of~$f$ at~$x$.
The image of a critical point is called a \emph{critical value}.
We denote the sets of critical points and critical values by $\mc C(f)\subset N $ and $\mc D(f)=f(\mc C(f))\subset P $, respectively, 
and refer to~$\mc C(f)$ as the 
	\emph{critical locus} 
and to~$\mc D(f)$ as the \emph{critical image} or \emph{discriminant}.
Points in the complements of~$\mc C(f)$ and~$\mc D(f)$ are called \emph{regular} points and values, respectively.

\subsubsection*{Germs, singularities, and local models.}
In order to study local structure of maps, it is convenient to use the language of germs.
Let $f\colon N\ra P$ be a smooth map and let~$S\subset N$ be a subset.
Recall that the \emph{germ} of~$f$ at~$S$ is the equivalence class of maps that are defined in a \nbhd of~$S$ and agree with~$f$ on a possibly smaller \nbhd.
The most important special cases are when~$S$ is a single point or a finite set of~$N$ such that~$f(S)$ is a point in~$P$.
In these situations we speak of \emph{mono-germs} and \emph{multi-germs}.
It is customary to denote the germ of~$f$ at~$S$ by $f\colon (N,S)\ra (P,f(S))$ which should not be confused with the notation for maps of pairs.
Two germs $f:(N,S)\to(P,f(S))$ and $f':(N',S')\to(P',f(S'))$ are called \emph{equivalent} if there are germs of diffeomorphisms $\phi:(N,S)\to (N',S')$ and $\psi:(P,f(S))\to(P',f(S'))$ which satisfy $f'\circ \phi = \psi\circ f$. 
%
A mono-germ $f\colon (N,x)\ra (P,y)$ is \emph{singular} if~$x$ is a critical point of some (hence any) representative.
With this understood, a \emph{singularity} is an equivalence class of singular mono-germs.
More generally, a \emph{multi-singularity} is an equivalence class of multi-germs such that each of its mono-germs is singular.
Obviously, by suitable choices of coordinates all (multi-)singularities become equivalent to (finite collections of) singular mono-germs $(\R^n,0)\ra(\R^p,0)$.
The set of germs $(\R^n,0)\ra(\R^p,0)$ is usually denoted by~$\mc E(n,p)$.
If a map $\mu\colon \R^n\ra\R^p$ with~$\mu(0)=0$ represents a given singularity, we will refer to~$\mu$ as a \emph{local model} for that singularity.
%
%

\subsection{Maps from 4--manifolds to surfaces}
	\label{ch:stable maps}
The first task is to give a precise version of the statement about the generic structure of maps from a 4--manifold~$X$ to a surface~$B$.
We first define the relevant multi-singularities in terms of local models 
	(which are maps from~$\R^4$ to~$\R^2$).
Here, and throughout the rest of the paper, we will denote coordinates on the source by~$(t,x,y,z)$ and on the target by~$(u,v)$.
The first pair of models describes the \emph{fold singularities}
	\begin{equation}\label{eq:fold models}
	\big(t,\, x,y,z\big) \mapsto \big( \; t,\; x^2 +y^2\pm z^2\big)
	\end{equation}
which can be thought of as a trivial family of 3--dimensional Morse singularities of index~$1$~or~$0$.
In the index~1 case, that is, when the sign is negative, we speak of \emph{indefinite} folds while the index~$0$ folds are called \emph{definite}.
If we superimpose two fold singularities in such a way that their discriminants intersect transversely, 
we obtain a multi-singularity which we call a \emph{transverse double fold}.
Next there are the (\emph{definite} and \emph{indefinite}) \emph{cusp singularities}
	\begin{equation}\label{eq:cusp models}
	\big(t,\, x,y,z\big) \mapsto \big( \; t,\; x^3+3tx +y^2\pm z^2\big).
	\end{equation}
With these definitions (together with the topology of $\Cinfty(X,B)$) in hand, we can make the statement (1) in \cpageref{pg:rough statements} precise as follows:
\begin{theorem}\label{T:maps to surfaces}
Let~$X$ be a closed 4--manifold and~$B$ a surface. 
Denote by~$\mc S^0(X,B)$ the set of smooth maps~$X\ra B$ whose only multi-singularities are folds, cusps and transverse double folds.
Then~$\mc S^0(X,B)$ is open and dense in~$\Cinfty(X,B)$.
\end{theorem}
As we will explain, $\mc S^0(X,B)$ agrees with the set of \emph{stable} maps which are classical objects in singularity theory. 
Since we are deliberately taking the point of view of singularity theory, we will henceforth use this terminology.
We note, however, that these maps have gained some popularity under the name of \emph{Morse \mbox{2--functions}} in the low dimensional topology community through the work of Gay and Kirby~\cites{GK_Morse2,GK_PNAS,GK_trisections}.
We now go on to explain how \cref{T:maps to surfaces} follows from the vastly more general theory of stable maps. 

\subsubsection{A glance at stability theory}
	\label{ch:stability}
In this interlude we first consider smooth maps between manifolds $N$ and~$P$ of arbitrary dimensions which we write as a pair~$(n,p)$.
As before, we assume that~$N$ is closed and equip $\Cinfty(N,P)$ with the Whitney $\Cinfty$~topology.
We begin with the central definition.
\begin{definition}[Stable maps]\label{D:stable maps}
A map is $f\colon N\ra P$ is called \emph{stable} if there is an open \nbhd $\mc U_f$ of~$f$ in $\Cinfty(N,P)$ 
such that each $g\in\mc U_f$ is equivalent to~$f$, that is, there are 
	$(\phi_g,\psi_g)\in\Diff(N)\times\Diff(P)$ 
such that $g=\psi_g\circ f\circ\phi_g\inv$.
\end{definition}
The group~$\Diff(N)\times\Diff(P)$ is commonly denoted by $\mc A=\mc A(N,P)$ in this context.
Stability of smooth maps was carefully studied in a highly influential series of papers by Mather (including~\cites{MatherV,MatherVI}) who attributes the above definition to Thom.
The theory has long matured and excellent textbook references are available, for example \cites{GG,Martinet,Gibson}.
Note that stability is an open condition by definition.
Moreover, a deep theorem of Mather states that stable maps are also dense in~$\Cinfty(N,P)$ provided that the dimensions $(n,p)$ lie in the so called range of ``nice dimensions'' which includes~$(n,2)$ for all~$n$ (see \cite{GG}*{VI.\parasign6} or~\cite{MatherVI}).
So in order to prove \cref{T:maps to surfaces} it is enough to show that the~$\mc S^0(X,B)$, which was defined by restricting the allowed multi-germs, agrees with the set of stable maps in dimensions~$(4,2)$.
This is done in two steps: the first is a characterization of stable maps in terms of their multi-germs, and the second is the classification of the multi-germs that appear in stable maps.
Observe that $\mc A$ acts on $\Cinfty (N,P)$ and $f$~is stable if and only if its $\mc A$--orbit is open.
Intuitively, one should think of the finite dimensional situation where a Lie group acts properly on a manifold.
In this case, the tangent spaces of orbits provide valuable information about the manifold.
%
%
Even though there are no completely satisfactory manifold structures on~$\Cinfty(N,P)$ and~$\mc A$, there are natural candidates for their tangent spaces.
The main theorem of~\cite{MatherV} states that $f$~is stable if and only if the (formal) tangent space to its $\mc A$--orbit has codimension~0 in the tangent space to~$\Cinfty(N,P)$; 
this is known as \emph{infinitesimal stability}
and can be interpreted as the surjectivity of the differential of the map~$\mc A\ra\Cinfty(N,P)$ given by evaluation of the action at~$f$ at the identity (see~\cite{GG}*{III.\parasign1}).
In general, the codimension of~$T_f\mc Af$ in~$T_f\Cinfty(N,P)$ is called the \emph{codimension of~$f$}.
A similar discussion applies to germs.
A group of diffeomorphism germs, usually also denoted by~$\mc A$, acts on the set of map germs and all relevant objects have (formal) tangent spaces;
for more details see \cite{Gibson} or~\cite{Wall}, for example.
This allows to define the \emph{codimension} of a multi-germ $f\colon (N,S)\ra(P,y)$ as the codimension of the tangent space to its orbit in its full tangent space, and the codimension~$0$ multi-germs are called \emph{(infinitesimally) stable}.%
	\footnote{In fact, there are two notions of codimension related to the $\mc A$--action on germs. We are interested in the $\mc A_e$--codimension which is denoted by $d_e(f,\mc A)$ in \cite{Wall}.}
The connection to the notions of codimension for maps and multi-germs is given by a formula of Mata-Lorenzo~\cite{Mata} which expresses the codimension of a map as the sum of codimensions of all its multi-germs.
In particular, a map is stable if and only if all its multi-germs are stable;
the latter condition is also known as \emph{local stability}.%
	\footnote{The equivalence of stability and local stability was already proved in~\cite{MatherV}*{p.314}.}
We end this digression with a word of warning.
One should be aware that the compactness assumption on~$N$ is crucial for the whole discussion.
In fact, the relations of stability with its infinitesimal and local versions become delicate issues when the source is not compact.
For state of the art accounts on these matters we refer the interested reader to~\cite{duPlessisWall} and~\cite{duPlessis_Vosegaard}.

\subsubsection{The proof of \cref{T:maps to surfaces} (Sketch)}
Now let us go back to the situation of \cref{T:maps to surfaces} where $(n,p)=(4,2)$.
The theory outlined in the previous section reduces the problem to showing that the stable multi-singularities in dimensions~$(4,2)$ are exactly the folds, cusps and transverse double folds.
As a first step, one can show directly using the definition that these germs are infinitesimally stable, which is an instructive exercise.
In the other direction, another exercise in the definitions shows that all stable (mono-)singularities in dimensions~$(4,2)$ have corank~1.
Then the normal forms for corank~1 singularities obtained by Morin~\cite{Morin} show that we are only dealing with folds and cusps.
As a last step, we have to discuss multi-singularities.
It is a basic fact that for any stable multi-singularity $f\colon (N,S)\ra (P,y)$ we have~$|S|\leq p$.
In particular, for~$p=2$ any critical value is covered by at most two critical points.
We leave it to the interested reader to show that 
	(a)~a bi-singularity that involves a cusp cannot be stable, and
	(b) that the discriminants of a stable bi-singularity consisting of two folds must intersect transversely.
Although some of the verifications that we have left out are tedious, they are all elementary. 

\subsection{Homotopies of maps from 4--manifolds to surfaces}
	\label{ch:stable homotopies}
The purpose of this section is briefly reviewing generic properties of homotopies from $4$--manifolds to surfaces, in particular making statement~(2) on \cpageref{pg:rough statements} precise.  
As in the previous section these properties are described in terms of local models. 
We thus first introduce another equivalence relation of germs appropriate for homotopies. 
Let $S\subset \R^n$ be a finite set and $f:(\mb{R}^n,S)\to (\mb{R}^p,f(S))$ a germ. 
A germ $F:(\mb{R}^{n+k},\{0\}\times S)\to \allowbreak (\mb{R}^{p+k},\{0\}\times f(S))$ is called an ($k$--parameter) \emph{unfolding} of $f$ if the restriction $F|_{\{0\}\times \mb{R}^n}$ is equal to $f$ and $pr_{\mb{R}^k} \circ F = pr_{\mb{R}^k}$, where $pr_{\mb{R}^k}:(\mb{R}^{m+k},\{0\}\times T) \to (\mb{R}^k,0)$ is the germ of the projection onto the former components. 
Two $k$--parameter unfoldings $F_0$ and $F_1$ of $f$ are said to be \emph{$k$--equivalent} if there exist germs of diffeomorphisms $\Phi:(\mb{R}^{n+k},\{0\}\times S)\to \allowbreak (\mb{R}^{n+k},\{0\}\times S)$, $\varphi:(\mb{R}^{p+k},\{0\}\times f(S)) \allowbreak \to (\mb{R}^{p+k},\{0\}\times f(S) )$ and $\phi:(\mb{R}^k,0)\to (\mb{R}^k,0)$ such that the following diagram commutes: 
\begin{equation}\label{E:diagram action G}
\begin{CD}
(\mb{R}^{n+k},\{0\}\times S) @> F_0 >> (\mb{R}^{p+k},\{0\}\times f(S)) @> pr_{\mb{R}^k} >> (\mb{R}^k,0) \\
@V \Phi VV  @V \varphi VV  @V \phi VV \\
(\mb{R}^{n+k},\{0\}\times S) @> F_1 >> (\mb{R}^{p+k},\{0\}\times f(S)) @> pr_{\mb{R}^k} >> (\mb{R}^k,0). 
\end{CD}
\end{equation}
An unfolding $F: (\mb{R}^{n+k},\{0\}\times S)\to (\mb{R}^{p+k},\{0\}\times f(S))$ of $f$ is said to be \emph{trivial} if $F$ is $k$--equivalent to the constant unfolding $\id_{\mb{R}^k}\times f$. 

Let $F:(\mb{R}^{n+k},\{0\}\times S)\to (\mb{R}^{p+k},\{0\}\times f(S))$ be an unfolding of a germ $f:(\mb{R}^n,S)\to (\mb{R}^p,f(S))$ and $h:(\mb{R}^l,0)\to (\mb{R}^k,0)$ a germ. 
We can obtain an \mbox{$l$--parameter} unfolding $h^\ast F$ of $f$ as follows: 
\[
h^\ast F (s,x) = (s, pr_{\mb{R}^p}\circ \hat{F}(\hat{h}(s),x)), 
\]
where $\hat{F}$ and $\hat{h}$ are representatives of $F$ and $h$, respectively. 
This unfolding is called the \emph{pull-back} of $F$ by $h$. 
An unfolding $F:(\mb{R}^{n+k},\{0\}\times S)\to (\mb{R}^{p+k},\{0\}\times f(S))$ of $f:(\mb{R}^n,S)\to (\mb{R}^p,f(S))$ is said to be \emph{versal} if any $l$--parameter unfolding of $f$ is $l$--equivalent to some pull-back of $F$ for any $l$. 
The connection to the theory outlined in \cref{ch:stability} is given by the following standard result.
\begin{proposition}[\cite{Martinet}*{p.189}]\label{T:versality and stability}
Let  $f\colon(\R^n,S)\ra(\R^p,0)$ be a multi-germ.
\begin{enumerate}[(i)]
	\item 
		$f$ is stable if and only if it is versal when considered as an unfolding of itself.
		In particular, any unfolding of a stable multi-germ is trivial.
	\item
		More generally, the codimension of~$f$ agrees with the minimal number of parameters needed to obtain a versal unfolding of~$f$.
\end{enumerate}
\end{proposition}
\begin{remark}\label{R:parameterized normal forms_codim0}
As an immediate consequence of \cref{T:versality and stability}~(i) is that whenever we have a normal form for a given stable multi-germ, we get a parametric normal form for free (that is, if the singularity appears embedded in a family of maps, then the whole family is locally equivalent to the constant family of the normal form).
	Moreover, \cref{T:versality and stability}~(i) globalizes to the statement that any family of stable maps is trivial (see~\cite{GG}*{V.\parasign2}).
\end{remark}
For the purpose of understanding homotopies of maps, the significance of the concept of versal unfoldings is illustrated by the following result: 
\begin{theorem}[{\cite{Chincaro}}]\label{T:generic homotopy1}
Let $X$ be a $4$--manifold, $B$ a surface and $J$ some interval. 
Let $\mc H^0(X,B)$ be the set of homotopies $F\colon J\times X\ra B$ with the following properties:
\begin{enumerate}[(a)]
	\item 
		For any $t\in J$ and finite subset $S\subset X$ the germ of $F$ at $\{t\}\times S$ is a versal unfolding of the germ of~$f_t$ at~$S$.
	\item
		Each map $f_t$ has codimension at most 1 (that is, it contains at most one multi-germ of codimension 1).
\end{enumerate}
Then $\mc H^0(X,B)$ is dense in $\Cinfty(J\times X,B)$ \wrt the Whitney topology.
\end{theorem}
In analogy with the situation of maps, we call the elements of $\mc H^0(X,B)$  \emph{stable homotopies}. 
In fact, this terminology is justified by the work of Chincaro~\cite{Chincaro} who develops a theory of stability for families of maps. 
As in the situation of maps, this immediately implies that~$\mc H^0(X,B)$ is also open in~$\Cinfty(J\times X, B)$.
Before giving an outline of a proof of \cref{T:generic homotopy1}, we quickly review the classification of versal $1$--parameter unfoldings of germs from $\mb{R}^4$ to $\mb{R}^2$. 
We begin with general remarks on versal unfoldings. 
First, it is easy to see that $\tilde{F}:(\mb{R}^{n+k},\{0\}\times S)\to \allowbreak (\mb{R}^{p+k},\{0\}\times f(S))$ is versal if and only if the restriction $\tilde{F}:(\mb{R}^{n+k},\{0\}\times (S\cap f^{-1}(q)))\allowbreak\to (\mb{R}^{p+k},(0,q))$ is versal for any $q\in f(S)$. 
Thus we can assume that $f(S)$ consists of a single point, say $0\in \mb{R}^p$, without loss of generality. 
Second, a multi-germ $f:(\mb{R}^n,S)\to (\mb{R}^p,0)$ has a versal \mbox{$1$--parameter} unfolding if and only if it has codimension at most one (in the sense of \cref{ch:stability}, see \cref{T:versality and stability}).
Furthermore, in this case the versal unfolding of~$f$ is unique up to equivalence.
In particular, we can classify $1$--parameter versal unfoldings once we obtain the classification of codimension--$1$ germs (and these versal unfoldings). 
Lastly, since the codimension of $f:(\mb{R}^n,S)\to (\mb{R}^p,0)$ is equal to that of $f\sqcup s:(\mb{R}^n\sqcup \mb{R}^n,S\sqcup \{0\})\to (\mb{R}^p,0)$, where $s:(\mb{R}^n,0)\to (\mb{R}^p,0)$ is a submersion-germ, in what follows we only deal with singular germs.  

\begin{figure}[t!]
\centering
\subfigure[birth/death]{\includegraphics[height=12mm]{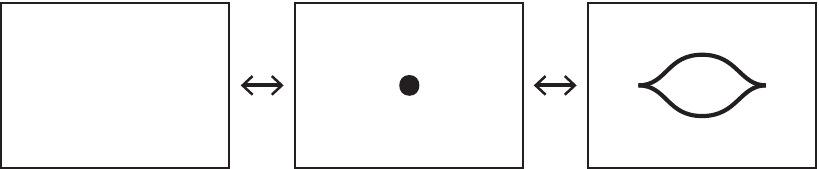}}
\hspace{1em}
\subfigure[fold/cusp merge]{\includegraphics[height=12mm]{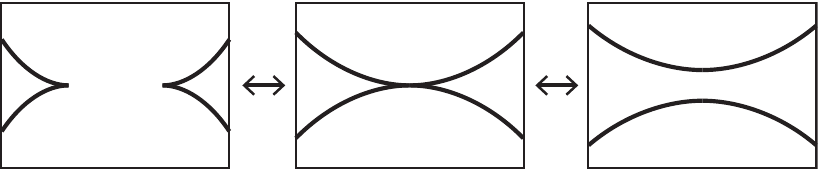}}

\subfigure[flip]{\includegraphics[height=12mm]{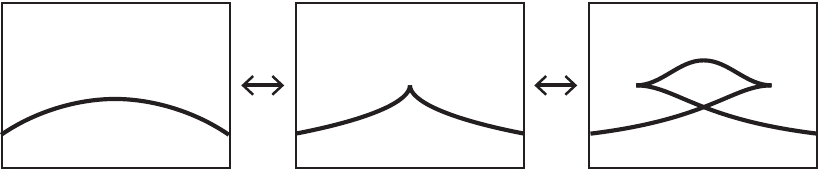}}
\hspace{1em}
\subfigure[cusp-fold crossing]{\includegraphics[height=12mm]{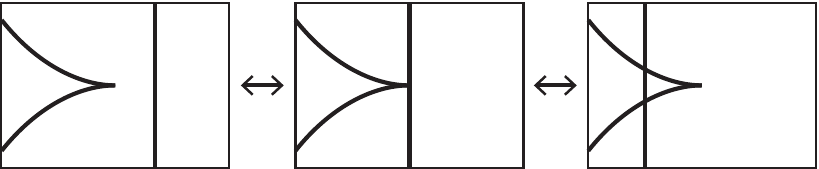}}

\subfigure[$R_2$--move]{\includegraphics[height=12mm]{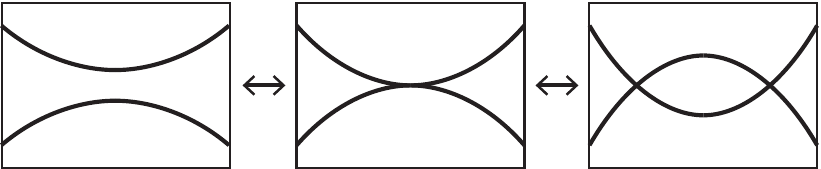}}
\hspace{1em}
\subfigure[$R_3$--move]{\includegraphics[height=12mm]{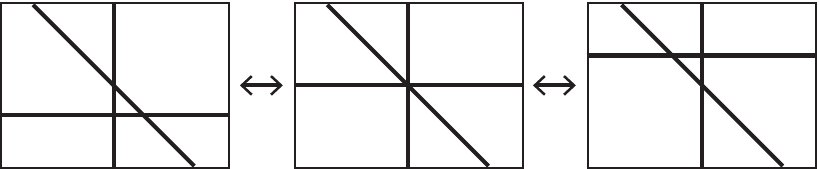}}
\caption{The discriminants of versal unfoldings.}
\label{fig_discriminant_homotopy}
\end{figure}

We now turn our attention to the specific dimension pair $(n,p)=(4,2)$. 
The germs with codimension $0$ are nothing but stable germs given in~\cref{ch:stable maps}.
As for germs with codimension $1$, one can obtain a list of mono-germs using \cite{Rieger_Ruas}*{Lemma~1.1} together with the classification of mono-germs between planes with small $\mathcal{A}$--codimensions due to Rieger~\cite{Rieger}. 
The codimension--$1$ multi-germs can be determined using an algorithm due to Cooper, Mond and Wik Atique~\cite{Cooper_Mond_WikAtique}*{Theorem~5.22 and Remark~5.23}.
We eventually obtain the following list of $1$--parameter versal unfoldings of germs from $\mb{R}^4$ to $\mb{R}^2$: 

\begin{theorem}\label{T:1-parameter versal unfoldings}

Any versal $1$--parameter unfolding of a germ $f:(\mb{R}^4, S)\to (\mb{R}^2,0)$ is $1$--equivalent to either a trivial unfolding of a stable germ or one of the six types of unfoldings of germs with codimension $1$, whose discriminants are shown in \cref{fig_discriminant_homotopy}.

%
%

\end{theorem}

\noindent
The former three germs in \cref{fig_discriminant_homotopy}, named birth/death, fold/cusp merge and flip, are mono-germs, while the other ones are multi-germs. 
We can find local models of the three mono-germs in \cite{Lekili}, for example.
As for the multi-germs, one can easily obtain local models of them following an algorithm in \cite{Cooper_Mond_WikAtique}. 

%
%

\begin{proof}[Proof of \cref{T:generic homotopy1} (Sketch)]
As usual in differential topology, density is proved using some form of transversality.
In this case, (a)~can be rephrased as a transversality condition on some multi-jet extension of~$F$ relative to the parameter and the density follows from standard methods (see~\cite{Wall_2009}*{2.1,~2.2}\footnote{Wall discusses $\mathcal{K}_e$--versal unfoldings, but a similar statement holds for $\mathcal{A}_e$--versal unfoldings.} or~\cite{Chincaro}*{III.4.}).
Similarly, one can use transversality to show that the maps~$f_s$ generically avoid all multi-singularities of codimension at least two which proves the density of~$(b)$.
Roughly, one has to show that the unions of orbits of multi-singularities of higher codimension have codimension at least two.
More precisely, this boils down to estimating the codimension of algebraic sets~${}_s\Sigma \subset {}_sJ^r(4,2)$ for~$s\leq 3$ and sufficiently large~$r$ as in~\cite{MatherV},
where ${}_s\S$ is the union of $\mathcal{A}$--orbits of $s$--fold $r$-jets of multi-singularities with codimension at least two.
This calculation is not hard but tedious, thus for reasons of brevity we leave the details to the really interested reader.
\end{proof}

\subsubsection{Stability of unfoldings}
	\label{ch:stable unfoldings}
Since we need to deal with families of homotopies in~\cref{ch:merge homotopies}, we introduce a notion of unfoldings for unfoldings of germs and equivalence relation for them, which is a generalization of \emph{$(r,s)$--equivalence} for unfoldings of function germs introduced by Wassermann \cite{Wassermann}. 

Let $F:(\mb{R}^{n+k},0)\to \allowbreak (\mb{R}^{p+k},0)$ be a $k$--parameter unfolding of a germ $f$. 
An unfolding $\Bf{F}:(\mb{R}^{n+k+l},0)\to \allowbreak (\mb{R}^{p+k+l}, 0)$ of $f$ is called an ($l$--parameter) \emph{unfolding} of $F$ if the restriction $\Bf{F}|_{\{0\}\times \mb{R}^{n+k}}$ is equal to $F$. 
Two $l$--parameter unfoldings $\Bf{F}_0, \Bf{F}_1:(\mb{R}^{n+k+l},0)\to \allowbreak (\mb{R}^{p+k+l}, 0)$ of a $k$--parameter unfolding $F$ of $f$ are said to be \emph{$(k,l)$--equivalent} if there exist germs of diffeomorphisms $\Bf{\Phi}$, $\Phi$, $\varphi$ and $\phi$ on $(\mb{R}^{n+k+l},0)$, $(\mb{R}^{p+k+l},0)$, $(\mb{R}^{k+l},0)$ and $(\mb{R}^{l},0)$, respectively, such that the following diagram commutes: 
\begin{equation}\label{E:diagram action G_e}
\begin{CD}
(\mb{R}^{n+k+l},0) @> \Bf{F}_0 >> (\mb{R}^{p+k+l},0) @> pr_{\mb{R}^{k+l}}>> (\mb{R}^{k+l},0) @> pr_{\mb{R}^l}>> (\mb{R}^{l},0) \\
@V \Bf{\Phi}VV @V \Phi VV @V \varphi VV @V\phi VV \\
(\mb{R}^{n+k+l},0) @> \Bf{F}_1 >> (\mb{R}^{p+k+l},0) @> pr_{\mb{R}^{k+l}}>> (\mb{R}^{k+l},0) @> pr_{\mb{R}^l}>> (\mb{R}^{l},0).
\end{CD}
\end{equation}
An unfolding $\Bf{F}$ of an unfolding $F$ is said to be \emph{trivial} if $\Bf{F}$ is $(k,l)$--equivalent to the constant unfolding $\id\times F$. 
In this paper we need the following theorem: 

\begin{theorem}[\cite{Martinet}]\label{T:stability versal unfoldings}
If $F$ is a versal unfolding of a germ~$f$, then every unfolding of~$F$ is trivial as an unfolding of~$F$ (and thus also as an unfolding of~$f$).
\end{theorem}

\begin{remark}

\cref{T:stability versal unfoldings} \emph{does not} follow directly from uniqueness of versal unfoldings (i.e.~\cite{Martinet}*{p.190, Theorem 1.2}). 
Indeed, the uniqueness of versal unfoldings only guarantees that any $l$--parameter unfolding $\Bf{F}$ of a $k$--parameter versal unfolding $F$ of $f$ is $(k+l)$--equivalent to the constant unfolding of $F$. 
However, the construction of diffeomorphisms in \cite{Martinet}*{Ch.~XIV} gives rise to $(k,l)$--equivalence between $\Bf{F}$ and the constant unfolding of $F$. 

\end{remark}

\begin{remark}\label{R:parameterized normal forms_codim1}
\cref{T:stability versal unfoldings} provides an extension of the parameterized normal forms for stable germs in families mentioned in~\cref{R:parameterized normal forms_codim0}.
Whenever a versal unfolding which has a normal form is embedded in a higher dimensional family of maps, then there is a normal form for the whole family given as a trivial product with the original normal form.
\end{remark}

\subsection{Wrinkled fibrations}
	\label{ch:definition wrinkled fibrations}
A map $f\in \mc S^0 (X,B)$ is called a \emph{wrinkled fibration}%
	\footnote{\Wfs are also known as 
	\emph{fiber-connected, indefinite Morse $2$--functions}~\cite{GK_Morse2}.} 
if~$f$ does not have definite folds and all fibers are connected. 
Note that a wrinkled fibration cannot have any definite cusps since these require definite folds.
It is easy to see that wrinkled fibrations are open maps, in particular they are surjective. %
In the case~$B=S^2$ a wrinkled fibration $f\colon X \to S^2$ over $S^2$ is said to be \emph{simple} if 
	the critical set $\Crit(f)$ is connected, 
	the restriction $f|_{\Crit(f)}$ is injective, and 
	$f$~has at least one cusp. 

As for the existence of \wfs, Saeki~\cite{Saeki} first proved that any smooth map $f:X\to S^2$ is homotopic to a \wf, and later Gay and Kirby~\cite{GK_Morse2} generalized this result:~they showed that $f:X\to B$ is homotopic to a \wf if and only if $f_\ast\pi_1(X,x)$ has finite index in $\pi_1(B,f(x))$. 
They further proved that any two homotopic \wfs can be connected by a stable homotopy $f_t$ such that each fiber of $f_t$ is connected for all~$t$ and $f_t$ is a \wf for all but finitely many values of~$t$. 
The latter statement was also proved by Williams~\cite{Williams1} for~$B=S^2$ who also shows that any map $f:X\to S^2$ is homotopic to a \swf. 
Moreover, using this fact he introduced \emph{surface diagrams} which describe closed $4$--manifolds by sequences of simple closed curves in closed surfaces. 
The simple closed curves in such a diagram represent \emph{vanishing cycles} of indefinite folds in a simple wrinkled fibration, which reflects configuration of singularities in the fibration. 
We will discuss vanishing cycles in detail in \cref{ch:vanishing sets chapter} and surface diagrams will be studied in \cref{ch:applications}.

%% file: merging.tex
\section{Elimination of Cusps I: Merge Homotopies and Joining Curves}
	\label{ch:merge homotopies}
In this section we will prove our first main result which we restate for convenience.
\restateCMvsJC*
We begin by giving precise definitions of all involved objects.
Recall from \cref{ch:singularity theory lingo} that we consider $\Cinfty(J\times X,B)$ where $J=[-1,1]$ as the space of homotopies of maps from~$X$ to~$B$ and that we can think of a homotopy~$F$ as a 1--parameter family of maps~$(f_s\colon X\ra B)_{s\in J}$.%
	\footnote{Since we are assuming that~$X$ is compact, the map~$s\mapsto f_s$ is continuous in the Whitney $\Cinfty$~topology.
	This is not true for non-compact~$X$! In the non-compact situation, continuity is equivalent to the homotopy being constant outside of some compact subset of~$X$.}
We will be interested in following special class of homotopies.
\begin{definition}[Merge homotopies]\label{D:merge homotopies}
A stable homotopy~$F=(f_s)$ is called a \emph{merge homotopy}, or simply a \emph{merge}, if the following two conditions are satisfied:
\begin{enumerate}[(a)]
	\item 
		All but one multi-germs that appear in the maps~$f_s$ are stable except for a single beak-to-beak point in~$f_{s_0}$ for some~$s_0\neq\pm1$.
	\item
		The numbers of cusps of~$f_{-1}$ and~$f_1$ differ by two.
\end{enumerate}
More precisely, we call $F$ a \emph{cusp merge} if the number of cusps decreases, and a \emph{fold merge} otherwise.
We denote by $\CM$ and $\FM$ for the subspaces of $\Cinfty(J\times X,B)$ formed by the cusp and fold merges.
Moreover, we write~$\CM(f)$ and $\CM(f,g)$ for the cusp merges starting with~$f$ or going from~$f$ to~$g$, respectively.
\end{definition}
In what follows we will mostly be concerned with cusp merges and only briefly comment on fold merges in \cref{ch:merging folds}.
Recall from \cref{ch:stable homotopies} for any~$F\in\CM$ the passage of through the beak-to-beak point is governed by the \emph{cusp merge model}
	\begin{equation}\label{eq:cusp merge model}
	\mu_s(t,x,y,z) = \big(t, x^3-3(t^2+s)x+y^2\pm z^2\big), \quad s\in\R.
	\end{equation}
A direct calculation shows that the critical locus of~$\mu_s$ is 
cut out by the equations $x^2-t^2=s$ and $y=z=0$ and mostly consists of folds except for 
	pairs of cusps for~$s<0$ located at~$(\pm\sqrt{-s},0,0,0)$ 
	which approach 
	a beak-to-beak point in the origin for~$s=0$. 
In particular, for $s>0$ all critical points are folds.
Note that for $s\in [-1,0]$, the cusps and the beak-to-beak trace out the line segment
	\begin{equation}\label{eq:model joining curve}
	L_0=\Set{(\tau,0,0,0)}{\tau\in[-1,1]}\subset\R^4
	\end{equation}
which will soon play an important role.

\subsection{Elementary cusp merge homotopies}
	\label{ch:elementary merge homotopies}
While the cusp merge model gives a precise picture of cusp merge homotopies near their beak-to-beak singularities, we do not have any control over the behavior further away.
And although all other multi-singularities are stable, they might move in complicated patterns.
It is therefore more convenient to work with a special class of cusp merge homotopies whose effect is purely local.
%
%
%
%
\begin{definition}[Elementary cusp merges]
	\label{D:elementary cusp merges}
A cusp merge $F\in\CM(f)$ is called \emph{elementary} if there are open subsets $U\subset V\subset X\times [-1,\varepsilon]$ (for sufficiently small $\varepsilon>0$) such that
\begin{enumerate}[(i)]
	\item
		the image of $V$ under the projection to $X$ is contained in a compact set $K$ diffeomorphic to the $4$--ball and $F$ is constant in $X\setminus K$,	
	\item
		there are coordinates on~$U$ and~$\widetilde{F}(U)$, called \emph{merge coordinates}, in which $\widetilde{F}$ coincides with the unfolding associated with the model map~$\mu_s$ in \eqref{eq:cusp merge model}, 
	\item 
		there are diffeomorphisms $\Phi$ of $X\times [-1,\varepsilon]$ and $\Psi$ of $B\times [-1,\varepsilon]$ with support in~$V$ and~$\widetilde{F}(V)$, respecitvely, both preserving the parameter levels, such that on $\big(X\times[-1,\varepsilon]\big)\setminus U$ we have
			$\widetilde{F} = \Psi\circ (f_{-1}\times \id)\circ \Phi$. 
\end{enumerate}
We write $\CM_0$, $\CM_0(f)$ and~$\CM_0(f,g)$ for the spaces of elementary cusp merges with or without fixed endpoints.
\end{definition}

In order to justify the terminology, we observe that an arbitrary cusp merge can be deformed into an isotopy followed by an elementary cusp merge and another isotopy.
Recall that an isotopy is a homotopy of the form $f_s=\psi_s\circ f\circ \phi_s\inv$ where $(\phi_s,\psi_s)\in\Diff(X)\times\Diff(B)$ with~$\phi_{-1}=\id_X$ and~$\psi_{-1}=\id_B$; 
let $\mf{Iso}(f,g)$ be the space of isotopies from~$f$ to~$g$.
\begin{lemma}\label{T:reduction to elemetary cusp merges}
For any cusp merge~$F\in\CM(f,g)$ there are 
	stable maps~$f',g'$ as well as 
	$\mc I\in \mf{Iso}(f,f')$, 
	$F'\in\CM_0(f',g')$, and
	$\mc J\in \mf{Iso}(g',g)$
such that $F$ can be deformed within~$\CM(f,g)$ to the concatenation $\mc I\ast F'\ast \mc J$.
\end{lemma}
In principle, this result could be proved directly using similar methods as in \cref{ch:normalization}.
However, since this would be rather lengthy and technical and \cref{T:reduction to elemetary cusp merges} only serves as a motivation, we content ourselves with sketching a proof using some general results on stratifications due to Cerf~\cite{Cerf}.
\begin{proof}[Proof of \cref{T:reduction to elemetary cusp merges} (sketch)]
Let us write~$\mc S^k$ for the subspace of~$\Cinfty(X,B)$ given by the maps of codimension~$k$ (as defined in \cref{ch:stability}).
This is completely analogous to the stratifications on spaces of real valued functions used by Cerf~\cite{Cerf}.
Conveniently, many of Cerf's general results about stratifications~\cite{Cerf}*{I.\parasign1--3} carry over to the surface valued setting. 
In particular, the results about maps of codimension at most one in \cref{ch:stable maps,ch:stable homotopies} show that we can apply the ``elementary path lemma''%
	\footnote{Named ``lemme des chemins élémentaires'' by Cerf.}%
~\cite{Cerf}*{p.20} to the stratification of~$\mc S^0\cup\mc S^1$ by the sets~$(\mc S^0,\mc S^1)$.
This reduces \cref{T:reduction to elemetary cusp merges} to the following problem:
	for any $f\in\mc S^1$ with a beak-to-beak singularity 
	we have to find an elementary cusp merge homotopy~$F=(f_s)_{s\in J}$ with~$f_0=f$.
But such a homotopy can easily be obtained using the normal form for beak-to-beak singularities, which allows to identify $f$ with~$\mu_0$ near its beak-to-beak, and extending~$f$ to a homotopy using a truncated version of the cusp merge model satisfying~$\mu_s=\mu_0$ outside of a \nbhd of the origin which can be obtained as in \cref{E:making compactly supported}.
\end{proof}

\subsection{Joining curves}
	\label{ch:joining curves}
We now focus on the second ingredient in \cref{T:main theorem:cusp merges and joining curves}.
Let~$F=(f_s)\in\CM_0(f)$ be an elementary cusp merge.
If we choose merge coordinates for~$F$, we can transfer the line segment~$L_0\subset\R^4$ defined  in~\eqref{eq:model joining curve} to an embedded arc~$L_F\subset X$ with endpoints on two cusps of~$f=f_{-1}$, say $p,q\in\Crit(f)$, which are eliminated by~$F$.
If we want to emphasize which cusps are eliminated we will use the more refined notation~$\CMel(f;p,q)$.
Since~$L_0$ can be understood as the trace of cusps in the merge model, we see that~$L_F$ is filled out by cusps of the map~$f_s$ for~$s<0$ and the beak-to-beak of~$f_0$.
Since the positions of the critical points of the map~$f_s$ are clearly independent of the choice of merge coordinates, so is~$L_F$.
We will call $L_F$ the \emph{(unparameterized) joining curve} of~$F$.
Alternatively, we could describe~$L_F$ without reference to a model as follows.
According to \cref{R:parameterized normal forms_codim0}, any cusp $c\in\Crit(f)$ will begin to evolve smoothly in~$X$ to a family of cusps~$c_s\in\Crit(f_s)$ as the parameter~$s$ increases from~$-1$ onward.
For $s<0$ no problems occur since each~$f_s$ is stable, but
for $s=0$ two things can happen.
In either case, the curve~$c_s$, $s<0$, has a limit~$c_0\in\Crit(f_0)$.
However, $c_0$~can either be a cusp, in which case $c_s$ will continue to be a cusp for all~$s$, or $c_0$~is the beak-to-beak point of~$f_0$ and the curve cannot be continued for~$s>0$.
The latter occurs for exactly the cusps~$p,q\in\Crit(f)$ and $L_F$~is traced out by the curves~$p_s$ and~$q_s$ for~$s\in[-1,0]$.
In fact, looking at the model again, we see that 
	\begin{equation}\label{eq:joining curve via cusp trace}
	\lambda_F\colon [-1,1]\lra X,
	\quad
	\lambda_F(\tau)=
	\begin{cases}
	p_{{-}\tau^2}, & \text{for $\tau\leq 0$} \\
	q_{{-}\tau^2}, & \text{for $\tau\geq0$}
	\end{cases}
	\end{equation}
constitutes a smooth parametrization of~$L_F$; 
we therefore call $\lambda_F$ the \emph{parameterized joining curve} of~$F$.
Note that the definition of~$\lambda_F$ implicitly involves the non-canonical choice of~$p$ as the first cusp.
We could just as well have chosen~$q$ and we would have obtained the reversed curve~$\lambda_F\rev(\tau)=\lambda_F(-\tau)$.
However, this choice turns out to be irrelevant, see \cref{R:order of cusps in joining curves}.
Soon we will see that~$\lambda_F$ (or even $L_F$) almost contain enough information to reconstruct~$F$ up to deformation within~$\CMel(f;p,q)$.
What is missing is a suitable notion of framing. 
For simplicity we will only focus on \emph{indefinite cusp merges},
	by which we mean cusp merges involving indefinite cusps, 
since these are the most important from the point of view of \wfs and broken Lefschetz fibrations.
However, our arguments can be modified to handle general cusp merges, see \cref{ch:general cusp merges}.
In the indefinite case, the framing will take the form of the line field spanned by the $\del_y$~vector field in some merge coordinates.
Unfortunately, the independence of the framing of the choice of coordinates is not as easy to see as in the case of~$L_F$ or~$\lambda_F$.
We therefore embark on a small digression before we continue the study of elementary cusp merges.

\subsubsection{Tangent spaces of indefinite cusps and beak-to-beaks}
	\label{ch:cusp observations}
Suppose first that a map $f\colon X\ra B$ has an indefinite cusp at~$p\in X$. 
%
%
%
For concreteness, we fix \emph{cusp coordinates} $(t,x,y,z)$ around~$p$ and~$(u,v)$ around~$f(p)$ such that $f$~is given by~$(t,x^3+3tx+y^2-z^2)$.
However, we will try to keep the discussion as intrinsic as possible.
Let $K_p=\ker(df_p)$ and $Q_p=\coker(df_p)$.
Our first observation is that the image of~$df_p$, which is 1--dimensional, has a preferred orientation induced by the orientations of~$X$ and~$B$.
Intuitively, this is determined by the direction in which the cuspidal tip of the discriminant points in~$B$.
In the coordinates this direction agrees with the ray~$D_p\subset T_{f(p)}B$ spanned by the coordinate vector field~$\del_u$.
Moreover, the results in~\cite{Levine}*{(4.2)} show that $D_p$ is independent of the choice of coordinates and therefore well-defined.
We will refer to~$D_p$ as the \emph{direction of the cusp}.
Following Levine, we say that a tangent vector~$v\in T_pX$ \emph{points downward} or \emph{upward} at~$p$ if $df_p(v)\in D_p\setminus\{0\}$ or $df_p(-v)\in D_p\setminus\{0\}$, respectively.
Next we note that the orientation of~$B$ together with the direction of the cusp induces an orientation on~$Q_p$, which is also 1--dimensional, by extending a non-zero vector in~$D_p$ to an oriented basis of $T_{f(p)}B$.
We therefore have a notion of positivity in~$Q_p$.
In coordinates, we can identify~$Q_p$ with the span of~$\del_v$.
Lastly, we consider the \emph{intrinsic second derivative} of~$f$ which gives a symmetric bilinear map
	\begin{equation*}
	\delta^2_pf\colon K_p\times K_p \lra Q_p
	\end{equation*}
and can be considered as a generalization of the Hessian for functions.
In fact, in coordinates $K_p$ is spanned by $\del_x$, $\del_y$, and $\del_z$ and in this basis $\delta^2_pf$~appears as the Hessian of the function $x^3+y^2-z^2$ using the identification~$Q_p\cong\R\del_v$.
We can therefore consider $\delta^2_pf$ as a symmetric bilinear form which provides a decomposition of~$K_p$ into three disjoint sectors $K^+_p$, $K^-_p$, and~$K^0_p$ defined as the set of vectors $v\in K_p$ for which $\delta^2_pf(v,v)$ is positive, negative, and zero, respectively.
From the coordinate description of~$\delta^2_pf$ we see that $\del_x\in K^0_p$, $\del_y\in K^+_p$, and~$\del_z\in K^-_p$.
Moreover, $K^+_p$ and~$K^-_p$ each have two connected components which are convex and interchanged by the multiplication with any negative number.
In particular, the image of~$K^\pm_p$ in~$\PP(K_p)$, henceforth denoted by~$\PP(K^\pm_p)$, is contractible.
The above discussion of~$\delta^2_pf$ goes through verbatim if~$f$ has a beak-to-beak at~$p$ and thus has the local model $(t,x^3-3t^2x+y^2-z^2)$.
However, this time there is no preferred orientation for~$Q_p$ but~$\delta^2_pf$ is still defined and singles out a well defined null sector~$K^0_p$ whose complement in~$K_p$ has four contractible components which are paired by comparing the sign of~$\delta^2_pf$ \wrt a fixed orientation of~$Q_p$.

\subsubsection{Abstract joining curves and framings}
	\label{ch:joining curves abstractly}
We now give an abstract definition of joining curves which is essentially due to Levine~\cite{Levine}*{(4.4)}.
However some adjustments were necessary since we are working with stable maps while Levine allows arbitrary multi-singularities made of folds and cusps.
Throughout this subsection we consider a stable map $f\colon X\ra B$ and a pair of cusps~$p,q\in\Crit(f)$.
\begin{definition}[Joining curves]\label{D:joining curves}
A curve~$\lambda\colon [-1,1]\ra X$ with $\lambda(-1)=p$ and $\lambda(1)=q$ is called a \emph{joining curve} for~$f$ (from~$p$ to~$q$) if
\begin{enumerate}[(i)]
	\item 
		both $\lambda$ and~$f\circ\lambda$ are smooth embeddings, 
	\item
		$\lambda$ points {downward} at~$p$ and {upward} at~$q$, and
 	\item
		the interior of~$\lambda$ does not meet any singular fibers.
\end{enumerate}
The joining curves for~$f$ form a subspace of~$\Cinfty([-1,1],X)$ which we denote by~$\JC(f)$, and by~$\JC(f;p,q)$ indicates fixed cusps at the endpoints.
\end{definition}
%
%
%
It is easy to see that~$\JC(f;p,q)$ is non-empty if and only if $f(p)$ and~$f(q)$ can be connected by an arc of regular values, say~$R\subset B$, and $p$~and~$q$ are in the same connected component of~$f\inv(R)$.
%
%
Moreover, Levine shows that if~$\JC(f;p,q)$ is non-empty and if $p$~and~$q$ form a ``matching pair'' (see~\cite{Levine}*{p.284}), then they can be eliminated by a cusp merge homotopy.
We will not explain the matching pair condition, we only observe that it is satisfied for any two indefinite cusps so that the existence of a joining curve guarantees that the pair can be eliminated.
Next we have to discuss a notion of framings for joining curves which are slightly different for definite and indefinite cusps.
Again, we only treat the indefinite case and refer to \cref{ch:general cusp merges} for some remarks about other cases.
Consider the subspace of~$TX$ given by~$\ker(df)=\cup_p\ker(df|_p)$.
Although $\ker(df)$ is not a vector bundle over~$X$ it is still meaningful to speak of sections and pullbacks.
Moreover, since the fibers are vector spaces (albeit of varying ranks) we can take the fiber-wise projectivization which we denote by~$\bs X=\PP(\ker(df))\subset\PP(TX)$ and we let $\pi\colon\bs X\ra X$ be the obvious projection.
\begin{definition}[Framed joining curves]\label{D:framed joining curves}
Let~$\lambda\in\JC(f;p,q)$ be a joining curve.
A \emph{framing} for $\lambda$ is a map~$\blam\colon J\ra\bs X$ such that~$\pi\circ\blam=\lam$ and $\bs\lambda(\pm1)\in\PP\big(K^\pm_{\lambda(\pm1)}\big)$.
Any map $\bs\lambda\colon J\ra\PP(\ker(df))$ with the latter property and $\pi\circ\bs\lambda\in\JC(f)$ is called a \emph{framed joining curve}.
We denote by $\JCfr(f;p,q)$ and~$\JCfr(f)$ the spaces of framed joining curves with or without fixed endpoints.
\end{definition}
In what follows it will be understood that whenever we discuss a framed joining curve~$\bs\lam$, the underlying joining curve is denoted by~$\lam=\pi\circ\bs\lam$.
\begin{remark}\label{R:framed joining curves unraveled}
With the arguments in \cref{ch:mapping class group interpretation} in mind, we take a moment to unravel \cref{D:framed joining curves}.
Let~$\lam\in\JC(f)$ be a joining curve and let $\S_\tau$ be the fiber of~$f$ containing~$\lambda(\tau)$, that is, $\S_\tau=f\inv\big(f(\lambda(\tau))\big)$.
For~$\tau\neq\pm1$, $\S_\tau$~is a smooth submanifold of~$X$ whose tangent bundle is the restriction of~$\ker(df)$ to~$\S_\tau$.
Therefore a framing of~$\lambda$ amounts to the choice of lines in~$T_{\lam(\tau)}\Sigma_\tau$ with limits in the negative sector $K^-_{p}$ over the negative end~$p=\lam(-1)$ and in~$K^+_q$ over~$q=\lam(1)$.
\end{remark}

\subsection{From \texorpdfstring{$\CM_0(f)$}{elementary cusp merges} to \texorpdfstring{$\JCfr(f)$}{framed joining curves}}
	\label{ch:CM to JC}
Now that we have proper definitions for the objects in the statement of \cref{T:main theorem:cusp merges and joining curves}, we have to understand how they are related.
A first step in this direction was already taken in the beginning of \cref{ch:joining curves} where we discussed the parameterized and unparameterized joining curves~$\lambda_F$ and~$L_F$ of an elementary cusp merge homotopy~$F\in\CMel(f)$.
We leave the following easy verifications to the reader.
\begin{lemma}\label{T:CM to JC}
The curve $\lambda_F$ defined in~\eqref{eq:joining curve via cusp trace} is a joining curve for~$f$ in the sense of \cref{D:joining curves}.
Moreover, in merge coordinates for~$F$ the line field spanned by~$\del_z$ constitutes a framing for~$\lambda_F$.
\end{lemma}
As mentioned before, there is no canonical choice of merge coordinates and therefore no canonical framing.
However, we we can get around this issue as follows.
\begin{lemma}\label{T:framing ambiguity}
There is a connected space of preferred framings of~$\lambda_F$ induced by~$F$ which contains the framings coming from merge coordinates.
In particular, an elementary cusp merge homotopy~$F\in\CMel(f;p,q)$ gives rise to a well defined element~$\bs\lambda_F\in\pi_0\big(\JCfr(f;p,q)\big)$.
\end{lemma}
\begin{proof}[Proof of \cref{T:framing ambiguity}]
For simplicity, we write~$\lambda=\lambda_F$.
As usual, we consider~$F$ as a family~$(f_s)$ of maps.
The discussion of the tangent spaces of cusps and beak-to-beaks in \cref{ch:cusp observations} results in a collection of 3--dimensional vector spaces
	\begin{equation}\label{eq:cusp kernel bundle}
	K_{\tau} = 
	\ker\big( d(f_{-|\tau|})_{\lambda(\tau)} \big) \subset T_{\lambda(\tau)}X.
	\end{equation}
In merge coordinates these are spanned by the vector fields~$\del_x$, $\del_y$ and~$\del_z$ and therefore fit together to a smooth vector sub-bundle~$K\subset \lambda^{\ast}TX$.
The crucial observation is that for~$\tau\neq0$ the sectors $K^\pm_{\tau}$ and $K^\mp_{-\tau}$ defined by the cusps of~$f_{-|\tau|}$ as in \cref{ch:cusp observations} limit to the same sector of~$K_0$ as $\tau$~approaches zero from either side.
This can again be checked in the model.
The change of signs is caused by the fact that~$\lambda$ {points} downward at one cusp and upward at the other.
As a consequence, we obtain a sub-bundle $\wt{\mc F}_\lambda$ of the projectivization~$\PP(K)$ such that for $\tau >0$ the fibers of $\wt{\mc F}_\lambda$ over $\tau$ and $-\tau$ respectively agree with $\PP(K^+_{\tau})$ and $\PP(K^-_{-\tau})$.
Note that $\wt{\mc F}_\lambda$ has contractible fibers. 
In order to make the connection to framings we consider the intersection of~$K$ with the kernel of~$df_{-1}$ along~$\lambda$.
First we observe that~$K$ and~$\lambda^\ast\ker(df_{-1})$ agree over the end-points of~$\lambda$. 
Over the interior, another local calculation shows that~$\ker(df_{-1})$ is spanned by~$\del_y$ and~$\del_z$ along~$\lambda$, thus $\lambda^*\ker(df_{-1})$ forms a subset of~$K$. 
Finally, after passing to projectivizations we arrive at a subset
	\begin{equation*}
	\mc F_\lambda = \wt{\mc F}_\lambda \cap \PP(\lambda^*\ker(df_{-1})) \subset \PP(K).
	\end{equation*}
The proof is finished by the following consequences of the (admittedly somewhat convoluted) definitions and the arguments in \cref{ch:cusp observations}:
%
		the map~$\mc F_\lambda\ra [-1,1]$ has contractible fibers,
		sections of~$\mc F_\lambda$ can be considered as framings of~$\lambda$,
		and any choice of merge coordinates give rise to a section of~$\mc F_\lam$ via~$\del_z$.
\end{proof}
We arrive at the main conclusion of the present subsection.
\begin{proposition}\label{T:CM->JC up to homotpy}
The assignment $F\mapsto\bs\lambda_F$ descends to a map
	\begin{equation*}
	\pi_0\big(\CMel(f;p,q)\big) \lra \pi_0\big(\JCfr(f;p,q)\big).
	\end{equation*}
\end{proposition}
\begin{proof}
We have to show that the~$\bs\lambda_F$ does not change in smooth families.
So let~$\bs F=(F_r)_{r\in[0,1]}$ be a smooth 1--parameter family where~$F_r\in\CMel(f;p,q)$ and let~$\lambda_r$ be the joining curve of~$F_r$.
Since each $F_r$ is stable as a homotopy, it follows from the \cref{R:parameterized normal forms_codim1} that the family~$(\lambda_r)_{r\in[0,1]}$ depends smoothly on~$r$.
Moreover, the discussion of framings in the proof of \cref{T:framing ambiguity} is only notationally more complicated in the presence of a parameter.
\end{proof}

\subsection{From \texorpdfstring{$\JCfr(f)$}{framed joining curves} to \texorpdfstring{$\CM_0(f)$}{elementary cusp merges}}
	\label{ch:JC to CM}
So far we have been able to extract framed joining curves from elementary cusp merge homotopies.
Now we try to go back and construct homotopies with prescribed joining curves and framings.
As mentioned earlier, the construction is essentially due to Levine~\cite{Levine}.
\begin{lemma}[Parametric Levine construction]\label{T:Levine with parameters}
Let~$\bs\lambda_r\in\JCfr(f;p,q)$ be a smooth family of framed joining curves.
Then there exists a smooth family~$F_r\in\CMel(f;p,q)$ of elementary cusp merges such that $\bs\lambda_r$ is the framed joining curve of~$F_r$.
\end{lemma}
\begin{proof}[Proof of \cref{T:Levine with parameters}]
We first consider a single framed joining curve~$\bs\lam$, that is, the case of a constant homotopy.
Let~$L\subset X$ denote the image of the underlying joining curve~$\lambda$ and let~$L_0\subset\R^4$ as in~\eqref{eq:model joining curve}.
The proof has two main steps.
\begin{description}
	\item[Step 1]
		Identify $f$ with the initial map~$\mu_{-1}$ of the cusp merge model~\eqref{eq:cusp merge model} in \nbhds of~$L$ and~$L_0$.
	\item[Step 2]
		For any \nbhd $\nu L_0$ of~$L_0$ in $\R^4$ construct an elementary cusp merge homotopy with initial map~$\mu_{-1}$ which is constant outside of~$\nu L_0$.
\end{description}
Once this is done, we can promote the homotopy constructed in the second step to an elementary cusp merge homotopy from~$f$ using the identification obtained in the first step in the obvious way.
Let us call this homotopy~
	\begin{equation*}
	F_{\bs\lambda}\in\CM_0(f;p,q).
	\end{equation*}
Step 1 is treated carefully in~\cite{Levine}*{(4.6),(4.8)} to which we refer for more details.
The basic observation is that the map~$\mu_{-1}\colon\R^4\ra\R^2$ can be divided into three parts:
	an inner part where~$|t|<1$ consisting of a tubular \nbhd of~$L_0$ fibered over a tubular \nbhd of~$\mu_{-1}(L_0)$,
	and two outer parts where $|t|\geq1-\epsilon$ which are both easily matched with the cusp model.
So the idea is to implant these structures into a \nbhd of~$L$ in~$X$.
To begin with, we choose cusp coordinates around the points~$p$ and~$q$.
We can assume that the {$\del_z$}~vector fields of both coordinate patches take values in the tangent lines specified by~$\bs\lam$ near the ends of~$L$.
Let us write~$\del_z^{(p)}$ and~$\del_z^{(q)}$ in order to distinguish between the two.
Next we use a partition of unity argument to extend~$\del_z^{(p)}$ to 
	a vector field~$\nu_z$ on~$X$ which 
	takes values in~$\bs\lam$ along~$L$ and 
	agrees with~$\pm\del_z^{(q)}$ in the other coordinate patch.
If the sign is negative, then we change the coordinates near~$q$ by the symmetry $(t,x,y,z)\mapsto(t,x,-y,-z)$.
After this is done, we can also find vector fields~$\nu_x$ and~$\nu_y$ on~$X$ that agree with the $\del_x$ and $\del_y$~fields in the coordinate patches and are point-wise linearly independent near~$L$.
Lastly, we choose a vector field~$T$ on~$X$ that agrees with~$\del_t^{(p)}$ near~$p$ and with $-\del_t^{(q)}$ near~$q$ and has~$L$ as a flow line.
(In particular, they form an honest \emph{normal framing} for~$L$ as a submanifold of~$X$.)
If the extensions $T$, $\nu_x$, $\nu_y$, and $\nu_z$ were chosen carefully enough, then we can use their flows to connect the two cusp coordinate systems to a single coordinate patch containing~$L$ in which $f$~appears as~$\mu_{-1}$ and~$L$ corresponds to~$L_0$.
At this point, we already have constructed enough of~$F_{\bs\lam}$ to form its joining curve whose image agrees with~$L$ by construction, and possibly after rescaling the coordinates in the $t$--direction we can assume that the parametrization agrees with~$\lambda$ as well.

As for Step~2, we refer to \cref{E:making compactly supported}.  
We can obtain an elementary cusp merge homotopy with the desired property by taking sufficiently small $K_\varepsilon$ and a family $S_{s,t}$ sufficiently close to $(t,0,0,0)$ as in \cref{E:making compactly supported} so that the support of the resulting homotopy is contained in the given neighborhood $\nu L_0$. 
Lastly, in order to treat the case of non-constant homotopies of joining curves, it is enough to have a parametric version of the first step.
But a closer investigation of Levine's arguments shows that the essential ingredients are the normal form for a cusp and the tubular \nbhd theorem which both have parametric versions.
\end{proof}
\begin{remark}\label{R:order of cusps in joining curves}
It turns out that the Levine construction applied to~$\bs\lam\in\JCfr(f;p,q)$ and its reverse~$\bs\lam\rev(\tau)=\bs\lam(-\tau)$, which lies in~$\JCfr(f;q,p)$ give rise to the same cusp merge homotopy.
The reason is that the map~$\mu_{-1}$ admits a symmetry that interchanges the two cusps at~$(\pm1,0,0,0)$.
For example, one can take the linear diffeomorphisms $(t,x,y,z)\mapsto(-t,-x,z,y)$ on~$\R^4$ and~$(u,v)\mapsto(-u,-v)$ on~$\R^2$.
\end{remark}
\begin{remark}\label{R:Levine model}
Note that Levine takes a slightly different approach to the second step in the above proof (see \cite{Levine}*{(4.9)}).
Unfortunately, his version does not lead to elementary cusp merge homotopies in the sense of \cref{D:elementary cusp merges}, as it results in a slightly more complicated model than~$(\mu_s)$ from~\eqref{eq:cusp merge model}.
Since we found the model~$(\mu_s)$ much more convenient for computations, we were forced to develop the machinery in \cref{ch:normalization}.
\end{remark}
The main goal of this section is to establish the following.
\begin{proposition}\label{T:JC->CM up to homotopy}
The parametric Levine construction induces a map
	\begin{equation*}
	\pi_0\big(\JCfr(f;p,q)\big) \lra \pi_0\big(\CMel(f;p,q)\big)
	\end{equation*}
\end{proposition}
\begin{proof}
The Levine construction involves several choices and we have to make sure that for each~$\bs\lam\in\JCfr(f;p,q)$ we get a well defined element of~$\pi_0\big(\CMel(f;p,q)\big)$.
Let~$F,F'\in\CMel(f;p,q)$ be two elementary cusp merges resulting from different choices in the Levine construction.
The ambiguities in the first step of the construction affect the behavior of the homotopies near the image of~$\bs\lam$, while the ambiguity in the second step concerns the extension to all of~$X$.
The first ambiguities can be eliminated using our results about the symmetries of cusps in~\cref{T:symmetries of folds and cusps} combined with the standard fact that 
	a parameterized tubular \nbhd of a submanifold is determined by a normal framing.
These results show that we can deform~$F'$ so that it agrees with~$F$ near the image of~$\bs\lam$ through homotopies arising from Levine constructions.
Instead of dealing with the second type ambiguities directly, we observe that the desired result follows from \cref{T:isotopy lemma} below.
\end{proof}
\begin{lemma}\label{T:isotopy lemma}
Suppose that $F,F'\in\CM_0(f)$ agree near their joining curves.
Then they can be deformed into each other within~$\CM_0(f)$.
\end{lemma}

\begin{proof}[Proof]
Let~$L\subset X$ be the image of the common joining curves of~$F$ and~$F'$.
By assumption we can find nested neighborhoods $U\subset U'$ of $L$ such that $F$ and $F'$ are constant in $U'$ and $U$ is contained in a common set of merge coordinates of $F$ and $F^\prime$.
We take other neighborhoods $V\subset U$ and $V'\supset U'$ of $L$ and $U'$, respectively. 
For convenience, we assume that~$U, f(U), V$ and~$f(V)$ are all open balls with compact closures and smooth boundaries.
As we construct $1$--parameter families of diffeomorphisms in the proof of \cref{T:normalization lemma} we can take codimension $0$ compact submanifolds $W\subset  [-1,\varepsilon]\times V$ and $K\subset f(V)$ which contain $ [-1,\varepsilon]\times L$ and $f(L)$, respectively, an open neighborhood $W'\subset [-1,\varepsilon]\times V$ of $\Pa W\setminus \left( \{-1,\varepsilon\}\times V\right)$, level preserving diffeomorphisms $\Phi, \Phi'$ of $[-1,\varepsilon]\times X$ and $\Psi, \Psi'$ of $[-1,\varepsilon]\times B$ such that 
\begin{itemize}

\item
$\wt{F}=\Psi \circ (\id\times f) \circ \Phi\inv$ and $\wt{F}'=\Psi' \circ (\id\times f) \circ {\Phi'}\inv$ on $W' \cup W^c$, 

\item
Supports of $\Phi$ and $\Phi'$ are contained in $V'\times [-1,\varepsilon]$, while those of $\Psi$ and $\Psi'$ are contained in $[-1,\varepsilon]\times\nu\Pa K$ where $\nu\Pa K$ is a neighborhood of $\Pa K$, 

\item
$\Phi = \Phi'$ on $V\times [-1,\varepsilon]$ and $\Psi = \Psi'$ on $f(V)\times [-1,\varepsilon]$. 

\end{itemize}
\noindent
We can consider the composition  $(\Phi\inv\Phi',\Psi\inv\Psi')$ as a map from~$[-1,\varepsilon]$ into the space of pairs $(\phi,\psi)\in\Diff(X)\times\Diff(B)$ such that $\phi|_{V}$ and $\psi|_{f(V)}$ are the identity maps and $(\psi\circ f \circ\phi\inv)|_{X\setminus V}=f|_{X\setminus V}$ that sends~$-1$ to~$(\id_{X},\id_{B})$.
In other words, $(\Phi\inv\Phi',\Psi\inv\Psi')$ corresponds to a path in the space just mentioned based at~$(\id_{X},\id_{B})$.
Since we do not have to keep the endpoint of this path fixed, we can deform it to the constant path.
But this means that we can find a 1--parameter family of pairs~$(\Phi_r,\Psi_r)$, $r\in[0,1]$, that connects~$(\Phi,\Psi)$ to~$(\Phi',\Psi')$ with the properties that $\Phi_r$ and~$\Psi_r$ are independent of~$r$ in~$[-1,\varepsilon]\times V$ and $[-1,\varepsilon]\times f(V)$ and the support of $\Phi_r$ (resp.~$\Psi_r$) is contained in $V'\times[-1,\varepsilon]$ (resp.~$\nu\Pa K\times[-1,\varepsilon]$).
Using this we obtain a 1--parameter family of homotopies~$F_r$ with unfoldings 
	\begin{equation*}
	\wt{F}_r=
	\begin{cases}
	\Psi_r \circ (\id\times f) \circ \Phi_r\inv & \text{on $W'\cup W^c$} \\
	\wt{F}=\wt{F}' & \text{on $W$}
	\end{cases}
	\end{equation*}
which provides a deformation of~$F=F_0$ into~$F'=F_1$.
It remains to show that each~$F_r$ is an elementary cusp merge homotopy.
But this is clear from the construction.
Indeed, on~$W$ and $\widetilde{F}(W)$ the common merge coordinates for~$F$ and~$F'$ can be used for each~$F_r$.
Moreover, $\Phi_r$ and $\Psi_r$, which are supported on $V'\times[-1,\varepsilon]$ and $\nu\Pa K\times[-1,\varepsilon]$, satisfy the second condition in the definition of elementary cusp merges.
\end{proof}

\subsection{The proof of \texorpdfstring{\cref{T:main theorem:cusp merges and joining curves}}{Theorem A}}
According to \cref{T:CM->JC up to homotpy,T:JC->CM up to homotopy}, we have natural maps between the sets~$\pi_0\big(\CM_0(f)\big)$ and~$\pi_0\big(\JCfr(f)\big)$.
We will now show that they are mutually inverse.
We briefly recall the constructions:
\begin{itemize}
	\item 
		Given $\bs\lambda\in\JCfr(f)$ we can find coordinates around the image of~$\lambda$ in which $f$ is represented by~$\mu_{-1}$ and~$\bs\lambda$ corresponds $\bs\lambda_0$, that is,  $\lambda_0(\tau)=(\tau,0,0,0)\in\R^4$ framed by the lines spanned by~$\del_y$. 
		A compactly supported version of the merge model then gives rise to~$F_{\bs\lambda}\in\CM_0(f)$.
	\item
		Given~$F\in\CM_0(f)$ we choose merge coordinates for~$f$ and obtain~$\bs\lambda_F\in\JCfr(f)$ by pulling back~$\bs\lambda_0$ to~$X$.
\end{itemize}
From \cref{T:Levine with parameters} we already know that $\bs\lambda_{F_{\bs\lambda}}$ is equal to~$\bs\lambda$.
So we only have to show that $F_{\bs\lam_F}$ and~$F$ represent the same element of~$\pi_0\big(\CM_0(f)\big)$.
But this follows from the same arguments as in the proof of \cref{T:JC->CM up to homotopy}.
Indeed, since $F_{\bs\lam_F}$ and~$F$ have the same framed joining curve, we can deform them to agree near the joining curve and then apply \cref{T:isotopy lemma}.
This finishes the proof of \cref{T:main theorem:cusp merges and joining curves}.

\subsection{More general cusp merges}
	\label{ch:general cusp merges}
While we have stated and proved \cref{T:main theorem:cusp merges and joining curves} only for the case of indefinite cusp merges in dimensions~$(4,2)$, our arguments can easily be modified to obtain an analogous result for general cusp merges in dimensions~$(n,2)$ with~$n\geq 2$.
Essentially the only difference arises in the discussion of framings of joining curves.
Indeed, the general cusp merge model in dimensions~$(n,2)$ takes the form
	\begin{equation*}
	(s;t,x,y_1,\dots,y_k,z_1,\dots,z_l)
	\mapsto
	(t,x^3-3(t^2+s)x +\sum_i {y_i}^2 - \sum_i{z_i}^2)
	\end{equation*}
where $k+l=n-2$ and the framing has to keep track of  $\del_{y_1},\dots,\del_{y_k},\del_{z_1},\dots,\del_{z_l}$.
With an appropriate notion of framings, the parametric Levine construction and all other ingredients in the proof of \cref{T:main theorem:cusp merges and joining curves} are available in the general case.

%% file: parallel_transport_short.tex
\section{Parallel Transport and Vanishing Sets}
	\label{ch:vanishing sets chapter}
Before we can address \cref{T:main theorem:cusp merges and MCG} we have to develop a theory of parallel transport.
We start in the setting of general smooth maps and then specialize to wrinkled fibrations on 4--manifolds.
Besides the relevance of \cref{T:main theorem:cusp merges and MCG}, our goal is to provide a general framework for discussing vanishing sets of \wfs which, as far as we know, has not been available.

\subsection{Connections for smooth maps}
	\label{ch:connections}
For the moment, we consider an arbitrary smooth map $f\colon N\ra P$ between connected manifolds of dimensions~$n\geq p$.
We need a tool to compare different fibers of~$f$.
The following definition should be thought of as a hybrid between connections in fiber bundles and gradients of (Morse) functions.
\begin{definition}[Connections]\label{D:connection}
A \emph{connection} for a smooth map~$f\colon N\ra P$ is a subset~$\H\subset TN$ obtained as the point wise orthogonal complement of~$\ker(df)$ \wrt some Riemannian metric on~$N$.
\end{definition}
Observe that if~$f$ is a fiber bundle, then we recover the usual notion of connection, while if~$f$ is a real-valued function, then a connection is simply the point wise span of the gradient of~$f$ \wrt some metric.
Recall that the central idea in the construction of parallel transport in fiber bundles is to consider lifts of curves in~$P$ to curves in~$X$ via~$f$.
We will follow a similar procedure; for simplicity, we restrict our attention to a special class of curves in~$P$.
\begin{definition}[Reference arcs and~$\mc H$--lifts]\label{D:reference arcs}
Let~$f\colon N\ra P$ be a smooth map.
\begin{enumerate}[(a)]
	\item 
		A \emph{reference arc} for~$f$ is an embedded arc~$\gamma\colon[0,1]\ra P$ such that~$\gamma(0)$ and~$\gamma(1)$ are regular values of $f$.
	\item
		Let~$\H$ be a connection for~$f$.
		A map~$\tilde{\g}\colon O\ra N$ is called an \emph{$\H$--lift} of a curve~$\g\colon[0,1]\ra P$ if
			$O\subset[0,1]$ is a relatively open subinterval,
			$f\circ \tilde{\gamma}=\g|_O$, and 
			the velocity vectors of $\tilde{\gamma}$ are contained in~$\H$.
\end{enumerate}
\end{definition}
We will occasionally blur the notational distinction between a reference arc~$\gamma$ and its image in~$P$.
For example, we usually write~$f\inv(\g)$ instead of~$f\inv(\g([0,1]))$.
Moreover, we will often abbreviate the fibers over~$\gamma$ by~$\S_t=f\inv(\g(t))$.
The following technical result is the key ingredient for the parallel transport construction.
\begin{proposition}\label{T:basic dynamics}
Let~$f\colon N\ra P$ be a proper map equipped with a connection~$\H$.
Furthermore, let~$\gamma\colon[0,1]\ra P$ be a reference arc for $f$ and $p\in\Sigma_t$ a regular point for some fixed $t\in [0,1]$.
\begin{enumerate}[(i)]\itemsep3pt
	\item 
		There exists a unique~$\H$--lift of~$\g$, denoted by
			\begin{equation*}
			\tilde{\gamma}_{t,p}^\H\colon O_{t,p}\lra N,
			\quad 
			O_{t,p}\subset[0,1],
			\end{equation*}
		where $O_{t,p}$ is a relatively open interval containing~$t$ such that~$\tilde{\gamma}_{t,p}^\H(t)=p$ and all other $\H$--lifts of~$\g$ with this property are restrictions of~$\tilde{\gamma}_{t,p}^\H$.
	\item
		If~$\H$ depends smoothly on some auxiliary parameters, then so does~$\tilde{\gamma}_{t,p}^\H$.
	\item
		If each fiber of~$f$ along~$\g$ contains at most finitely many critical points of~$f$, then $\tilde{\gamma}_{t,p}^\H$ limits to a critical point on each open end of~$O_{t,p}$.
\end{enumerate}
\end{proposition}
\begin{proof}[Proof of \cref{T:basic dynamics}]
This is well known for submersions and the usual proofs of~(i) and~(ii) are easily adapted to take critical points into account, so we shall be brief.
Since $f|_{M\setminus\Crit_f}$ is a submersion (albeit non-proper), the set $Y=f\inv(\g)\setminus\Crit_f$ is a non-compact smooth $(n-p+1)$--manifold with boundary~$\S_0\amalg\S_1$.
For~$q\in\S_t$ let~$\Gamma(q)\in\H$ be the unique element that is mapped to the velocity vector~$\dot{\gamma}(t)$ under~$df$.
It is easy to see that $\Gamma$~is tangent to~$Y$ and can thus be considered as a vector field on~$Y$.
The claims~(i) and~(ii) now follow from standard properties of the flow of~$\Gamma$.
In order to prove~(iii), let~$\tilde{\g}$ be some \mbox{$\H$--lift} of~$\g$ defined on an open interval~$(a,b)\subset[0,1]$.
It suffices to show that~$\tilde{\gamma}$ can be extended to the closed interval~$[a,b]$.
Since the arguments for~$a$ and~$b$ are virtually the same, we only focus on~$b$.
Let~$\{b_n\}\subset(a,b)$ be a sequence that converges to~$b$.
Since $f\inv(\gamma)$ is compact and~$\g$ is continuous, the sequence $\tilde{\gamma}(b_n)$~has an accumulation point~$\tilde{b}\in N$ with~$f(\tilde{b})=\gamma(b)$;
in other words, we have~$\tilde{b}\in\S_b$.
Note that if $\tilde{b}$~is a regular point of~$f$, then $\tilde{\gamma}$ must be the restriction of the maximal lift~$\tilde{\gamma}_{b,\tilde{b}}^\H$ which obviously extends $\tilde{\gamma}$ to~$b$.
On the other hand, if $\tilde{b}$~is a critical point of~$f$, then the same reasoning provides an accumulation point of the sequence~$\{\tilde{\gamma}(b_n)\}$ and we only have to show that it is unique.
We argue by contradiction and assume that there are two or more accumulation points.
The above arguments show that all accumulation points have to be critical points of~$f$ in the fiber~$\S_b$ which are isolated by assumption.
Let~$U\subset\S_b$ be the union of pairwise disjoint open \nbhds~$U_i$ of all critical points in~$\S_b$.
Since~$\tilde{\g}$ must enter and leave some of $U_i$ infinitely many times, we can find a subsequence of~$\{{c}_k\}\subset\{{b}_n\}$ with $c_k\ra b$ such that~$\tilde{\g}(c_k)\in\S_b\setminus U$. 
Now, $\S_b\setminus U$ is also compact so that $\{\tilde{g}(c_k)\}$ must have an accumulation point.
But this would have to be a regular point of~$f$ which is impossible.
\end{proof}
Note that the finiteness assumption in~(iii) of \cref{T:basic dynamics} is automatically satisfied if~$f$ is stable.
In particular, the conclusion always holds for \wfs.

\subsection{Vanishing sets and parallel transport}
	\label{ch:parallel transport}
We continue with the notation of~\cref{T:basic dynamics}.
For a critical point~$q\in f\inv(\g)\cap\Crit(f)$ we say that 
	\emph{$\tilde{\gamma}_{t,p}^\H$~runs into~$q$} (or \emph{emerges from~$q$}) 
if its left (or right) limit is $q$. 
We can now define what we call the \emph{vanishing sets} of the triple~$(\g,\H;q)$ as
	\begin{align*}
	V_0(\gamma,\H;q)&= \big\{ p\in\S_0 \big|\; \text{$\tilde{\gamma}_{0,p}^\H$ runs into~$q$} \big\} \subset\S_0\\
	V_1(\gamma,\H;q)&= \big\{ p\in\S_1 \big|\; \text{$\tilde{\gamma}_{1,p}^\H$ emerges from~$q$} \big\} \subset\S_1.
	\end{align*}
This is reminiscent of the theory stable and unstable manifolds of gradients of critical points of functions which is already a delicate subject for degenerate critical points (see \cref{eg:Takens} below).
However, our situation can be even more complicated since~$f\inv(\g)$ does not have to be a manifold.
So we should be alarmed and proceed with caution.
In particular, the following example shows that the vanishing sets can exhibit a delicate dependence on the choice of connection.
\begin{example}\label{eg:Takens}
In the case of a real-valued function the vanishing sets can be considered as the intersections of stable and unstable manifolds associated with critical points with level sets.
In particular, for a reference arc whose preimage contains a single non-degenerate critical point we recover the well known ascending and descending spheres prominently featured in Morse theory and one can show that they are independent of the choice of connection up to ambient isotopies within the levels.
However, for degenerate critical points the topological type of the stable and unstable manifolds can depend on the choice of metric, it can even jump in 1--parameter families of metrics (see~\cite{Takens} for an example).
\end{example}
Continuing with the previous discussion, let us write~$V_0\subset\S_0$ and~$V_1\subset\S_1$ for the unions of~$V_i(\g,\H;q)$ over all~$q\in f\inv(\g)\cap\Crit(f)$.
Of course, the unions $V_i$ depend even more drastically on~$\H$ than the vanishing sets alone as is easily seen by considering two critical points of a Morse function.
Regardless, we can define the \emph{parallel transport} along~$\gamma$ \wrt a fixed connection~$\H$ by
	\begin{equation*}
	\PT^\H_\gamma\colon \S_0\setminus V_0 \ra \S_1 \setminus V_1
	,\quad\quad 
	p\mapsto\tilde{\gamma}_{0,p}^\H(1).
	\end{equation*}
Note that this is essentially the same construction as in the case of fiber bundles, the price to pay for the presence of critical points is that the parallel transport is only partially defined.
The following observations are either obvious or follow from the standard theory of ordinary differential equations.
\begin{itemize}
	\item 
		$\PT_\g^\H$ is a diffeomorphism.
	\item
		If two reference arcs $\g$ and $\g'$ have the same image, then they give rise to the same vanishing sets. 
		Furthermore, $\PT_\g^\H$ is equal to $\PT_{\g'}^\H$.
	\item
		If~$\g$ is an arc of regular values, then~$\PT_\g^\H\colon\S_0\ra\S_1$ is everywhere defined and its isotopy class is independent of~$\H$ and depends only on the homotopy class of~$\g$ relative to its endpoints in the connected component of~$P\setminus f(\Crit(f))$ containing~$\g$.
	\item
		The whole discussion generalizes immersed, piecewise smooth, and closed curves.
		For parallel transport along a closed curve of regular values we use the term \emph{monodromy}.
\end{itemize}
\begin{remark}\label{R:isotopies of vanishing sets}
As an addendum to the third statement, we note that if~$\g$ is an arc of regular values, then for any diffeomorphism~$\phi\in\Diff(\S)$ that is isotopic to the identity one can find a connection~$\H'$ which agrees with~$\H$ outside an arbitrarily small \nbhd of~$f\inv(\g)$ such that~$\PT_\g^{\H'}=\PT_\g^\H\circ\phi$.
(This follows as in the proof of~\cite{Nicolaescu}*{Lemma~2.28}, for example.)
Applied to regular parts of general reference arcs this provides useful flexibility for rearranging the vanishing sets.
\end{remark}

\subsection{Vanishing Sets in Wrinkled Fibrations}
	\label{ch:wrinkled fibrations_vanishing sets}
We now leave the general theory behind and focus on \wfs on $4$--manifolds as defined in~\cref{ch:definition wrinkled fibrations}.
Recall that these are stable maps~$f\colon X\ra B$ from a closed 4--manifold~$X$ to a closed surface~$B$ which have no definite critical points and whose fibers are connected.
The discriminant $\mc D(f)=f(\Crit(f))$ consists of finitely many cusp and double fold values, smooth arcs of simple fold values between them, and possibly also embedded circles of simple fold values.
We refer to the connected components of~$B\setminus\mc D(f)$ as the \emph{regions} of~$f$.
The different fibers over a fixed region are all diffeomorphic and our goal is to understand the relation of fibers over different regions using parallel transport.
Clearly, we can perturb an arbitrary reference arc so that it avoids the cusps and double folds and only intersects $\mc D(f)$ transversely in simple fold values.
Moreover, any such reference arc can be subdivided into parts which contain only one critical value.
We study these first.

\subsubsection{Parallel transport across folds}
	\label{ch:parallel transport across folds}
A reference arc $\g\colon[0,1]\ra B$ for a \wf $f\colon X\ra B$ is called a \emph{fold reference arc} if it meets the discriminant of $f$ in a single fold value and is transverse to the corresponding fold arc.
The vanishing sets and parallel transport along a fold reference~$\g$ are easily understood.
Indeed, it follows from the indefinite fold model that the preimage~$f\inv(\g)$ is a smooth manifold and that~$\gamma\inv\circ f\colon f\inv(\g)\ra[0,1]$ is a Morse function with a single critical point of index 1 or 2 depending on the direction in which~$\g$ crosses the fold arc.
For convenience, we write~$\S=\S_0$ and~$\S'=\S_1$ and we always assume that the index is~$2$ so that, according to \cref{eg:Takens}, the vanishing sets \wrt any connection~$\H$ are a simple closed curve $c\subset\S$ and a pair of points $\{p,q\}\subset \S'$.
In analogy with the theory of Lefschetz fibrations $c$~is usually called \emph{vanishing cycle}.
Parallel transport along~$\g$ \wrt $\H$ gives a diffeomorphism
	\begin{equation*}
	\PT_\g^\H\colon \S\setminus c \lra \S'\setminus\{p,q\}
	\end{equation*}
which can be considered as an identification of~$\S'$ as the surface obtained from~$\S$ by surgery on~$c$ as follows.%
	\footnote{So $\S'$ is in some sense ``derived'' from~$\S$, whence the notation.}
Note that the surgery of~$\S$ on~$c$ can be identified with the endpoint compactification of~$\S\setminus c$ while the endpoint compactification of~$\S'\setminus\{p,q\}$ is canonically identified with~$\S'$.
Moreover, $\PT_\g^\H$ extends to a diffeomorphism of the endpoint compactifications.
As a consequence, we see that the vanishing cycle~$c\subset\S$ must be non-separating (otherwise~$\S'$ would be disconnected, but fibers of \wfs are by definition connected)
and~$\S'$ has genus one lower than~$\S$.
Note that the specific vanishing sets and parallel transport diffeomorphisms depend on both~$\g$ and~$\H$. 
It is therefore important to understand this dependence.
\begin{lemma}\label{T:vanishing set isotopies}
Let~$f\colon X\ra B$ be a \wf.
For $s\in[0,1]$ we consider smooth families of 
	connections~$\H_s$
	and 
	fold reference arcs~$\gamma_s\colon[0,1]\ra B$ with common endpoints.
Then the vanishing sets $c_s\subset\S$ and~$\{p_s,q_s\}\subset\S'$ evolve by ambient isotopies.
Moreover, all ambient isotopies of~$c_0\subset\S$ and~$\{p_0,q_0\}\subset\S'$ can be realized by changing each~$\H_s$ in an arbitrarily small \nbhd of~$f\inv(\g_s)$.
\end{lemma}
\begin{proof}
For brevity we denote the unit interval by~$I=[0,1]$.
For each~$s\in I$ we have a smooth manifold $Y_s=f\inv(\g_s)$ and a Morse function $f_s=\g_s\inv\circ f\colon Y_s\ra I$.
We claim that the disjoint union
	$Y=\amalg_s \big(\{s\}\times Y_s\big)$ 
is a smooth submanifold of~$I\times X$.
To see this, we observe that~$Y$ can be expressed as the transverse preimage of $\amalg_s(\{s\}\times\g_s(I))$, which is a smooth submanifold of~$I\times B$, under the product map~$\id_I\times f$.
Next  we consider the map $F\colon Y\ra I^2$ given by~$F(s,y)=(f_s(y),s)$. 
Possibly after rescaling the reference arcs we can assume that $\g_s(\tfrac{1}{2})$ is the unique critical value along~$\g_s$.
It follows that the critical locus of $F$~is an arc of indefinite folds which is mapped to~$\g_I\big(\{\tfrac{1}{2}\}\big)$.
Another investigation of the fold model shows that~$F$ composed with the projection $pr_2\colon I^2\ra I$ onto the second factor is a submersion.
We can therefore integrate a suitable gradient vector field for~$pr_2\circ F$ to find a diffeomorphism~$\phi\colon I\times Y_0\ra Y$ such that~$F\circ\phi=\id \times f$.
As a consequence, we get a family of diffeomorphisms~$\phi_s\colon Y_0\ra Y_s$ satisfying $f_s\circ\phi_s=f_0$.
Now let~$g_s$ be a family of Riemannian metrics on~$X$ that induce the connections~$\H_s$.
Restricting~$g_s$ to~$Y_s$ yields a gradient for~$f_s$ (with values in~$\H_s$) which we pull back to a gradient~$v_s$ for~$f_0\colon Y_0\ra I$ via~$\phi_s$.
According to \cref{eg:Takens}, the vanishing sets of~$\g_s$ \wrt $\H_s$ are the same as the ascending and descending spheres of the unique critical points of~$f_0$ \wrt~$v_s$.
This immediately shows that the vanishing sets evolve by ambient isotopies and for the realization of arbitrary isotopies we proceed as in \cref{R:isotopies of vanishing sets}.
\end{proof}

\subsubsection{Parallel transport through cusps}
	\label{ch:parallel transport through cusps}
\begin{figure}
	\includegraphics[width=35mm]{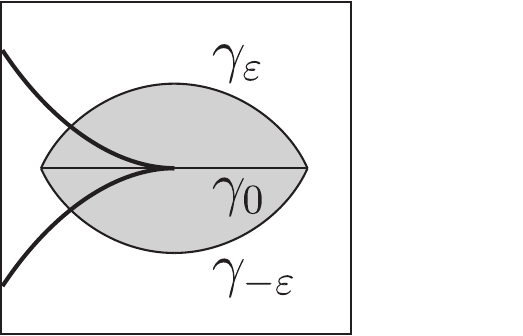}
	\caption{A cusp reference arc~$\g_0$ embedded in a family~$\g_s$ of fold reference arcs for~$s\neq0$.}
	\label{F:cusp reference arcs}
	\end{figure}

While for many purposes it is enough to consider fold reference arcs, in some situations one is forced to consider reference arcs that interact with the critical values in more complicated ways. 
For example, in \cref{ch:mapping class group interpretation} we will have work with reference arcs that contain a single cusp value at which they are tangent to the direction of the cusp -- we will refer to these as \emph{cusp reference arcs}.
Note that any cusp reference arc~$\g_0$ can be embedded in a family~$(\g_s)_{s\in[-\eps,\eps]}$ such that~$\g_s$ is a fold reference arc for~$s\neq0$ (see \cref{F:cusp reference arcs}).
The following is common folklore in this situation:
\begin{enumerate}
	\item 
		The vanishing sets of~$(\g_0,\H)$ are a single point~$p\in\S'$ and the union of two \sccs $c\cup d\subset\S$ that intersect transversely in one point.
	\item
		As $s$ approaches~$0$ from above the vanishing sets of $(\g_s,\H)$ and~$(\g_{-s},\H)$ converge to~$c,d\subset\S$ and~$r\in\S'$.
\end{enumerate}
However, we have come to realize that this might be more subtle than expected.
In fact, \cref{eg:Takens} indicates that~(1) might actually fail for arbitrary connections.
Since our arguments in \cref{ch:mapping class group interpretation} rely on~(1) we include a proof for a reasonable class of connections which are ``standard'' near the cusps.
\begin{lemma}\label{T:PT through a cusp}
The conclusions (1) and (2) above are valid if for each cusp of~$f$ there are model coordinates in which~$\H$ is induced by the standard metric on~$\R^4$.
\end{lemma}
\begin{proof}
We only give the details for~(1) and note that~(2) can be proved similarly.
By the assumption on~$\H$ it is enough to consider the indefinite fold model
	\begin{equation*}
	f\colon\R^4\ra\R^2,
	\quad
	(t,x,y,z) \mapsto (t,x^3+3tx+y^2-z^2)
	\end{equation*}
and the reference arc $\gamma_0 \colon [-1,1]\ra\R^2$ given by~$\gamma_0(\tau)=(\tau,0)$.
Our goal is to determine the vanishing sets $V_\pm$ of~$\gamma_0$ in the fibers
	\begin{equation*}
	\S_\pm=f\inv(\pm1,0)=\set{t=\pm1}\cap\set{x^3\pm 3x +y^2-z^2=0 } \subset \R^4
	\end{equation*}
\wrt the horizontal distribution~$\H_0$ induced by standard Euclidean metric on~$\R^4$.
(It is easy to see that~$\S_-$ is a once punctured torus while~$\S_+$ is diffeomorphic to the plane.)
We fix points $p_\pm=(\pm1,X_\pm,Y_\pm,Z_\pm)$ in~$\S_\pm$ and denote by
	$\tilde{\gamma}_0\colon[-1,1]\setminus\{0\}\ra \R^4$
the unique $\H_0$--lift of~$\gamma_0$ on the specified domain with~$\tilde{\gamma}_0(\pm1)=p_\pm$.
Note that~$V_\pm$ consists exactly of those~$p_\pm$ for which $\tilde{\gamma}_0(\tau)$ converges to the origin as $\tau$~goes to zero from the left and the right.
A direct computation shows that 
	\begin{equation*}
	\H_0= \left<\tfrac{\partial}{\partial t}, 3(x^2 +t)\tfrac{\partial}{\partial x} + 2y \tfrac{\partial}{\partial y} -2z \tfrac{\partial}{\partial z}  \right>
	\end{equation*}
and that $\tilde{\gamma}_0$, which we henceforth write as
	$\tilde{\gamma}_0(\tau)=\big(\tau, X(\tau), Y(\tau), Z(\tau) \big)$,
is the solution of the system of differential equations
	\begin{equation}\label{eq:cusp ODE system}
	\frac{dX}{d\tau} = -3X(X^2+\tau) \big/ \varXi,
	\hspace{5mm}
	\frac{dY}{d\tau} = -2XY \big/ \varXi,
	\hspace{5mm}
	\frac{dZ}{d\tau} = 2XZ \big/ \varXi
	\end{equation}
where $\varXi(\tau)=9(X^2+\tau)^2+4Y^2+4Z^2$.
There is one obvious solution given by~$X=Y=Z\equiv0$ which shows that~$V_\pm$ contains the point~$(\pm1,0,0,0)$.
We next show that~$V_+$ contains no other points.
This follows from an inspection of \eqref{eq:cusp ODE system} on the interval~$(0,1]$.
The first equation shows that $X\equiv 0$ if and only if~$X(\tau)=0$ for some~$\tau>0$, and in that case the other equations force~$Y$ and~$Z$ to be constant.
In particular, the only solution with~$X_+=0$ that contributes to the vanishing set is the obvious one.
On the other hand, if~$X_+\neq0$, then $X$~is nowhere zero; in fact, we have~$|X|\geq|X_+|$.
To see this, observe that if~$X_+>0$, then the \rhs of the first equation in~\eqref{eq:cusp ODE system} is strictly negative everywhere.
It follows that $X$~is monotonically decreasing and therefore bounded below by~$X_+$.
The case~$X_+<0$ is completely analogous.
As a consequence, no solution with~$X_+\neq0$ can contribute to~$V_+$ which therefore only consists of the point~$(1,0,0,0)$, as claimed.
It remains to determine~$V_-$ for which we restrict our attention to the interval~$[-1,0)$.
Again, \eqref{eq:cusp ODE system} shows that the only solution with~$X_-=0$ that contributes to~$V_-$ is the obvious one, and that if~$X_-\neq0$, then~$X$ is nowhere zero.
Arguing similarly as above we can show that if~$X_->0$, then~$|Z|\geq|Z_-|$ and if $X_-<0$, then~$|Y|\geq|Y_-|$.
To summarize, we have shown that $V_-$~is contained in the union of $\S_- \cap \set{x\geq0, \, z=0}$ and $\S_- \cap \set{x\leq0, \, y=0}$ which is easily seen to be a pair of \sccs in~$\S_-$ intersecting transversely in~$(-1,0,0,0)$.
Finally, we note that if~$p_-\in\S_- \cap \set{x\geq0, \, z=0}$, the limit of~$\tilde{\gamma}$ must lie in the intersection~$\S_0\cap \set{x\geq0, \, z=0}$ which contains only the origin.
Similar arguments apply to~$p_-\in\S_- \cap \set{x\leq0, \, y=0}$ and we conclude that 
	\begin{equation*}
	V_-=
	\big( \S_- \cap \set{x\geq0, \, z=0} \big)
	\cup 
	\big( \S_- \cap \set{x\leq0, \, y=0} \big)
	\end{equation*}
has the desired structure.
\end{proof}
\begin{remark}\label{R:adaptedness}
We do not know whether the assumption on~$\H$ in \cref{T:PT through a cusp} is necessary but it seems likely that a proof for general~$\H$ would have to involve more elaborate tools from the theory of dynamical systems.
\end{remark}

\subsubsection{Vanishing cycles in simple \wfs}
	\label{ch:vanishing cycles in swfs}
Using the techniques developed in this chapter we can prove another folklore fact about \swfs.
Let~$w\colon X\ra S^2$ be such a map.
We parameterize the equator of~$S^2$ by~$S^1$ and denote the meridian through~$\theta\in S^1$ by~$\mu_\theta$.
By a suitable reparametrization of~$S^2$ we can assume that the critical values of~$w$ are arranged near the equator such that~$\mu_\theta$ is a fold reference arc for all but finitely many values of~$\theta$ for which it is a cusp reference arc where the cusp points toward the south pole.
Let~$\theta_1,\dots,\theta_l\in S^1$ be the exceptional values, cyclically ordered according to the orientation of~$S^1$, and let~$\S$ and~$\S'$ be the fibers over the north and south pole, respectively.
\begin{lemma}\label{T:vanishing cycles in swfs}
If $\S$ has genus at least two, there is a connection for~$w$ such that the vanishing cycles of~$\mu_\theta$ are
	constant for each~$\theta\in(\theta_i,\theta_{i+1})$, say~$c_i\subset\S$, and
	the union~$c_{i-1}\cup c_i$ for~$\theta=\theta_i$.
\end{lemma}
The data~$(\S;c_1,\dots,c_l)$ is known as the \emph{surface diagram} of~$w$ and it was already observed in~\cite{Williams1} that~$\SD$ contains enough information to recover~$w$ up to equivalence (also see~\cite{Behrens}).
We will come back to these diagrams in \cref{ch:applications}.
\begin{proof}
Let~$\H$ be an arbitrary connection for~$w$ that is standard near the cusps and let~$c_i\cup d_i\subset\S$ be the vanishing set for~$\theta_i$ guaranteed by \cref{T:PT through a cusp}.
As in~\cref{R:isotopies of vanishing sets} and \cref{T:vanishing set isotopies} we can modify~$\H$ so that~$d_i=c_{i-1}$ and that the vanishing cycles for~$\theta$ near~$\theta_i$ are constantly~$c_i$ on one side and~$c_{i-1}$ on the other.
The claim now follows from \cref{T:vanishing cycles constant} below.
\end{proof}
\begin{lemma}\label{T:vanishing cycles constant}
Let $f\colon X\ra B$ be a \wf and let $(\g_s)_{s\in[0,1]}$ be a family of fold reference arcs with common endpoints such that the map~$(s,t)\mapsto\g_s(t)$ is an embedding for~$t\neq0,1$.
Let~$\H$ be a connection for~$f$ and let~$c_s\subset\S$ be the vanishing cycle of~$\g_s$.
If~$c_0=c_1$ and $\S$~has genus at least two, then~$\H$ can be modified in~$f\inv(\cup_{s\in(0,1)}\g_s)$ such that $c_s=c_0$ for all~$s$.
\end{lemma}
\begin{proof}[Proof of \cref{T:vanishing cycles constant}]
This follows from a variation of the proof of \cref{T:vanishing set isotopies} together the extra input that the space of homotopically non-trivial \sccs in a surface with negative Euler characteristic has simply connected components~(see~\cite{Ivanov}*{p.535ff}).
We can therefore find a 2--parameter family of diffeomorphisms $\phi_{s,t}\in\Diff(\S)$ such that~$\phi_{s,0}(c_0)=c_s$ and~$\phi_{s,1}(c_0)=c_0$ for all~$s$.
Finally, we can deform~$\H$ to a 1--parameter family of connections~$\H_t$ within~$f\inv(\cup_{s\in(0,1)}\g_s)$ such that the vanishing cycle of~$(\g_s,\H_t)$ is~$\phi_{s,t}(c_0)$.
\end{proof}

%% file: MCG_interpretation.tex
\section{Elimination of Cusps II: A Mapping Class Group Interpretation}
	\label{ch:mapping class group interpretation}
We can now prove our second main result which we first recall.
\restateCMvsMCG*
We begin with a more precise description of the problem.
Let $f\colon X\ra B$  be a wrinkled fibration with a pair of (indefinite) cusps~$p,q\in\Crit(f)$.
We consider an elementary cusp merge~$F=(f_s)\in\CMel(f;p,q)$ obtained by the Levine construction with a framed joining curve~$\blam_F\in\JCfr(f;p,q)$.
It will be convenient to call the composition~$\g_F=f\circ\lam_F$ the \emph{joining arc} of~$F$.
By a slight abuse of notation, we use the same notation for a slightly longer curve $\g_F\colon J_\eps\ra B$ where $J_\eps=[-1-\eps,1+\eps]$ whose ends run into the regions enclosed by the cusps of~$f$.
Since we can choose the support of~$F$ as an arbitrarily small \nbhd of the joining curve, we can assume that the fibers over the endpoints
	\begin{equation}\label{eq:reference fibers high genus}
	\S_1=f_s\inv\big(\g_F(-1-\eps)\big)
	\eqand
	\S_2=f_s\inv\big(\g_F(1+\eps)\big)
	\end{equation}
are independent of~$s$.
Observe that we can consider the extended version of~$\g_F$ as a reference arc for the various maps in the family~$(f_s)$.
We will do this for the initial map~$f=f_{-1}$ and a map~$f_\delta$ where~$\delta>0$ is small enough so that~$(f_s)$ still evolves according to the cusp merge model.
For~$f$ we find that~$\g_F$ is a concatenation of two cusp reference arcs for~$f$ in the sense of \cref{ch:parallel transport through cusps}.
So if we fix a connection~$\H$ for~$f$ which is standard near the cusps, then according to \cref{T:PT through a cusp} the vanishing sets of the nearest cusp of each endpoint of~$\g_F$ are the unions of \sccs $c_i,d_i\subset\S_i$ that intersect transversely in one point.
However, in the case of~$f_\delta$ it is easy to see from the cusp merge model that~$\g_F$ is an arc of regular values.
As a consequence, the parallel transport along~$\g_F$ \wrt any connection for~$f_F$ gives rise to a well defined isotopy class of diffeomorphisms 
	\begin{equation*}
	\varphi_F\in \pi_0\big(\Diff(\S_1,\S_2)\big).
	\end{equation*}
In order to prove \cref{T:main theorem:cusp merges and MCG} we will establish a relation between the a priori unrelated objects~$\varphi_F$, the framed joining curve~$\blam_F$, and the vanishing cycles~$c_i,d_i\subset\S_i$.
For that purpose, we also consider two fibers just outside of the cusps of~$f$
	\begin{equation}\label{eq:reference fibers low genus}
	\S'_1=f\inv\big(\g_F(-1+\eps)\big)
	\eqand
	\S'_2=f\inv\big(\g_F(1-\eps)\big)
	\end{equation}
and observe that~$\blam_F$ specifies tangent lines~$L_i\in\PP(T\S'_i)$.
However, before we can continue this discussion we need to make some general considerations.
We will return to the setup described above in \cref{ch:parallel transport before and after}.

\subsection{The difference of framed joining curves}
	\label{ch:difference of joining curves}
As above, let~$f\colon X\ra B$ be a \wf with a pair of indefinite cusps~$p,q\in\Crit(f)$.
Observe that if~$\JCfr(f;p,q)$ is non-empty, then $f$~maps the interior of the image of any joining curve from~$p$ to~$q$ into a connected component~$R\subset B\setminus \mc D(f)$.
We let~$X_R=f\inv(R)$ and write~$f_R\colon X_R\ra R$ for the restriction of~$f$ to~$X_R$.
In addition, let~$\bs X_R=\PP(\ker(df_R))$ and let~$\bs f_R\colon \bs X_R\ra R$ be the composition of~$f_R$ and the projection~$\pi\colon\bs X_R\ra X_R$.
It easily follows from the definitions that $\pi_0\big(\JCfr(f;p,q)\big)$ admits a free and transitive action of~$\pi_1(\bs X_R)$.
For concreteness, let us fix $\blam\in\JCfr(f;p,q)$ and some fiber~$\S'$ of~$f$ over the interior of the image of~$\bs f_R\circ\blam$ in~$B$.
As mentioned before, $\blam$ gives rise to a tangent line~$L\subset T_x\S'$ which we take as a base point for both~$\bs X_R$ and~$\PP(T\S')$.
Since~$R$ consists of regular values of~$f$, both~$f_R$ and $\bs f_R$ are fiber bundles over~$R$ with fibers~$\S'$ and~$\PP(T\S')$, respectively.
In particular, we can express~$\pi_1(\bs X_R;L)$ as a semi-direct product of~$\pi_1(\PP(T\S');L)$ and~$\pi_1(R;y)$ where~$y=f(x)$.
Thus, if we fix the $\pi_1(R)$~component, that is, if we only consider those~$\blam'\in\JCfr(f;p,q)$ with~$\bs f_R\circ\blam'=\bs f_R\circ\blam$, then get an element
	\begin{equation}\label{eq:difference of joining curves}
	\delta(\blam,\blam')\in\pi_1(\PP(T\S');L)
	\end{equation}
which measures the difference of the homotopy classes of~$\blam$ and~$\blam'$.
It turns out that~$\delta(\blam,\blam')$ has an interpretation in terms of mapping class groups that will be the key to the proof of \cref{T:main theorem:cusp merges and MCG}.

\subsection{A brief digression on mapping class groups}
	\label{ch:mapping class groups}
In this section we will prove some abstract results about certain mapping class groups that are closely related to our problem.
Let~$\S'$ be a closed, orientable surface of genus~$g'$ and let~$\MCG(\S')=\pi_0(\Diff^+(\S'))$ be its mapping class group.
Moreover, let~$L\subset T_x\Sigma'$ be a tangent line and let~$\Diff^+(\S';L)$ consist of all~$\phi\in\Diff^+(\S')$ such that~$d\phi(L)=L$.
We consider the group~$\MCG(\S';L)=\pi_0(\Diff^+(\S';L))$ and the ``forgetful map''
	\begin{equation*}
	F_L\colon\MCG(\S';L)\lra\MCG(\S')
	\end{equation*}
induced by the inclusion of~$\Diff(\S';L)$ in~$\Diff^+(\S')$.
We have the following analogue of the Birman exact sequence (see~\cite{Farb_Margalit}*{Ch.4.2.3}, for example).
\begin{lemma}[Generalized Birman sequence]\label{T:Birman sequence}
There is an exact sequence
	\begin{equation*}
	\cdots\ra 
	\pi_1(\Diff^+(\S'),\id) \ra 
	\pi_1(\PP(T\S'),L) \overset{\mathrm{push}}{\lra}
	\MCG(\S';L) \overset{F_L}{\lra} \MCG(\S')\ra 1.
	\end{equation*}
Moreover, $\pi_1(\Diff^+(\S'),\id)$ is trivial for~$g'\geq2$.
\end{lemma}
The map~$\mathrm{push}$ comes from a variation of the usual point pushing construction which not only drags a point around but also keeps control of a tangent direction.
\begin{proof}
We consider the map~$\mathrm{ev}_L\colon\Diff^+(\S')\ra\PP(T\S')$ defined by~$\mathrm{ev}_L(\phi)=d\phi(L)$. 
An easy modifications of the arguments in~\cite{Farb_Margalit}*{Ch.4.2.3} can be modified to show that~$\mathrm{ev}_L$ is a fiber bundle with fiber~$\Diff^+(\S';L)$.
The desired exact sequence then follows from the long exact sequence of homotopy groups and the vanishing of~$\pi_1(\Diff^+(\S'),\id)$ follows from the work of Earle and Eells~\cite{Earle-Eells}.
\end{proof}
Now let~$\S$ be another closed, orientable surface of genus~$g=g'+1$ and let~$c,d\subset\S$ be a pair of \sccs that intersect transversely in a single point.
We consider the subgroup~$\MCG(\S)(c,d)\subset\MCG(\S)$ consisting of all mapping classes that have representatives~$\phi\colon\S\ra\S$ such that~$\phi(c)=c$ and~$\phi(d)=d$.
In order to make the connection with the previous discussion, let~$\S'_{c,d}$ the surface given by the endpoint compactification of~$\S\setminus(c\cup d)$ with its distinguished endpoint~$x\in\S'_{c,d}$.
\begin{lemma}\label{T:Psi}
For any tangent line~$L\subset T_x\S'_{c,d}$ there is a canonical isomorphism
	\begin{equation*}
	\wt{\Psi}_L\colon \MCG(\S'_{c,d};L) \overset{\cong}{\lra} \MCG(\S)(c,d).
	\end{equation*}
\end{lemma}
\begin{proof}
Let~$T\subset~\S$ be a closed \nbhd of~$c\cup d$ that is diffeomorphic to a one-holed torus and let~$\S^\circ$ be the closure of~$\S\setminus T$.
Note that we can also consider~$\S^\circ$ as a subsurface of~$\S'_{c,d}$ where its complement is an open disk containing~$x$ whose closure we denote by~$D\subset\S'_{c,d}$.
The two inclusions of~$\S^\circ$ induce homomorphisms of~$\MCG(\S^\circ)$ to~$\MCG(\S)(c,d)$ and~$\MCG(\S'_{c,d};L)$.%
	\footnote{When discussing mappings class groups of surfaces with boundary we always assume that all diffeomorphisms restrict to the identity on the boundary.}
It is easy to see that the groups $\MCG(T)(c,d)$ and~$\MCG(D;L)$ are both infinite cyclic generated by~$\Delta_{c,d}=(t_c t_d)^3$ and the right-handed half twist~$h_x$ around~$x\in D$, respectively. 
Moreover, standard arguments yield short exact sequences
	\begin{equation*}
	1\ra 
	\big\langle (t_{\partial\S^\circ}\inv, h_p^2) \big\rangle \ra
	\MCG(\S^\circ)\times \MCG(D ;L) \ra 
	\MCG(\S_{c,d};L) \ra 1
	\end{equation*}
and
	\begin{equation*}
	1\ra 
	\big\langle (t_{\partial\S^\circ}\inv, \Delta_{c,d}^2) \big\rangle \ra
	\MCG(\S^\circ)\times \MCG(T)(c,d) \ra 
	\MCG(\S)(c,d) \ra 1.
	\end{equation*}
since $\Delta_{c,d}$~and~$h_p$ both square to a right-handed boundary parallel Dehn twist.
The desired isomorphism~$\Psi_L$ is therefore induced by the identity map~$\MCG(\S^\circ)$ and the isomorphism of~$\MCG(D;L)$ and~$\MCG(T)(c,d)$.
\end{proof}
By combining the maps from~\cref{T:Birman sequence,T:Psi} we obtain a homomorphism
	\begin{equation}\label{eq:Psi_L}
	\Psi_L=\wt{\Psi}_L\circ\mathrm{push}\colon \pi_1(\PP(T\S'_{c,d});L)\lra \MCG(\S)
	\end{equation}
that factors through the subgroup~$\MCG(\S)(c,d)$.
Since~$\wt{\Psi}_L$ is an isomorphism, it follows from~\cref{T:Birman sequence} that $\Psi_L$ is injective if~$\S$ has genus at least three and its kernel can easily be computed from the exact sequence.
More importantly, in the remainder of this section we will focus on the image of~$\Psi_L$.
Let~$\S'_c$ be the surfaces obtained by surgery on~$c$.
We have a subgroup $\MCG(\S)(c)$ of~$\MCG(\S)$, defined exactly as~$\MCG(\S)(c,d)$ with the single difference that only one curve is preserved, which is the source of the so called \emph{surgery homomorphism}
	\begin{equation*}
	\Phi_c\colon\MCG(\S)(c)\lra\MCG(\S'_c)
	\end{equation*}
obtained by choosing a representative of a given mapping class that fixes~$c$ set wise and extending it to the endpoint compactification.
Of course, the same discussion applies to~$d\subset\S$.
Surgery homomorphisms have appeared several times in the context of \wfs and broken Lefschetz fibrations, see~\cites{Baykur2,Behrens,HayanoR2}.
\begin{proposition}\label{T: image of Psi_l}
The image of~$\Psi_L$ is given by the subgroup
	\begin{equation}\label{eq:kernel intersection}
	\mc K(c,d)=\ker\Phi_c\cap\ker\Phi_d \subset\MCG(\S).
	\end{equation}
\end{proposition}
\begin{proof}
Let $\S^\circ$ be as in the proof of \cref{T:Psi}.
Observe that the three surfaces~$\S'_c$, $\S'_d$ and~$\S'_{c,d}$ are all obtained from~$\S^\circ$ by filling the boundary with a disk. 
Their mapping class groups are therefore canonically isomorphic to~$\MCG(\S^\circ)/\scp{t_{\del\S^\circ}}$.
Using these canonical isomorphisms we obtain a diagram
	\begin{equation*}
	\xymatrix{
	\MCG(\S)(c)
		\ar[rr]^{\Phi_c}
	&
	&
	\MCG(\S'_c)
		\ar[d]^{\cong}
	\\
	\MCG(\S)(c,d)
		\ar[u]^{\mathrm{incl}}
		\ar[d]_{\mathrm{incl}}
	&
	\MCG(\S'_{c,d};L)
		\ar[l]_{\wt{\Psi}_L}^\cong
		\ar[r]^{F_L}
	&
	\MCG(\S'_{c,d})
	\\
	\MCG(\S)(d)
		\ar[rr]^{\Phi_d}
	&
	&
	\MCG(\S'_d)
		\ar[u]_{\cong}
	}
	\end{equation*}
whose commutativity follows from the definitions.
The diagram shows that~$\wt{\Psi}_L\inv$ maps~$\mc K(c,d)$ isomorphically onto~$\ker F_L$ which agrees with $\mathrm{push}\big(\pi_1(\PP(T\S_{c,d}),L)\big)$ by \cref{T:Birman sequence}.
This proves the claim.
\end{proof}
Observe that the definition of~$\mc K(c,d)$ in~\eqref{eq:kernel intersection} makes sense for pairs of non-separating with arbitrary intersection patterns.
In the case where $c$ and~$d$ are disjoint this group was studied by the second author who obtained a generating set in~\cite{HayanoR2}*{Theorem~3.4}.
In our situation, where $c$ and~$d$ intersect in a single point, we can also obtain a generating set.
Again, let~$\S^\circ\subset\S$ be the subsurface obtained by deleting a regular \nbhd of~$c\cup d$. %
Let~$a$ be a \scc in the interior of~$\S^\circ$ and let~$\zeta\subset\S^\circ$ be an embedded arc connecting~$a$ and~$\del\S^\circ$ whose interior is disjoint from~$a$.
Then the band sum~$a\#_\zeta\del\S^\circ$ is a \scc in~$\S$ and we define
	\begin{equation}\label{eq:push lifts}
	\vartheta_{a,\zeta}:= t_a \, t_{a\#_\zeta \del\S^\circ}\inv \in\MCG(\S).
	\end{equation}
\begin{lemma}\label{T:generating the image of Psi}
For~$g(\S)\geq3$ the image of~$\tilde{\Psi}_L$ is generated by~$\Delta_{c,d}$ and the elements~$\vartheta_{a,\zeta}$ described in~\eqref{eq:push lifts}.
For~$g(\S)\leq2$ the image is generated by~$\Delta_{c,d}$ which has infinite order for~$g(\S)=2$ and order two for~$g(\S)=1$.
\end{lemma}
\begin{proof}
Note that $\pi_1(\PP(T\S'_{c,d});L)$ is generated by a loop~$\delta$ around~$\PP(T_x\S'_{c,d})\cong S^1$ and embedded loops in~$\S'_{c,d}$ based at~$x$. 
The generalized point pushing map sends~$\delta$ to the half-twist~$h_x$ or its inverse, depending on the orientation of~$\delta$.
For the loops in~$\S'_{c,d}$ we get the usual point pushing map.
Possibly after reversing~$\delta$ we get~$\Psi_L(\delta)=\Delta_{c,d}$ which is easily seen to have infinite order for~$g(\S)\geq2$ and order two for~$g(\S)=1$.
For an embedded loop~$\alpha\colon (S^1,1)\ra(\S'_{c,d},x)$ it is known that 
	\begin{equation*}
	\mathrm{push}(\alpha)= t_{a_l}\inv \, t_{a_r} \in\MCG(\S'_{c,d},x)
	\end{equation*}
where~$a_{l/r}\subset\S^\circ$ are left and right push-offs of the image of~$\alpha$ which can be chosen in the interior of~$\S^\circ$.
We can therefore also consider them as curves in~$\S$ for which we use the notation~$\tilde{a}_{l/r}$ to avoid confusion.
Now, by definition we have 
	\begin{equation*}
	\wt{\Psi}_L\big(t_{a_l}\inv \, t_{a_r}\big)=t_{\tilde{a}_l}\inv \, t_{\tilde{a}_r}\in\MCG(\S).
	\end{equation*}
and it is easy to see that $\tilde{a}_r$ can be obtained as a band sum of~$\tilde{a}_l$ and~$\del\S^\circ$ so that~$\Psi_L(\alpha)$ has the form~\eqref{eq:push lifts}.
It remains to show that~$t_{\tilde{a}_l}\inv \, t_{\tilde{a}_r}$ is the trivial mapping class when~$g(\S)\leq2$.
But this is trivial for~$g(\S)=1$ and for~$g(\S)=2$ it is easy to see that~$\tilde{a}_l$ and~$\tilde{a}_r$ are isotopic in~$\S^\circ$.
\end{proof}

\subsection{Parallel transport before and after merging cusps}
	\label{ch:parallel transport before and after}
We now return to the setting described in the beginning of \cref{ch:mapping class group interpretation} and prove \cref{T:main theorem:cusp merges and MCG}.
As before let~$f\colon X\ra B$ be a \wf with two cusps~$p,q\in\Crit(f)$ and let~$R\subset B\setminus\mc D(f)$ be the region containing the images of joining curves for~$p$ and~$q$.
Given a fixed joining arc~$\g\colon J\ra B$ between~$f(p)$ and~$f(q)$ we denote by~$\CM_0^{[\g]}(f;p,q)$ the set of elementary cusp merges whose joining is homotopic in~$R$ to~$\g$  relative to its endpoints.
Combining the considerations in \cref{ch:difference of joining curves} with \cref{T:main theorem:cusp merges and joining curves}, we see that $\pi_0\big(\CM_0^{[\g]}(f;p,q)\big)$ has a free and transitive action of~$\pi_1(\PP(T\S'_i);L_i)$ with~$\S'_i$ and~$L_i$ as in~\eqref{eq:reference fibers low genus}.
Moreover, \cref{ch:mapping class groups} already indicates how to connect the groups~$\pi_1(\PP(T\S'_i);L_i)$ and~$\mc K(c_i,d_i)$.
The precise relation will be established by carefully choosing connections and investigating the resulting parallel transport diffeomorphisms.
\begin{proof}[Proof of \cref{T:main theorem:cusp merges and MCG}]
We consider two elements $F,G$ of $\CM_0^{[\g]}(f;p,q)$.
Up to deformations within~$\CM_0^{[\g]}(f;p,q)$ we can make a series of assumptions on~$F$ and~$G$.
Later it will become clear that they are irrelevant.
First, we assume that the joining arcs~$\g_F$ and~$\g_G$ are not only homotopic to but actually agree with the fixed arc~$\g$.
Second, by \cref{T:main theorem:cusp merges and joining curves} we can assume that both homotopies are obtained from the Levine construction applied to framed joining curves~$\blam_F,\blam_G\in\JCfr(f;p,q)$ that agree near~$p$ and~$q$.
We can therefore choose the same extension of~$\g\colon J_\eps\ra B$ and, since we can control the supports of~$F$ and~$G$ simultaneously, we can consider the same fibers~$\S_i$ and~$\S'_i$ as in \eqref{eq:reference fibers high genus} and~\eqref{eq:reference fibers low genus} such that~$\blam_F$ and~$\blam_G$ span the same lines~$L_i\subset T_{x_i}\S'_i$.
Next we choose connections for the maps~$f_s$ and~$g_s$ in the homotopies as follows.
The proof of~\cref{T:Levine with parameters} shows that we can apply the Levine construction such that~$F$ and~$G$ are obtained by implanting the same truncated version~$(\tilde{\mu}_s)$ of the merge model~$(\mu_s)$ into~$f$, and such that we obtain merge coordinates for both homotopies that agree near~$p$ and possibly differ near~$q$ by the involution~$(t,x,y,z)\mapsto(t,x,-y,-z)$.
Let~$U_F$ and~$U_G$ be the neighborhoods of the joining curves where~$F$ and~$G$ agree with~$(\tilde{\mu}_s)$ which we can assume to have the same intersections with the preimages of the interior regions of the cusps of~$f$.
Moreover, by shortening~$\g$ we can assume that~$U_F\cap\S_i=U_G\cap \S_i$ is non-empty.
The Euclidean connections for the truncated model maps~$\tilde{\mu}_s$ give rise to connections for the maps~$f_s$ and~$g_s$ near the images of~$\lam_F$ and~$\lam_G$ at least for those~$s$ where~$F$ and~$G$ evolve according to the cusp merge model.
We extend these connections arbitrarily to the rest of~$X$.
For our purposes, we only need the resulting connections~$\H'_F$ and~$\H'_G$ for the initial map~$f$ as well as~$\H_F$ and~$\H_G$ for maps~$f_\delta$ and~$g_\delta$ where~$\delta>0$ is small enough.
Observe that~$\H'_F$ and~$\H'_G$ agree near~$p$ and~$q$ so that both give rise to the same vanishing cycles~$c_i,d_i\subset\S_i$, assuming that the connections are standard near the cusps so that \cref{T:PT through a cusp} applies.
Clearly, this can be verified in the model.
One readily checks that for~$s<0$ the maps~$(t,x,y,t)\mapsto(\pm(s-t^2),x,y,z)$ and~$(u,v)\mapsto(\pm(s-t^2),v)$ take the origin of the cusp model to the cusps of~$\mu_s$ and also map the Euclidean connections of both maps into each other.
Using this setup we can now start the actual proof of~\cref{T:main theorem:cusp merges and MCG}.
Since the situation is symmetric in~$\S_1$ and~$\S_2$, we will only consider the case of~$\S_1$ .
By taking parallel transport along~$\g$ we obtain diffeomorphisms
	\begin{equation*}
	\wt{\varphi}'(F,G)=
	\big(\PT^{\H'_G}_\g \big)\inv\circ\PT^{\H'_F}_\g
	\in\Diff(\S'_1)
	\end{equation*}
	\begin{equation*}
	\wt{\varphi}(F,G)=
	\big(\PT^{\H_G}_\g \big)\inv\circ\PT^{\H_F}_\g
	\in\Diff(\S_1).
	\end{equation*}
We are ultimately interested in the mapping class~$\varphi(F,G)\in\MCG(\S_1)$ represented by~$\wt{\varphi}(F,G)$.
However, we first show that~$\wt{\varphi}'(F,G)$ preserves the tangent line $L_1\in\PP(T\S'_1)$ and represents a mapping class
	\begin{equation}\label{eq:parallel transport before}
	\varphi'(F,G)=\mathrm{push}\big(\delta(\blam_F,\blam_G)\big)\in\MCG(\S'_1;L_1)
	\end{equation}
with~$\delta(\blam_F,\blam_G)\in\pi_1(\PP(T\S'_1);L_1)$ as in \cref{ch:difference of joining curves}.
To see this we trivialize~$f$ over~$\g(J'_\eps)$ where $J'_\eps=[-1+\eps,1-\eps]$ using parallel transport \wrt~$\H'_F$ so that we can identify~$f$ over~$\g(J'_\eps)$ with the projection~$\S'_1\times J'_\eps\ra J'_\eps$.
Note that in this trivialization $\blam_F$ appears as the constant curve~$\tau\mapsto(L_1,\tau)\in\PP(T\S'_1)\times J'_\eps$ while $\blam_G$ takes the form~$\tau\mapsto(\tau,\delta(\tau))$ where~$\delta$ is a closed loop in~$\PP(T\S'_1)$ representing~$\delta(\blam_F,\blam_G)$.
Moreover, the parallel transport between the various fibers~$\S'_1\times\{\tau\}$ give rise to an isotopy of~$\S'_1$ and it is easy to see that this isotopy realizes the mapping class~$\mathrm{push}\big(\delta(\blam_F,\blam_G)\big)$.
Next we show that
	\begin{equation}\label{eq:parallel transport after}
	\varphi(F,G)=\wt{\Psi}_{L_1}\big(\varphi'(F,G)\big)\in\MCG(\S_1)
	\end{equation}
with~$\wt{\Psi}_{L_1}$ as in \cref{T:Psi}.
Observe that parallel transport \wrt either~$\H'_F$ or~$\H'_G$ provides diffeomorphisms between~$\S_i\setminus(c_i\cup d_i)$ and~$\S'_i\setminus\{x_i\}$.
(Since the connections agree near the cusps, they give the same diffeomorphisms.)
Now we take closed disks~$D_i\subset\S'_i$ containing the base points~$x_i$ of~$L_i$ in their interior.
Parallel transport \wrt~$\H'_F$ or~$\H'_G$ maps them to one~holed tori~$T_i\subset\S_i$ containing~$c_i\cup d_i$.
Note that we can assume that the~$D_i$ are contained in the part of~$\S_i'$ where~$F$ and~$G$ evolve according to the cusp merge model~$(\mu_s)$. 
Next we observe that for~$t_0>0$ the Euclidean parallel transport for~$\mu_s$ maps  $(-t_0,x,y,z)\in\mu_s\inv(-t_0,0)$ to $(t_0,x,y,z)\in\mu_s\inv(t_0,0)$ provided that~$(t_0,0)$ is a regular value and the points are not contained in the vanishing sets.
From this we can deduce several things.
First, it follows that we can assume that~$D_1$ and~$D_2$ are mapped into each other by the parallel transport \wrt both~$\H'_F$ and~$\H'_G$.
Second, with this choice~$D_1$ and~$T_1$ are preserved by~$\wt{\varphi}'(F,G)$ and~$\wt{\varphi}(F,G)$, respectively.
Third, $\wt{\varphi}'(F,G)$ acts on~$D_1$ either as the identity or as a 180~degree rotation (depending on whether the merge coordinates for~$F$ and~$G$ agree near~$q$).
Moreover, in the first case $\wt{\varphi}(F,G)$ acts on~$T_1$ as the identity, and in the second as a hyperelliptic involution which has the same effect on~$c_1$ and~$d_1$ as the \mbox{$\Delta$--twist}~$\Delta_{c_1,d_1}$.
Comparing this with the definition of the map~$\wt{\Psi}_{L_1}$ we obtain~\eqref{eq:parallel transport after}.
Finally, observe that \cref{T:main theorem:cusp merges and MCG} is a consequence of~\eqref{eq:parallel transport before},~\eqref{eq:parallel transport after}, and \cref{T: image of Psi_l} and the fact that~$\varphi(F,G)$ and the isotopy classes of~$c_i$ and~$d_i$ are independent of the various choices made in the above constructions.
\end{proof}

\subsection{The case of fold merges}
	\label{ch:merging folds}
After our extensive study of cusp merges we now take a brief look at fold merges. 
As it turns out, the parallel transport problem is much simpler in this case.
Let~$f\colon X\ra B$ be a \wf and assume that a square in~$B$ intersects the critical image in two parallel fold arcs. 
We choose a connection~$\H$ for~$f$, a reference fiber~$\S$ that maps into the middle region of the square and one fold reference arc for each fold arc, say~$\g_c$ and~$\g_d$ with corresponding vanishing cycles~$c,d\subset\S$. 
In order to perform a fold merge it is necessary and sufficient that $c$ and~$d$ intersect transversely in a single point for some choice of~$\H$.
We will call this structure a \emph{fold merge configuration} and refer to~$\S$ as the \emph{inner fiber}.
In this situation, it is convenient to label the \emph{outer fibers}, that is, the fibers over the other endpoints of~$\g_c$ and~$\g_d$, by~$\S'_c$ and~$\S'_d$, respectively, since the parallel transport provides an identification with the surfaces obtained by surgery on~$c$ and~$d$.
We choose a regular \nbhd~$T$ of~$c\cup d$ as in the proof of \cref{T:Psi}.
The parallel transport along~$\g_c$ then maps~$T\setminus c$ onto a twice punctured disk~$D_c\subset\S_c$ in which~$d$ appears as an arc~$d'$ connecting the two punctures  which we can think of as marked points. 
The analogous statement holds with the roles of~$c$ and~$d$ interchanged.
The crucial observation is that parallel transport yields a diffeomorphism
	\begin{equation*}
	(\PT^\H_{\g_c})\inv \circ \PT^\H_{\g_d} \colon \S_c\setminus D_c \lra \S_d\setminus D_d
	\end{equation*}
between the complements of these disks and up to isotopy there is a unique extension to a diffeomorphism~$\S_c\ra\S_d$.
Moreover, one can show that the parallel transport \emph{after} performing a fold merge will be in this isotopy class.
So up to isotopy, there is a unique identification of the outer fibers that can be obtained after merging folds which can already be described \emph{before} actually merging folds. 

%% file: applications.tex
\section{Applications to surface diagrams of 4--manifolds}
	\label{ch:applications}
In this last section shift our attention to \swfs as defined in \cref{ch:definition wrinkled fibrations}.
As indicated in the introduction, these maps can be used to obtain a combinatorial descriptions of closed 4--manifolds reminiscent of a Heegaard diagram.
After reviewing the structure of these diagram, we will discuss how they change are affected by Williams's basic homotopies that involve cusp merges (\cref{ch:surface diagram moves}), and we will give new examples of surface diagrams of $4$--manifolds derived from Lefschetz fibrations and Heegaard diagrams (\cref{ch:Lefschetz and Heegaard}).

\subsection{Surface diagrams: a review}
	\label{ch:surface diagrams}
Let~$w\colon X\ra S^2$ be a \swf.
As explained in \cref{ch:vanishing cycles in swfs}, the set of regular values has two connected components which we refer to as the \emph{higher} and \emph{lower genus regions} depending on the fiber genus which differs by one.
We choose one reference point in each region and connect them by fold reference arcs with pairwise disjoint interiors, one passing through each arc of fold values.
Note that the orientation of the higher genus region (induced from the standard orientation of~$S^2$) induces a cyclic order of the arcs of fold values.
After choosing a connection we obtain \sccs $c_1,\dots,c_l\subset\S$ as the vanishing cycles in the higher genus fiber.
Following Williams~\cite{WilliamsPNAS}, the data $\SD_w=(\S;c_1,\dots,c_l)$ is called the \emph{surface diagram} of~$w$.
It follows easily from the results of \cref{ch:parallel transport} that~$\SD_w$ is well defined up to isotopy%
	\footnote{For~$\S$ this means that there are canonical isotopy class of diffeomorphisms between any two higher genus reference fibers since the higher genus region is simply connected.}
(cf.~\cite{Behrens}).
Moreover, using \cref{T:vanishing set isotopies} we can and will assume that the vanishing cycles are transverse and realize the minimal number of intersections within their isotopy classes.
Note that the curves in~$\SD_w$ are not arbitrary.
They satisfy the \emph{intersection condition} that any two consecutive curves intersect in one point and the \emph{monodromy condition}
\begin{equation}\label{eq:monodromy condition SD}
\Phi_{c_1}(t_{t_{c_l}(c_1)}\circ t_{t_{c_{l-1}(c_l)}}\circ \cdots \circ t_{t_{c_1}(c_2)}) = 1. 
\end{equation}
It is shown in~\cite{Behrens}*{Section~4} that if~$\S$ has genus at least three, then any collection of curves satisfying the above properties determines a \swf up to left-right equivalence.
In particular, if the fiber genus is sufficiently high, then $w$ can be recovered from~$\SD_w$ as suggested by Williams~\cite{Williams1}*{p.1054}.
\begin{remark}\label{R:abstract SD theory}
The second author developed an abstract theory of surface diagrams that allows to discuss \swfs over more general base surface (cf.~\cite{Behrens}).
In particular, it turns out that intersection condition is enough to build a \swf over the disk such that the left hand side of~\eqref{eq:monodromy condition SD} describes the monodromy of the boundary fibration.
Those familiar with Lefschetz fibrations will probably recognize a familiar pattern here.
\end{remark}

A surface diagram depicts a \swf as seen from the higher genus perspective.
Sometimes it is also useful to take a look from the lower genus side.
For that purpose Williams introduced \emph{surgered surface diagrams}~\cite{Williams2}.
As the name suggests, they are the result of performing surgery on some vanishing cycle in a surface diagram. 
Instead of giving a formal definition we will only explain their structure.
Let~$\SD=(\S;c_1,\dots,c_l)$ be a surface diagram.
We parameterize a \nbhd~$\nu c_1$ by~$S^1\times(-1,1)$ such that the non-empty intersections~$c_i\cap\nu c_1$ correspond to lines of the form~$\{\theta\}\times(-1,1)$ 
and denote by~$\S'$ the surgered surface obtained by filling the two boundary components of~$\S\setminus\nu c_1$ with disks~$D_1$ and~$D_2$ to which we refer as the \emph{surgery disks}.
Since~$\S$ can be reconstructed from~$\S'$, the surgery disks and an identification of their boundaries, we can try to capture the whole surface diagram in~$\S'$.
In order to do so we have to understand how the vanishing cycles appear in~$\S'$.
We denote by~$c_i'=c_i\setminus\nu c_1$ the remains of~$c_i$ after the surgery.
Clearly, if~$c_i$ is disjoint from~$c_1$, then~$c_i'$ is a \scc in~$\S'$ which is disjoint from the surgery disks. 
However, if~$c_i$ intersects~$c_1$, then~$c_i'$ is a collection of embedded arcs in~$\S'$ whose 
endpoints lie on the boundaries of the surgery disks while the interiors are disjoint from the disks.
Each surgery disk contains the same number of endpoints and the endpoints on the different disks are matched via the identification of the boundaries of the disks.
Finally, the curve~$c_1$ itself appears in~$\S'$ as the boundary of either disk which we also denote by~$c_1'$.
The collection~$\SD'=(\S';c_1',\dots,c_l')$ is then a surgered surface diagram.
Note that this way of picturing surfaces and curves therein is commonly used for Heegaard diagrams.
Of course, we could have performed surgery on another vanishing cycle and the resulting diagram might look rather different. 
If necessary we will use the more precise notation~$\SD'_{c_i}$ to indicate on which vanishing cycle we performed surgery, otherwise~$\SD'$ will always stand for~$\SD'_{c_1}$. 

\smallskip
When we draw surgered surface diagrams the identification of the boundaries of the surgery disks is implicitly encoded as follows. 
Up to isotopy such an identification is determined by specifying one point on each boundary. 
To obtain such a pair of points we put an adjacent vanishing cycle of the one on which we performed surgery in minimal position so that it will appear as a single arc in the surgered diagram and we take its endpoints.
Hence, as soon as we draw such an arc, we know how we have to identify the boundaries of the disks.

%% file: applications_moves.tex
\subsection{Homotopy moves for surface diagrams}
	\label{ch:surface diagram moves}
In this section we discuss how the surface diagrams of \swfs change in those of Williams's basic homotopies introduced in~\cite{Williams2} that involve cusp merges.
For brevity of notation, given a collection of non-separating \sccs~$c_1,\dots,c_k\subset\S$ we write~$\mc K(c_i)=\ker\Phi_{c_i}$ and~$\mc K(c_1,\dots,c_1)=\cap_i\mc K(c_i)$.
As explained in \cite{Williams2}, \emph{a multislide deformation} is a fold merge, resulting in an additional pair of cusps and a canonical merge arc, followed by a cusp merge along the canonical merge arc.
As an immediate consequence of \cref{T:main theorem:cusp merges and MCG} we obtain the following.
\begin{proposition}[Multislides]\label{T:multislides}
Let $w: X\rightarrow S^2$ be a simple wrinkled fibration with surface diagram $\SD=(\S;c_1,\dots,c_l)$ such that for some $1<k<l$ the curves $c_1$~and~$c_k$ intersect transversely in a single point. 
\begin{enumerate}[(i)]
	\item
		If $w^\prime:  X\rightarrow S^2$ is obtained from $w$ by a multislide deformation involving the two fold arcs with vanishing cycles $c_1$ and $c_k$, then there exists an element~$\phi\in\mc K(c_1,c_k)$ such that the surface diagram of $w^\prime$ is given by
			\begin{equation}\label{Eq:SD after multislide}
			\big( \S;c_1,\dots,c_k,\phi(c_{k+1}),\dots,\phi(c_l) \big).
			\end{equation}
	\item 
		For any $\phi\in\mc K(c_1,c_k)$, there exists a merge deformation from $w$ to a simple wrinkled fibration whose surface diagram is given by \eqref{Eq:SD after multislide}.
\end{enumerate}
\end{proposition}

\begin{figure}
\begin{center}
\includegraphics[width=110mm]{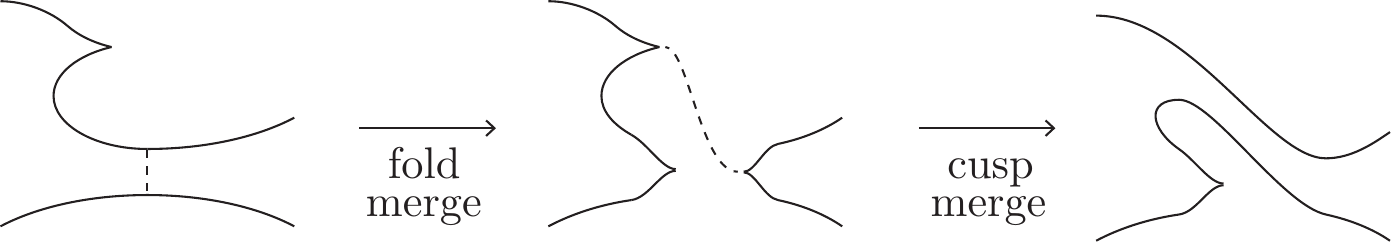}
\end{center}
\caption{The critical values throughout a shift deformation. 
}
\label{shiftdeformation}
\end{figure}

\begin{figure}[htbp]
\begin{center}
\includegraphics[width=80mm]{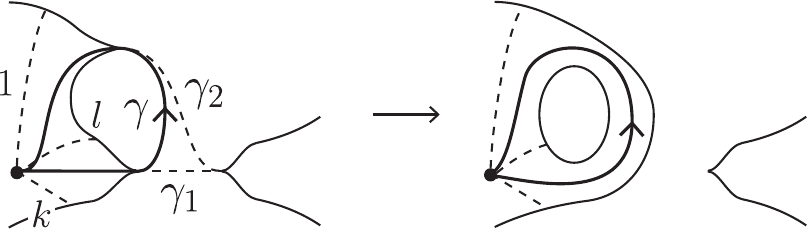}
\end{center}
\caption{Base diagrams of wrinkled fibrations. }
\label{Fig_another_cuspmerge}
\end{figure}

A slightly different combination of a fold merge followed by a cusp merge are the so called \emph{(generalized) shift deformations}.
The difference is that one of the cusps resulting from an initial fold merge is merged with a cusp adjacent to the other along a canonical merge arc, see \cref{shiftdeformation}. 
Observe that a surface diagram of the resulting fibration is uniquely determined once we specify a framed joining curve mapping to the canonical joining arc.
It turns out that there are two distinguished joining curves and we shall reserve the name \emph{shift deformation} for generalized shifts derived from these curves.
We will say more about this momentarily after proving the analogue of \cref{T:multislides} for generalized shifts.
\begin{proposition}[Generalized shifts]\label{T:shifts}
Let $w: X\rightarrow S^2$ be a simple wrinkled fibration with surface diagram $\SD=(\S;c_1,\dots,c_l)$ such that for some $1<k<l$ the curves $c_1$~and~$c_k$ intersect transversely in a single point. 
\begin{enumerate}[(i)] 
	\item 
		If $w':  X\rightarrow S^2$ is obtained from $w$ by a generalized shift deformation whose initial fold merge is applied to the fold arcs with vanishing cycles $c_k$ and $c_l$. 
		Then the surface diagram of $w'$ is given by
			\begin{equation}\label{Eq:SD after shift}
			\big( \S ; c_1,\dots,c_k,c_l,\chi(c_{k+1}),\dots,\chi(c_{l-1}) \big),
			\end{equation}
		where~$\chi\in \MCG(\S)$ satisfies 
			$\chi t_{c_l}^{-1}t_{c_k}^{-1} \in \mc K(c_l)$ and 
			$\chi t_{c_1}^{-1} t_{c_l}^{-1} \in \mc K(c_1)$. 
	\item 
		For any $\chi \in \MCG(\S)$ satisfying the conditions in~(i) there exists a generalized shift deformation from $w$ to a simple wrinkled fibration whose surface diagram is given by \eqref{Eq:SD after shift}. 
\end{enumerate}
\end{proposition}
\begin{proof}
\cref{shift_pointpath_forproof} describes a neighborhood of the fold arcs in~$S^2$ corresponding to the vanishing cycles $c_1$, $c_k$ and $c_l$ decorated with various reference points and oriented reference arcs between them.
%
%
Let~$\mathcal{H}$ be a connection for~$w$ and let~$\chi\in\MCG(\S)$ be the isotopy class of the diffeomorphism ${\PT^{\mathcal{H}}_{\gamma_1}}^{-1}\circ \PT^{\mathcal{H}}_{\gamma_0}: \S\rightarrow \S$. 
Clearly, the surface diagram of~$w^\prime$ is given by
	$\big( \S ; c_1,\dots,c_k,c_l,\chi(c_{k+1}),\dots,\chi(c_{l-1}) \big)$. 
Observe that the composition $\chi\, t_{c_l}^{-1}t_{c_k}^{-1}t_{c_l}$ preserves~$c_l$ up to isotopy and its image under~$\Phi_{c_l}$ can be considered as the monodromy along the concatenation of $\varepsilon_1, \varepsilon_2$ and $\varepsilon_3$. 
Since this concatenation is null-homotopic, we have $\chi t_{c_l}^{-1}t_{c_k}^{-1}t_{c_l}\in\mc K(c_l)$, hence also $\chi t_{c_l}^{-1}t_{c_k}^{-1}\in\mc K(c_l)$. 
Similar arguments involving the arcs $\delta_1, \delta_2$ and $\delta_3$ show that $\chi t_{c_1}^{-1}t_{c_l}^{-1}\in\mc K(c_1)$ which concludes the proof of~(i).

Now suppose that $\chi^\prime \in \MCG(\S)$ also satisfies the conditions in~(i). 
Then we have $\chi \chi^{\prime-1}\in\mc K(c_1,c_l)$ and we can use \cref{T:main theorem:cusp merges and MCG} to construct a new generalized shift that produced the surface diagram $\big( \S ; c_1,\dots,c_k,c_l,\chi^\prime(c_{k+1}),\dots,\chi^\prime(c_{l-1}) \big)$.
%
	\begin{figure}[t!]
	\begin{center}
	\includegraphics[width=100mm]{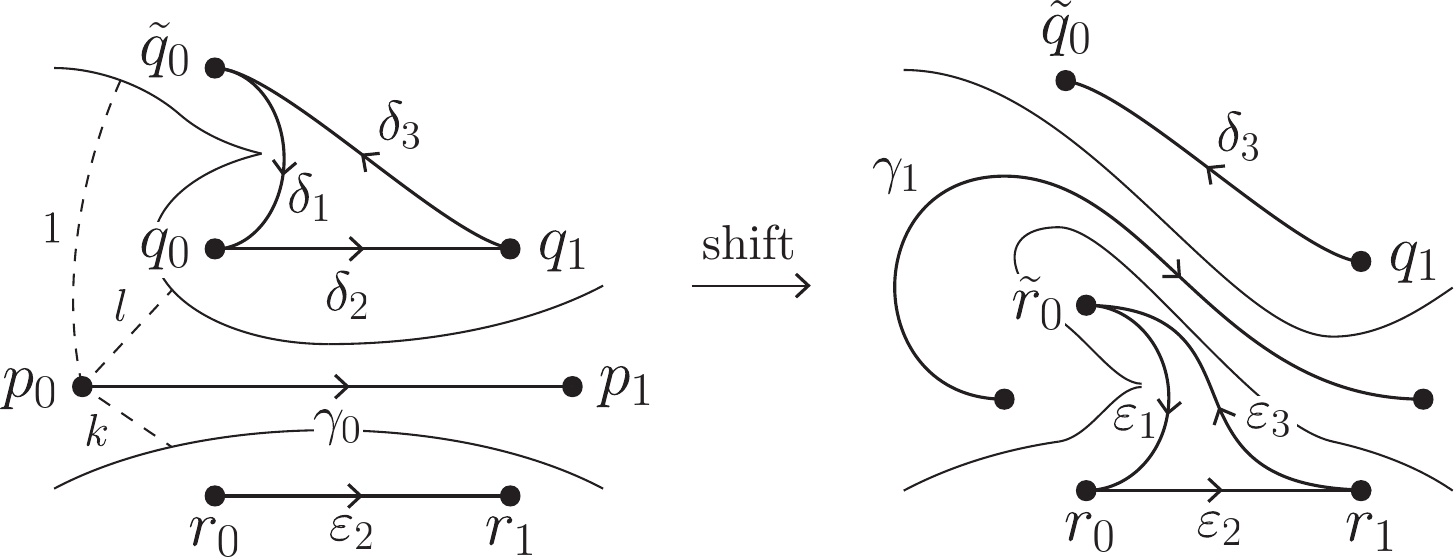}
	\end{center}
	\caption{The dashed curves describe reference paths which give vanishing cycles $c_1, c_k$ and $c_l$. }
	\label{shift_pointpath_forproof}
	\end{figure}
\end{proof}
Next we explain the actual \emph{shift deformations} alluded to earlier.
As before, let~$w:X\rightarrow S^2$ be a \swf with surface diagram $(\S; c_1,\ldots, c_l)$ such that $c_k$~and~$c_l$ intersect transversely in one point. 
Let~$w_0$ be the \wf obtained from~$w$ by a fold merge involving the fold arcs corresponding to~$c_k$ and~$c_l$.
Let~$\lam_0$ be the canonical joining curve for the two cusps created by the fold merge.
We consider a based loop~$\g$ in~$S^2$ as on the left of \cref{Fig_another_cuspmerge}.
Part of~$\g$ can be considered as a joining arc for two cusps of~$\hat{w}_0$.
We claim that there are two distinguished homotopy classes joining curves lifting~$\g$ to~$X$.
Indeed, each regular value along~$\g$ can be connected to the arc of folds to the left by a horizontal line and the vanishing set \wrt some fixed connection that is standard near the cusps is a pair of points in the corresponding fiber.
Using \cref{T:vanishing set isotopies,T:PT through a cusp} we can show that these points trace out two curves in~$X$ connecting the two cusps over~$\g$ which we can consider as joining curves, say~$\hat{\lam}_1$ and~$\hat{\lam}_2$.
Moreover, arguing as in the proofs of \cref{T:vanishing set isotopies,T:PT through a cusp} it follows that the homotopy classes of~$\hat{\lam}_1$ and~$\hat{\lam}_2$ are independent of the connection and the precise choice of horizontal lines.
By perturbing the concatenation of~$\lam_0$ with the reverse of~$\hat{\lam}_i$ we obtain two joining curves~$\lam_i$, $i=1,2$, and the cusp merge along either of them with an arbitrary framing completes a shift deformations resulting in simple \wfs $w'_i$.
Our goal is to understand the elements $\chi_i\in\MCG(\S)$ that arise from \cref{T:shifts}(i) applied to~$w_i$.
For that purpose, we consider  \wfs~$\hat{w}_i$ obtained from~$w_0$ by cusp merges along either~$\hat{\lam}_i$ with some framing.
The discriminant of~$\hat{w}_i$ has three components as shown on the right of \cref{Fig_another_cuspmerge}.
Let~$\varphi_i\in\MCG(\S)$ be the monodromy along~$\g$ \wrt~$\hat{w}_i$.
Clearly, $\varphi_i$~preserves~$c_l$ and maps~$c_k$ to~$c_1$ up to isotopy.  
Changing the framing of the joining curve if necessary, we can assume that~$\varphi_i$ preserves the orientation of $c_l$. 
We choose a reference point~$p_0\in S^2$ in the lower genus region of~$\hat{w}_i$ and a fold reference arc from~$p_0$ through the fold arc with vanishing cycle~$c_l$. 
This path gives a surgered surface diagram in the fiber~$\S'$ over~$p_0$ which we identify with the surgered surface~$\S'{c_l}$.
We can collapse the surgery disks to points and obtain a graph in~$\S'_{c_l}$ with two vertices $v_1,v_2\in \S'_{c_l}$ and edges given by the arcs~$c_k^\prime$ and~$c_1^\prime$ corresponding to~$c_k$ and~$c_1$.
These arcs form a loop which goes through both~$v_1$ and~$v_2$. 
By small perturbations we can obtain two loops~$\eta_1\subset \S'_{c_l}\setminus \{v_2\}$ and $\eta_2\subset \S'_{c_l}\setminus \{v_1\}$ where~$\eta_i$ is based at~$v_i$. 
We orient~$\eta_i$ so that it leaves~$v_i$ along $c_k^\prime$ and returns along $c_1^\prime$. 

\begin{lemma}\label{prop_monodromy_another_cuspmerge}

Let $\vec{c}_l\subset T\S$ be a nowhere vanishing section of $Tc_l$. 
After changing indices if necessary, there is an integer~$m$ such that $\varphi_i\Delta_{c_k,c_l}^{2m}\in \MCG(\S)(\vec{c}_l)$ is sent to the pushing maps along~$\eta_i$ by the map
	$\Phi_{c_l}^\ast: \MCG(\S)(\vec{c}_l)\rightarrow \MCG(\S'_{c_l};v_1,v_2)$.

\end{lemma}

\begin{proof}

The fold circle of~$\hat{w}_i$ in the source which is mapped into the inner region of~$\gamma$ bounds a disk. 
Thus $\varphi_i$~is contained in the kernel of the composition of $\Phi_{c_l}^\ast$ with the forgetful map of either~$v_1$ or~$v_2$. 
Possible after changing the labels we can achieve that $\varphi_i\in\ker(\mathrm{Forget}_{v_i}\circ \Phi_{c_l}^\ast)$.
We only give the arguments for~$\varphi_2$ since $\varphi_1$ can be dealt with similarly.
The image $\Phi_{c_l}^\ast(\varphi)$ is a pushing map along some loop~$\delta\in \pi_1(\S'_{c_l}\setminus \{v_2\}, v_1)$. 
Let~$\Pi_{v_1,v_2}$ be the set of homotopy classes of arcs connecting~$v_1$ to~$v_2$.  
The group $\pi_1(\S'_{c_l}\setminus \{v_2\}, v_1)$ acts on~$\Pi_{v_1,v_2}$ in the obvious way. 
Since~$\varphi$ maps~$c_k$ to~$c_1$, $\delta \cdot c_k^\prime$~is equal to~$c_1^\prime$. 
Thus, the element~$\eta_1^{-1}\delta$ is contained in the stabilizer of~$c_k^\prime$ under the action of~$\pi_1(\S'_{c_l}\setminus \{v_2\}, v_1)$, which is the infinite cyclic group generated by the loop obtained by connecting~$v_1$ to a sufficiently small circle around~$v_2$ using~$c_k^\prime$. 
This loop is the image of~$\Delta_{c_k, c_l}^2$ under the homomorphism~$\Phi_{c_l}^\ast$. 
\end{proof}

\begin{remark}\label{R:non-uniqueness shifts}

In general, the two pushing maps along~$\eta_1$ and~$\eta_2$ are different.
For example, it is easy to see that the two pushing maps give distinct elements of~$\MCG(\S'_{c_l};v_1,v_2)$ if $c_1$ and $c_k$ are disjoint. 

\end{remark}

Now let us return to the mapping classes $\chi_i\in \MCG(\S)$ associated via \cref{T:shifts} with the shift deformations resulting in~$w'_i$ for~$i=1,2$
%
\begin{lemma}\label{prop_MCGinterpretation_shift}

For some integer~$m$ we have 
	$\chi_i=\varphi_i \, t_{c_k}t_{c_l}t_{c_k}\Delta_{c_k,c_l}^m\in \MCG(\S)$.

\end{lemma}

\begin{proof}

We denote by~$P$ the subset of~$\MCG(\S)$ which maps~$c_k\cup c_l$ to~$c_l \cup c_1$, and by $P^\prime$ the set of isotopy classes of diffeomorphisms from $(\S'_{c_k, c_l},v)$ to $(\S'_{c_1,c_l},v^\prime)$. 
The sets $P$ and $P^\prime$ admit free and transitive actions of $\MCG(\S)(\{c_k,c_l\})$ and $\MCG(\S'_{c_k,c_l},v)$ defined by compositions as maps, respectively. 
Using similar definitions as for the surgery homomorphisms used in \cref{ch:mapping class groups} we obtain maps
	\begin{equation*}
	\Phi^\ast: P\rightarrow P^\prime 
	\eqand
	\Phi_{c_k,c_l}^\ast: \MCG(\S)(\{c_k,c_l\}) \rightarrow \MCG(\S'_{c_k,c_l};v)  
	\end{equation*}
which are compatible with the actions on~$P$ and~$P'$.%
	\footnote{that is, for $\xi\in P$ and $\varphi \in \MCG(\S)(\{c_k,c_l\})$ we have $\Phi^\ast(\xi \cdot \varphi)=\Phi^\ast(\xi)\cdot \Phi_{c_k,c_l}^\ast(\varphi)$}
%
%
Note that both~$\chi_i$ and~$\varphi_i$ are contained in~$P$.
Since the action of~$\MCG(\S)(\{c_k,c_l\})$ on~$P$ is transitive, we can find some~$\xi\in\MCG(\S)(\{c_k,c_l\})$ such that~$\xi \cdot \chi_i=\varphi_i$. 
The image~$\Phi^\ast(\chi_i)$ is the monodromy of~$\hat{w}_i$ along~$\gamma$ under a suitable connection, while~$\Phi^\ast(\varphi_i)$ is that along the curve obtained from~$\gamma$ by changing part of~$\gamma$ into the images~$\gamma_1$ and~$\gamma_2$ of joining curves (see \cref{Fig_another_cuspmerge}), in particular~$\Phi^\ast(\chi_i)=\Phi^\ast(\varphi_i)$. 
Since the action of~$\MCG(\S'_{c_k,c_l};v)$ on~$P^\prime$ is free, $\xi$~is contained in the kernel of~$\Phi_{c_k,c_l}^\ast$ which is the infinite cyclic group generated by~$t_{c_k}t_{c_l}t_{c_k}$. 
Since $\varphi_i$ preserves~$c_k$ up to isotopy while~$\chi_i$ does not, $\xi$~is equal to~$(t_{c_k}t_{c_l}t_{c_k})^{2m+1} = t_{c_k}t_{c_l}t_{c_k} \Delta_{c_k,c_l}^m$ for some integer~$m$. 
%
\end{proof}

\begin{remark}\label{R:}
According to \cref{prop_MCGinterpretation_shift}, we can obtain~$\chi_i$ once we determine~$\varphi_i$. 
The kernel of~$\Phi_{c_k,c_l}^\ast$ in \cref{prop_monodromy_another_cuspmerge} is generated by the Dehn twist along~$c_l$. 
Moreover, $\varphi_i$~maps~$c_k$ to~$c_1$. 
Thus, we can easily obtain~$\varphi_i$ from the vanishing cycles of~$w$ using \cref{prop_monodromy_another_cuspmerge}. 
\end{remark}

%
\begin{remark}[Handleslides and stabilizations]

For completeness, we would like to mention that the effect of the two remaining cases of Williams's basic homotopies, namely handleslides and stabilizations, was studied by the second author in~\cite{HayanoR2}
Indeed, the effect of stabilizations are completely described by \cite{HayanoR2}*{Theorems~6.5 and~6.7}) and the effects of handleslides can easily be understood using \cite{HayanoR2}*{Theorem~3.9}. 
\end{remark}

%% file: applications_LF_HD.tex
\subsection{Constructions of surface diagrams}
	\label{ch:Lefschetz and Heegaard}
In this section we describe some constructions of \swfs that involve cusp merges and show how to obtain surface diagrams. 
In order to do this we have to turn our theoretical understanding of connections in cusp merge into a method which is useful in practice. 

\begin{figure}[h]
	\includegraphics[width=40mm]{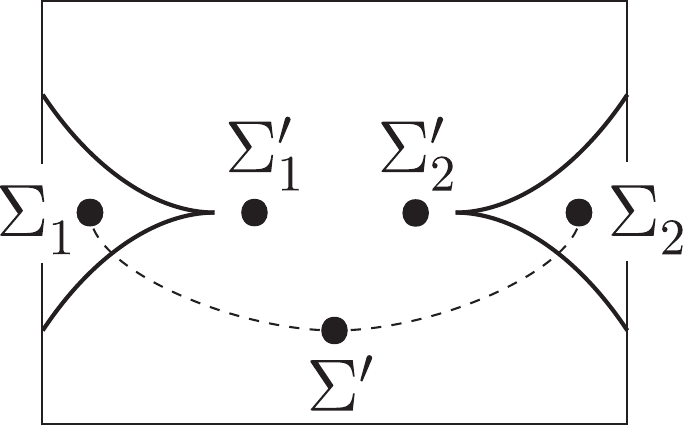}
	\caption{Matching surgered surface diagrams.}
	\label{F:matching surgered SDs}
	\end{figure}

\smallskip
Let~$f\colon X\ra B$ be a \wf that contains a cusp merge configuration as depicted in \cref{F:matching surgered SDs} and fix a connection~$\H$ that is standard near the cusps.
As indicated, we choose reference fibers~$\S_i$ and~$\S'_i$ near the cusps and another one~$\S'$ further away in the lower genus region.
We denote the vanishing cycles in~$\S_i$ by~$c_i$ and~$d_i$ where the~$c_i$ correspond to the lower fold arcs.

By parallel transport we can identify~$\S'$ with the surgery on either~$c_1\subset\S_1$ or~$c_2\subset\S_2$ and each possibility gives rise to a pair of surgery disks~$D_i^\pm\subset\S'$ and an arc~$d'_i$ connecting the disks.
Recall that with this data we can completely reconstruct the higher genus fiber~$\S_i$ from~$\S'$.
There is now an easy way to identify the higher genus fibers.
We simply have to find a \emph{matching isotopy}~$\bs\phi=(\phi_t)$ of~$\S'$ that moves the whole configuration~$D_2^-\cup d'_2\cup D_2^+$ to~$D_1^-\cup d'_1\cup D_1^+$.
Clearly, such an isotopy induces a diffeomorphism~$\tilde{\bs\phi}\colon\S_2\ra\S_1$.

To see that the identifications we obtain in this way are the same as the ones coming from cusp merges, we first observe that the choice of a matching isotopy is essentially equivalent to the choice of a path in~$\PP(T\S')$.
Indeed, it easy to see that from any path~$\bs\lam$ in~$\PP(T\S')$ that connects a tangent spaces of~$d'_2$ and~$d'_1$ at interior points, we can construct a matching isotopy which realizes~$\mathrm{push}(\bs\lam)$.
On the other hand, any matching isotopy gives rise to such a path as the trace of a tangent space to~$d'_2$.
For concreteness, let us fix tangent spaces~$T_i\in\PP(T\S')$ of~$d'_i$.
To make the connection to cusp merge homotopies, we denote the tangent lines in the lower genus fibers near the cusps by~$L_i\in\PP(\S'_i)$. 
We can assume that, after trivializing the lower genus region using~$\H$, the lines~$L_i$ correspond to the lines~$T_i$. 
As a consequence, given any path in~$\PP(T\S')$ from~$T_2$ to~$T_1$, we can construct a framed joining curve covering the obvious joining curve and it is obvious that the parallel transport obtained after performing the corresponding cusp merge agrees with the identification obtained from the matching isotopy associated to~$\bs\lam$.

We will discuss examples of this procedure in~\cref{ch:Lefschetz,ch:SD from HD} below.

\subsubsection{Surface diagrams for Lefschetz fibrations}
	\label{ch:Lefschetz}
\begin{figure}
\includegraphics[width=100mm]{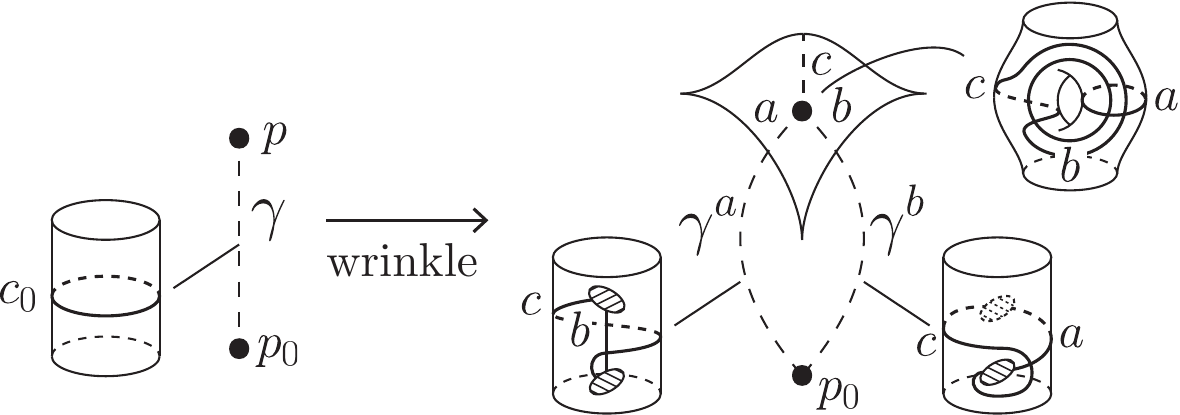}
\caption{Surgered surface diagrams in a wrinkling homotopy. }
\label{vanishingcycle_wrinkle}
\end{figure}

%
Let~$f\colon X\ra B$ be a smooth map from an oriented 4--manifold~$X$ to an oriented surface~$B$.
A critical point $p\in \mc C(f)$ is called a \emph{Lefschetz singularity} if there exist orientation preserving complex charts $\varphi:U\to \C^2$ and $\phi:V\to \C$ around $p$ and $f(p)$, respectively, such that
\begin{equation}\label{E:local model Lef}
\phi\circ f \circ \varphi^{-1}(z,w)= z^2+w^2. 
\end{equation}
Furthermore, $f$ is called a \emph{Lefschetz fibration} if all its critical points are Lefschetz singularities.
We also consider so called \emph{achiral Lefschetz singularities} which have the same model as \cref{E:local model Lef} but the chart $\phi$ does not preserve the orientation. 
A map with only Lefschetz and achiral Lefschetz singularities is called an \emph{achiral Lefschetz fibration}.
It is clear from the local models that the critical points of an achiral Lefschetz fibration are isolated.
In particular achiral Lefschetz fibrations are not stable.
But \cref{T:maps to surfaces} guarantees that we can perturb any achiral Lefschetz fibration so that it becomes a stable map. 
An explicit perturbation of the local model was given in~\cite{Lekili}:
	\begin{equation*}
	W\colon [0,\infty) \times \mathbb{C}^2 \to \mathbb{C}\colon (s,z,w) \mapsto z^2 + w^2 + s\Re(z).  
	\end{equation*}
This perturbation is  known as \emph{wrinkling} and changes a Lefschetz singularity into a circular singularity which consists of indefinite folds and three indefinite cusps (see \cref{vanishingcycle_wrinkle}). 
Observe that we can obtain a simple wrinkled fibration from a Lefschetz fibration over~$S^2$ by first wrinkling all the Lefschetz singularities and then connecting all the components of singular sets by cusp merges. 
According to~\cite{Lekili}, the vanishing cycles of the three folds in a fiber over the inner region of the triangle of critical values of the wrinkled Lefschetz singularity are related to the Lefschetz vanishing cycle as described in \cref{vanishingcycle_wrinkle}. 
Combining this diagram with the method of matching isotopies described in \cref{ch:surface diagrams} we can obtain a surface diagram of a Lefschetz fibration once we know its vanishing cycles.
We will carry out this procedure for the genus--$1$ achiral Lefschetz fibration on $S^4$ due to Matsumoto~\cite{Matsumoto_1982}.

\begin{example}\label{eg:Matsumoto fibration}
Matsumoto's achiral Lefschetz fibration $f\colon S^4\to S^2$ is constructed as follows.
Let $\pi\colon Y \to D^2$ be a singular fibration with a regular fiber $F\cong T^2$ which has a Lefschetz and an achiral Lefschetz singularities with the same non-separating vanishing cycle $c\subset F$. 
It is easy to see that the restriction $\pi|_{\partial Y}$ is trivial. 
We take a bundle isomorphism $\Phi\colon \partial D^2\times T^2 \to \partial Y$ so that the circle $\Phi(\partial D^2\times \{\ast\})$ represents a generator of $H_1(Y;\mathbb{Z})$. 
We put $X = Y \cup_{\Phi} (D^2\times T^2)$.  
It is not hard to see that $X$ is diffeomorphic to $S^4$ and admits an achiral Lefschetz fibration $f = \pi \cup \mathrm{pr} \colon X = Y\cup (D^2\times T^2)\to S^2$, where $\mathrm{pr}\colon D^2\times T^2 \to D^2$ is the projection onto the first component. 
Denote by $p_1\in S^2$ (resp.~$p_2\in S^2$) the image of the Lefschetz (resp.~achiral Lefschetz) singularity of $f$. 
After wrinkling the two singularities of~$f$ we obtain a wrinkled fibration $f_1\colon S^4\to S^2$ whose critical values and vanishing cycles in a regular fiber $f_1^{-1}(p_0)$ are described in the middle of \cref{fig_cycles_achiral_S4}. 
\begin{figure}[htbp]
\includegraphics[width=125mm]{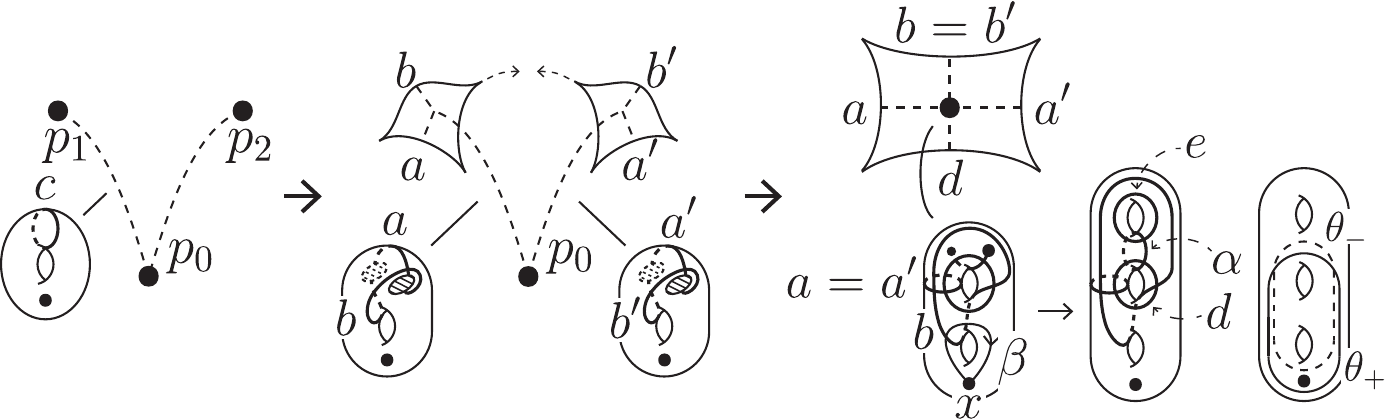}
\caption{The procedure for obtaining a genus--$2$ simple wrinkled fibration.}
\label{fig_cycles_achiral_S4}
\end{figure}
We apply a cusp merge to~$f_1$ along the dotted arrow in the figure. 
The resulting fibration~$f_2\colon S^4\to S^2$ is a simple wrinkled fibration. 
Vanishing cycles of~$f_2$ in a genus--$2$ regular fiber $\Sigma$ are described in the right of \cref{fig_cycles_achiral_S4}. 
We take a small disk $D_0\subset S^2$ in the higher-genus region. 
We can deduce from the construction of $f$ in the previous paragraph that an attaching loop of $f_2$ is $\beta\in \pi_1(\S, x)$ described in the right of \cref{fig_cycles_achiral_S4}. 
We replace a neighborhood of $d$ in $\Sigma$ with a once punctured torus as shown in \cref{fig_cycles_achiral_S4} and denote the resulting genus--$3$ surface by $\tilde{\Sigma}$. 
We take simple closed curves $\theta_+$ and $\theta_-$ in $\Sigma_3$ as shown in \cref{fig_cycles_achiral_S4}. 
Denote the mapping class $t_{t_a(d)}t_{t_b(a)}t_{t_{a}(b)}t_{t_{d}(a)}\mathrm{push}(\beta)t_{\theta_+}t_{\theta_-}^{-1}\in \MCG(\tilde{\Sigma};x)$ by $\theta$. 
The mapping class $\theta$ is contained in the group $\MCG(\Sigma;x)(d,e)$. 
Furthermore, we can easily verify that $\theta$ satisfies the condition $W_3'$ in \cite{HayanoR2}*{\parasign6}. 
Thus, by \cite{HayanoR2}*{Theorem~6.7} a surface diagram of a fibration obtained by stabilizing $f_2$ is as follows: 
\[
(\tilde{\Sigma};a,b,a,d,\tilde{\theta}(\alpha),e,\alpha,d), 
\]
where $\tilde{\theta}$ denotes the mapping class $t_{t_a(d)}t_{t_b(a)}t_{t_{a}(b)}t_{t_{d}(a)}\theta^{-1} = t_{\theta_-}t_{\theta_+}^{-1}\mathrm{push}(\beta)^{-1}$. 
The simple closed curves in this surface diagram are shown in the upper half of \cref{fig_SD_S4}. 
\begin{figure}[htbp]
\includegraphics[width=90mm]{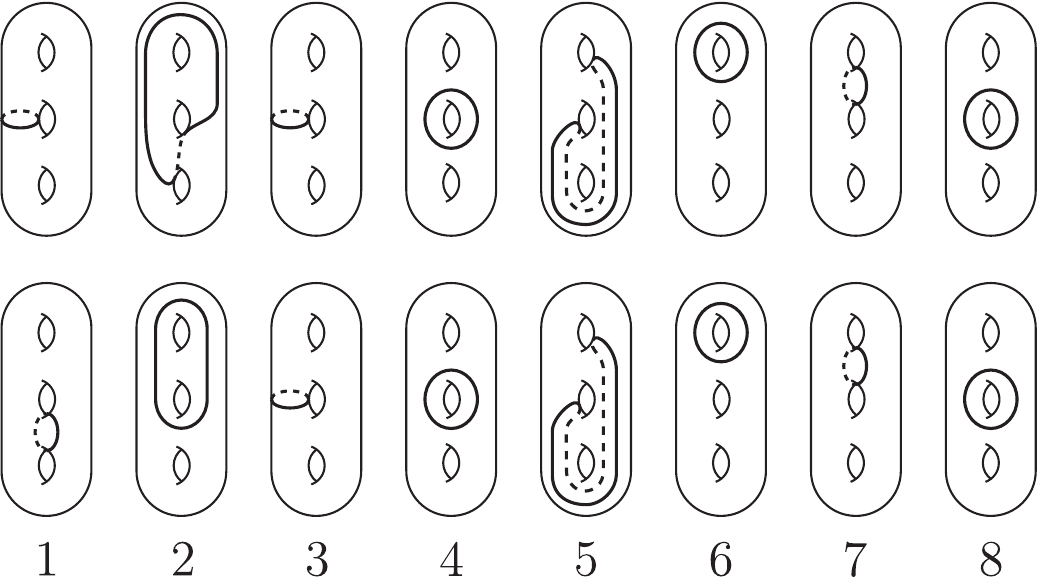}
\caption{Surface diagrams corresponding to the two elements in $\pi_4(S^2)\cong \mathbb{Z}/2\mathbb{Z}$: the upper (resp.~lower) one represents the element $1\in\pi_4(S^2)$ (resp.~$0\in\pi_4(S^2)$). }
\label{fig_SD_S4}
\end{figure}

\end{example}

\begin{remark}\label{R:two SDs on S4}

Note that Matsumoto's fibration represents the non-trivial element of~$\pi_4(S^2)\cong\mathbb{Z}/2\mathbb{Z}$. 
Indeed, it can be obtained as a perturbation of the map $h\circ Sh\colon S^4\to S^2$, where~$h\colon S^3\to S^2$ is the Hopf fibration and~$Sh\colon S^4\to S^3$ is its suspension, which is known to generate~$\pi_4(S^2)$ (see \cite{Matsumoto_1982}*{\parasign3}). 
On the other hand, we can also realize the constant homotopy class as follows. 
One readily checks that the projection~$\R^5\ra\R^2$ restricts to a stable map~$\pi\colon S^4\ra S^2$ whose critical locus is a single circle of definite folds mapped injectively into the plane.
Using a homotopy described in~\cite{Williams1}*{Example~2} we can trade the definite folds for indefinite ones. 
Applying stabilization twice to the resulting map, we eventually obtain a genus--$3$ \swf.
It is not hard to see that the surface diagram of this map is give by lower half of \cref{fig_SD_S4}.
It is interesting to observe that, although both diagrams in \cref{fig_SD_S4} describe the 4--manifold~$S^4$, they cannot be related by the moves derived from Williams's basic homotopies since the corresponding simple wrinkled fibrations are not homotopic. 
Yet the diagrams share a striking resemblance. 
This indicates that there might be a reasonable diagrammatic way to relate surface diagrams derived from non-homotopic maps which is still an open problem at the time of writing.
\end{remark}

\subsubsection{Surface diagrams from Heegaard diagrams}
	\label{ch:SD from HD}
\begin{figure}
\centering
\subfigure[The map $f_0$.]{\includegraphics[height=36mm]{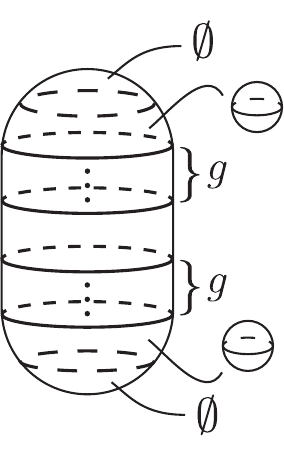}
\label{second_deformation1}}
\hspace{.2em}
\subfigure[The map $f_1$.]{\includegraphics[height=36mm]{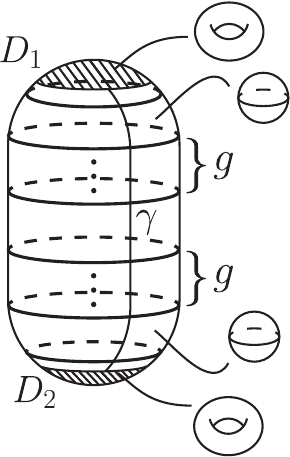}
\label{second_deformation2}}
\hspace{.2em}
\subfigure[]{\includegraphics[height=36mm]{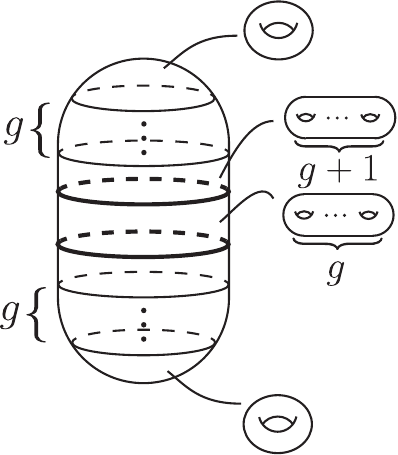}
\label{second_deformation3}}
\hspace{.2em}
\subfigure[The map $f_2$.]{\includegraphics[height=36mm]{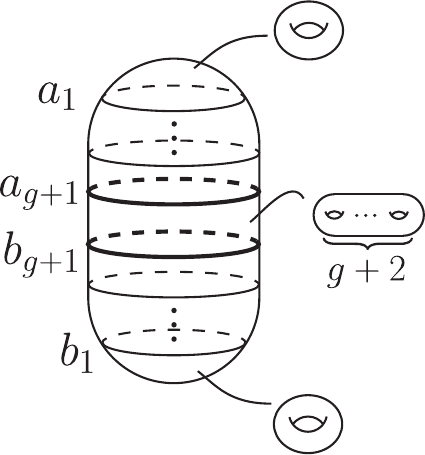}
\label{second_deformation4}}
\caption{The critical values and fibers in various stages of the algorithm.}
\end{figure}

Let~$M$ be a connected, closed, and oriented \mbox{3--manifold}.
In this section we will explain an algorithm to obtain a surface diagram of~$S^1\times M$ from a Heegaard diagram of~$M$.
It will follow from the construction that the corresponding simple wrinkled fibration is not surjective and is therefore homotopic to a constant map. 
Recall that a Heegaard diagram of~$M$ is a triple~$(\Sigma_g; \bs\alpha, \bs\beta)$, where 
	$\Sigma_g$ is a closed, oriented surface of genus~$g$, and 
	$\bs\alpha = \{\alpha_1,\ldots, \alpha_g\}$ and $\bs\beta=\{\beta_1,\ldots, \beta_g\}$ are $g$-tuples of mutually disjoint simple closed curves on $\Sigma_g$ both of which are linearly independent in~$H_1(\Sigma_g)$. 
Any Heegaard diagram of~$M$ can be obtained from a Morse function~$h: M\rightarrow \R$ with unique minima and maxima mapping to~$-g-1$ and~$g+1$, respectively, and the Heegaard surface at level~$0$ such the index~1 critical point with ascending sphere~$\alpha_i$ maps to~$-g-1+i$ and the index~2 critical point with descending sphere~$\beta_i$ maps to~$i$.
Using such a Morse function we obtain a stable map 
	\begin{equation*}
	\id \times h : S^1\times M \rightarrow S^1\times [-g-2,g+2]
	\end{equation*}
which we can promote to a stable map~$f_0: S^1\times M\rightarrow S^2$ by collapsing the boundary of the annulus~$S^1\times [-g-2, g+2]$.
Note that the critical set of~$f_0$ consists of $2g+2$ circles of folds that are mapped injectively into~$S^2$ as parallel copies of the equator.
The two outermost circles are definite while all others are indefinite.
Using the deformation explained in~\cite{Williams1}*{Example~2} we can replace the definite folds by indefinite ones (see \cref{second_deformation1,second_deformation2}) and we obtain a wrinkled fibration~$f_1: S^1\times M\rightarrow S^2$ (which in particular is surjective, unlike~$f_0$).
As a next step we want to move the outermost folds of~$f_1$ into the equatorial region as indicated in \cref{second_deformation2,second_deformation3,second_deformation4}.
This is possible by the following result.
\begin{lemma}\label{lem_identification_complement}
Let~$D_1,D_2$ be two disks in the polar regions of the complements of~$\mc D(f_1)$ and let~$\gamma\subset S^2$ be the arc shown in \cref{second_deformation2}. 
There exist diffeomorphisms
	\begin{align*}
	\Phi : & ~(S^1\times M)\setminus \Int\bigl(f_1^{-1}(D_1)\amalg f_1^{-1}(D_2)\bigr) \to S^1\times f_1^{-1}(\gamma) \text{ and}\\
	\phi : & ~S^2\setminus \Int(D_1\amalg D_2) \to S^1\times I,  
	\end{align*}
such that $(\id \times h_1)\circ \Phi = \phi \circ f_1$, where~$h_1$ is the restriction~$f_1|_{f_1^{-1}(\gamma)}$.
\end{lemma}

\begin{proof}

Let~$B_i$ be a disk in~$S^2$ which contains~$D_i$ and either of the outermost images of folds, and is disjoint from the other images of folds of~$f_1$. 
It is easy to see by construction of $f_1$ that there exist desired diffeomorphisms $\Phi$ and $\phi$ defined on $(S^1\times M) \setminus \Int(f_1^{-1}(B_1)\amalg f_1^{-1}(B_2))$ and $S^2 \setminus \Int(B_1\amalg B_2)$, respectively. 
Next, we observe that the manifold $(S^1\times M) \setminus \Int(f_1^{-1}(D_1)\amalg f_1^{-1}(D_2))$ is obtained from~$(S^1\times M) \setminus \Int(f_1^{-1}(B_1)\amalg f_1^{-1}(B_2))$ by a round $1$-handle attachment (for a precise definition of round handles, see \cite{Baykur2}, for example). 
This round $1$-handle is untwisted and is attached with $0$-framing since the monodromy of the higher-genus side of the inner-most fold is trivial. 
This implies the statement in \cref{lem_identification_complement}. 
\end{proof}

\begin{figure}[htbp]
\begin{center}
\includegraphics[width=120mm]{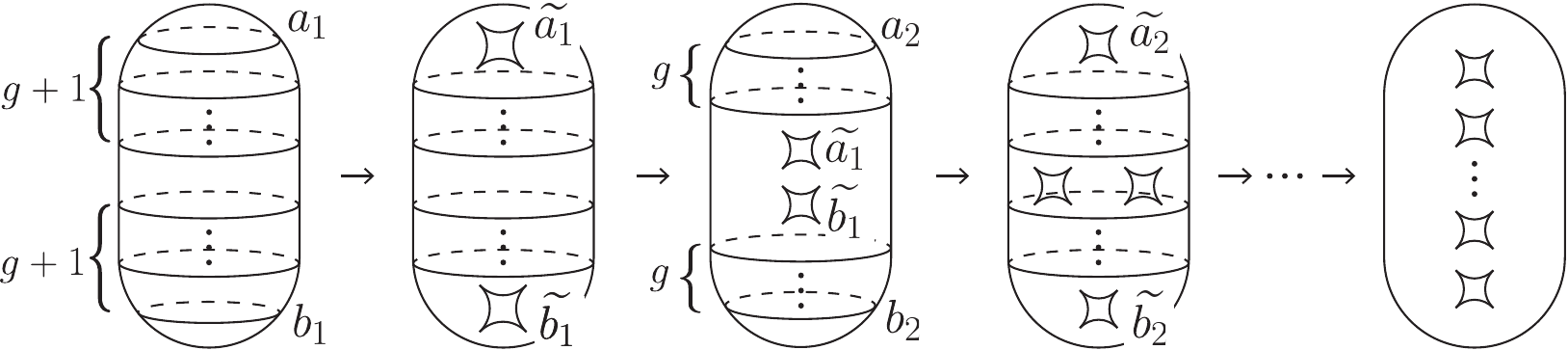}
\end{center}
\caption{Successive application of flip and slips. 
The far right figure is a base diagram of $f_3$. }
\label{flipslip_deformation}
\end{figure}

We denote the thus obtained map with critical values as in \cref{second_deformation4} by $f_2: S^1\times M \rightarrow S^2$ and label the images of the fold circles by $a_1,\ldots, a_{g+1}, b_1,\ldots, b_{g+1}$ as indicated. 
As shown in \cref{flipslip_deformation}, we next apply so called \emph{flip and slip} moves to the outermost folds~$a_1$ and~$b_1$ of~$f_2$ to obtain a new map having indefinite circles~$\widetilde{a_1}$ and~$\widetilde{b_1}$ with four cusps in place of~$a_1$ and~$a_2$.
These circles can be moved into the equatorial region and we can repeat this construction until all fold circles of~$f_2$ have been replaced by indefinite circles with four cusps as on the right of \cref{flipslip_deformation}.
We denote the resulting \wf by~$f_3$.
Finally, we can obtain a simple wrinkled fibration $w_h\colon S^1\times M\ra S^2$ by applying cusp merges to~$f_3$.
We will now explain how to perform the cusp merges in such a way that we are able to obtain a surface diagram of~$w_h$ from the information in the Heegaard diagram for~$M$.
First of all, we need to know surgered surface diagrams for all components~$\widetilde{a_i}$ and $\widetilde{b_i}$ of the critical value set of~$f_3$.
As explained in the previous paragraph, the fibration $f_3$ is obtained by successive application of flip and slips followed by deformations that move cusped circles to the equatorial region. 
Using the results of the second author about flip and slips~\cite{HayanoR2}*{Section~6}, it is therefore enough to understand surgered surface diagrams before the flip and slips.
The first flip and slips are applied to folds~$a_1$ and~$b_1$, whose lower-genus fibers are tori. 
According to \cite{HayanoR2}*{Theorem~6.7}, in order to obtain surgered surface diagrams of~$\widetilde{a_1}$ and~$\widetilde{b_1}$, we have to look at not only monodromies along the boundaries of the disks~$B_i$, $i=1,2$, which are trivial by \cref{lem_identification_complement}, but also sections over these disks.
According to \cref{lem_identification_complement} we can choose an
	\begin{equation*}
	S^1\times M\cong S^1\times f_1^{-1}(\gamma) \cup_{\varphi_1\amalg \varphi_2} (D_1\times T^2 \amalg D_2\times T^2)
	\end{equation*}
where $\varphi_i : \partial D_i\times T^2 \rightarrow S^1 \times \partial f^{-1}(\gamma)$ is a gluing diffeomorphism. 
By identifying~$\partial f_1^{-1}(\gamma)$ with a pair of tori we can regard~$\varphi_i$ as a self-diffeomorphism of~$S^1\times T^2$. 
Note that we can understand the behavior of the sections over~$B_1$ and~$B_2$ once we know the maps~$\varphi_1$ and~$\varphi_2$. 
Thus, the following lemma enables us to obtain surgered surface diagrams of~$\widetilde{a_1}$ and~$\widetilde{b_1}$.  

\begin{lemma}\label{lem_attachment}

For a suitable identification $\partial f_1^{-1}(\gamma)\cong T^2\amalg T^2$, there exists an embedding $\theta: S^1 \rightarrow T^2$ whose image intersects a vanishing cycle of the fold in one point such that the restriction $\varphi_i|_{S^1\times \{\ast\}}$ is isotopic to the map
\[
S^1\ni z \mapsto (z, \theta(z)) \in S^1\times T^2. 
\] 

\end{lemma}

\begin{proof}

Let $H_1$ be a solid torus. 
We take a small ball $B^3$ in the interior $\Int(H_1)$. 
By \cref{lem_identification_complement}, we can identify $f_1^{-1}(B_i\setminus \Int(D_i))$ with $S^1 \times (H_1 \setminus B^3)$. 
Thus, we can regard $f_1^{-1}(B_i)$ as the manifold obtained by attaching $D^2\times T^2$ to $S^1 \times (H_1 \setminus B^3)$. 
Manifold of this form were extensively studied in the proof of \cite{Hayanogenus1}*{Theorem 4.2}. 
It is easy to see from the construction of~$f_1$ that the manifold~$f_1^{-1}(B_i)$ is diffeomorphic to~$S^1\times D^3$.
Thus, by the arguments in~\cite{Hayanogenus1} the attaching map~$\varphi_i$ satisfies the desired condition. 
\end{proof}

The other flip and slips are applied to folds over~$a_i$ and~$b_i$ ($i \geq 2$) whose lower genus fibers have genus two or higher. 
Thus, in these cases, we can obtain surgered surface diagrams of~$\widetilde{a_i}$ and~$\widetilde{b_i}$ once we know the monodromies of the higher-genus sides of~$a_i$ and~$b_i$. 
These monodromies are no longer trivial, yet we can easily obtain them by looking at lifts of surgered diagrams of the previous folds~$\widetilde{a_{i-1}}$ and~$\widetilde{b_{i-1}}$ (see also the example below).

\begin{example}\label{ex_SD_3sphere}

Here we consider the genus~$0$ Heegaard diagram $(S^2,\emptyset,\emptyset)$ of~$S^3$.
We denote by~$f_i\colon S^1\times S^3\rightarrow S^2$,~$i=0,\dots,3$ the maps obtained in the various steps of the algorithm explained above.
It is easy to see that surgered surface diagrams of~$f_2$ are as described in \cref{surgered_diagram_sphere1}. 
\begin{figure}[htbp]
\begin{center}
\includegraphics[width=80mm]{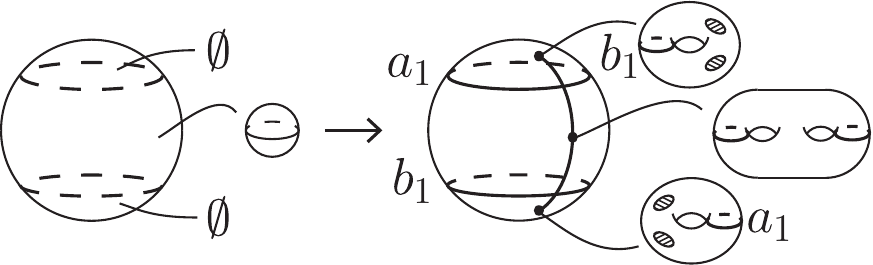}
\end{center}
\caption{Left: a base diagram of $f_0$. Right: that of $f_2$.  
}
\label{surgered_diagram_sphere1}
\end{figure}
We choose sections of~$f_2$ over the disks~$B_1$ and~$B_2$ as in the proof of \cref{lem_identification_complement}. 
These sections lift the monodromies of the higher-genus sides of the folds~$a_1$ and~$b_1$ to the mapping class group with a marked point~$x$. 
We can choose a homotopy from~$f_0$ to~$f_2$ and sections so that the monodromy along~$a_1$ (resp.~$b_1$) is equal to the pushing map along~$\theta_1$ (resp.~$\theta_2$), where~$\theta_i$ is the oriented loop described in \cref{surgered_diagram_sphere2a}. 
Denote by~$\Sigma_1$ a regular fiber of~$f_3$ which is mapped to the inside of~$\widetilde{a_1}$. 
We consider the simple closed curve $\theta_1^\pm\subset\Sigma_1$ as described in \cref{scc_for_flipslip} and put~$\varphi_1 = t_{\theta_1^+} t_{\theta_1^-}^{-1}$. 
By \cite{HayanoR2}*{Theorem~6.7} a surface diagram derived from $\widetilde{a_1}$ is $\mathfrak{S}_a= (\Sigma_1; d_1, \varphi_1^{-1}(\delta_1), a_1, \delta_1)$, where~$\delta_1$ and~$d_1$ are the simple closed curves in~$\Sigma_1$ described in \cref{surgered_diagram_sphere2b}. 
We can also verify in a similar manner that a surface diagram derived from the critical set~$\widetilde{b_1}$ is $\mathfrak{S}_b=(\Sigma_2; e_1,\epsilon_1, b_1, \varphi_2^{-1}(\epsilon_1))$, where we put $\varphi_2=t_{\theta_2^+} t_{\theta_2^-}^{-1}$. 

\begin{figure}[htbp]
\centering
\subfigure[The map $f_2$.]{\includegraphics[height=35mm]{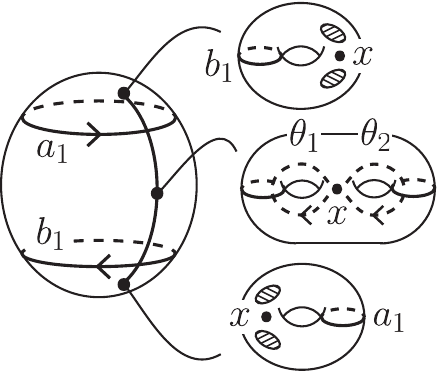}
\label{surgered_diagram_sphere2a}}
\hspace{.2em}
\subfigure[The map $f_3$.]{\includegraphics[height=35mm]{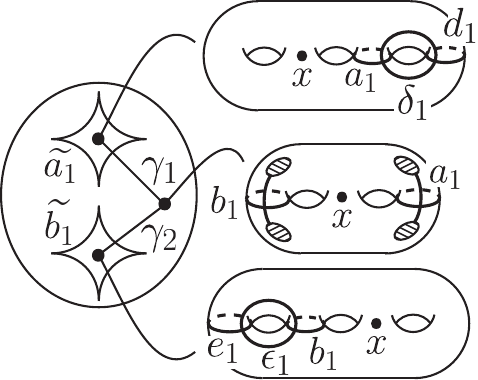}
\label{surgered_diagram_sphere2b}}
\hspace{.2em}
\subfigure[Curves in fibers.]{\includegraphics[height=35mm]{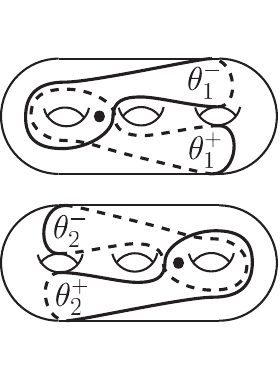}
\label{scc_for_flipslip}
}
\caption{Base diagrams and curves. }
\end{figure}

We consider arcs~$\gamma_1$ and~$\gamma_2$ in~$S^2$ as in \cref{surgered_diagram_sphere2b}. 
We denote by~$\Sigma_0$ the fiber on the intersection between~$\gamma_1$ and~$\gamma_2$. 
These paths give surgered surface diagrams which are described in \cref{surgered_diagram_sphere3}. 
\begin{figure}[htbp]
\centering
\subfigure[The diagrams $\mathfrak{S}_a$ and $\mathfrak{S}_b$.]{\includegraphics[height=35mm]{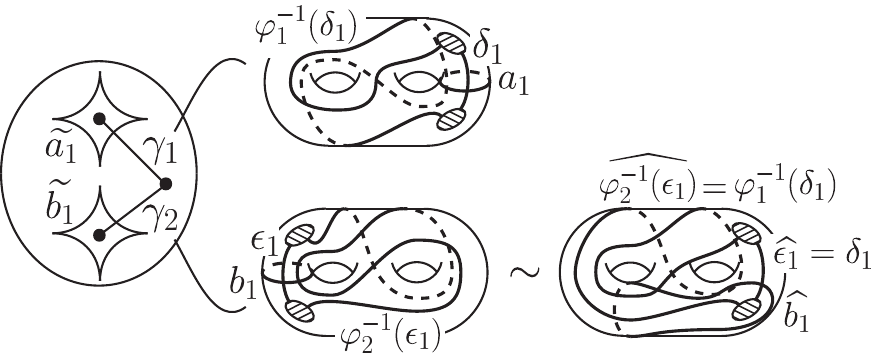}
\label{surgered_diagram_sphere3}}
\hspace{.2em}
\subfigure[A ``null-homotopic'' surface diagram of $S^1\times S^3$.]{\includegraphics[height=35mm]{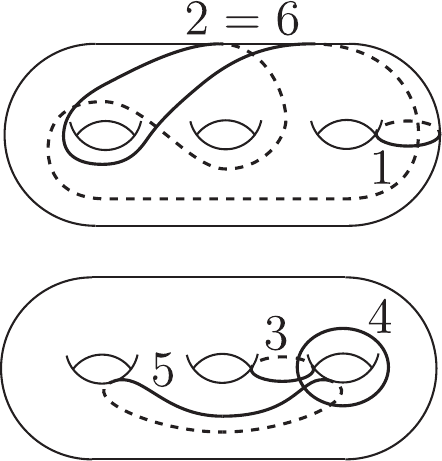}
\label{surface_diagram_sphere}}
\caption{(Surgered) surface diagrams.}
\end{figure}
We denote these diagrams by~$\mathfrak{S}_a^\prime$ and~$\mathfrak{S}_b^\prime$. 
To obtain a surface diagram of the fibration after applying a cusp merge we change the diagram~$\mathfrak{S}_b^\prime$ by a matching isotopy of~$\Sigma_0$ (as described in \cref{surgered_diagram_sphere3}) so that the surgery disks and the arc~$\epsilon_1$ in~$\mathfrak{S}_b^\prime$ coincide with the surgery disks and the arc~$\delta_1$ in~$\mathfrak{S}_a^\prime$.
We denote by~$\widehat{b_1}$, $\widehat{\epsilon_1}$ and~$\widehat{\varphi_2^{-1}(\epsilon_1)}$ the curves and arcs in~$\Sigma_0$ obtained after applying the above modification to~$b_1$, $\epsilon_1$ and~$\varphi_2^{-1}(\epsilon_1)$, respectively. 
We can regard~$\widetilde{b_1}$ and~$\widehat{\varphi_2^{-1}(\epsilon_1)}$ as simple closed curves in~$\Sigma_1$. 
Note that the curve~$\widehat{\varphi_2^{-1}(\epsilon_1)}$ is equal to~$\varphi_1^{-1}(\delta_1)$. 
We eventually obtain the following surface diagram $(\Sigma_1; d_1, \varphi_1^{-1}(\delta_1), a_1, \delta_1, \widehat{b_1}, \varphi_1^{-1}(\delta_1))$ of $S^1\times S^3$ which is shown in \cref{surface_diagram_sphere}. 

\end{example}

%% file: appendix_normalization_new.tex
\section{Truncating homotopies}
	\label{ch:normalization}
In this appendix we will give a way to make a given homotopy compactly supported, that is, we will find another homotopy which is constant outside a compact set and coincides with a given one in some open subset of the support. 
Although we merely need to find a compactly supported model of a cusp merge in this paper, the recipe we will give here can be applied in much more general situations (see \cref{E:making compactly supported}). 

For a subset $A \subset N\times J$ we denote the intersection $A \cap N\times \{-1\}$ by $A_{-1}$. 
Let $\pi:N\times J \to J$ be the projection. 

\begin{theorem}\label{T:normalization lemma}

Let $F = (f_t)$ be a homotopy of a map from $N$ to $P$, $K \subset P$ a closed set diffeomorphic to a $p$--ball and $R\subset \tilde{F}^{-1}(K\times J)$ a codimension $0$ compact submanifold which may have boundaries and corners.  
Suppose that the following conditions are satisfied: 
\begin{itemize}

\item there exists a diffeomorphism $\nu \Pa K \to \Pa K \times (-\delta, \delta)$, where $\nu \Pa K$ is a tubular neighborhood of $\Pa K$, such that $f_t$ is transverse to $K\times \{s\}$ and the restriction 
\[
f_t: f_t^{-1}(\Pa K \times \{s\}) \to \Pa K \times \{s\}
\] 
is locally stable for any $t$ and $s$ (see \cref{ch:stability} for the definition of local stability), 

\item the boundary $\Pa R$ is decomposed into three parts $\Pa_c R$, $\Pa_s R$ and $\Pa_e R$ so that 

\begin{itemize}

\item the set of corners of $R$ is the union of the intersections $\Pa_c R \cap \Pa_s R$, $\Pa_c R \cap \Pa_e R$ and $\Pa_s R \cap \Pa_e R$,

\item $\Pa_e R = R \cap \left(N \times \{-1, 1\}\right)$, 

\item $\Pa_c R$ is contained in $F^{-1}(\Pa K)$, 

\item the restriction $\wt{F}|_{\Pa_s R}:\Pa_s R \to K\times J$ is a surjective submersion. 

\end{itemize}

\end{itemize}
Then there exist two neighborhoods $R_1 \subset R_2$ of $\Pa R$ and a level preserving diffeomorphism $\Phi$ (resp.~$\phi$) of $N\times J$ (resp.~$P\times J$) such that the following conditions hold: 
\begin{itemize}

\item the support of $\Phi$ is contained in $R_2$ and that of $\phi$ is contained in $\nu \Pa K$, 

\item $f_{-1} = \phi_t^{-1} \circ f_t \circ \Phi_t$ on $R_1$, where $\Phi_t$ and $\phi_t$ are the maps satisfying the equations $\Phi(x,t) = (\Phi_t(x),t)$ and $\phi(y,t) = (\phi_t(y),t)$, respectively. 

\end{itemize} 

\end{theorem}

\begin{remark}\label{R:compactly supported model}

Using \cref{T:normalization lemma} we can obtain a compactly supported homotopy $g_t$ which coincides with $f_t$ in some open subset of the support as follows:
\[
g_t(x) = \begin{cases}
\phi_t \circ f_{-1} \circ \Phi_t^{-1} (x) & ((x,t) \in R_1\cup R^c) \\
f_{t}(x) & ((x,t) \in R)
\end{cases}
\]
 
\end{remark}

\begin{proof}[Proof of \cref{T:normalization lemma}]
Suppose that $\Pa_s R$ is not empty. 
The restriction $f_{-1}:\Pa_s R_{-1} \allowbreak \to K$ is a fiber bundle whose fiber $M$ is a closed manifold. 
Since $K$ is contractible $\Pa_s R_{-1}$ is isomorphic to $M\times K$ as a fiber bundle. 
The restriction $\pi: \Pa_s R \to J$ is also a fiber bundle with a fiber $\Pa_s R_{-1}$. 
Thus $\Pa_s R$ is isomorphic to $M \times K \times J$. 
Let $X_0$ be an everywhere non-zero vector field on a neighborhood of $\Pa_s R$ contained in $\ker(d \tilde{F})$. 
We can obtain the following bundle isomorphism by integrating $X_0$: 
\begin{equation*}
\Psi : \nu \Pa_s R \to M \times (-\varepsilon, \varepsilon) \times K\times {J},  
\end{equation*}
where $\varepsilon>0$ is sufficiently small and $\nu \Pa_s R$ is a neighborhood of $\Pa_s R$ in $\tilde{F}^{-1}(K\times J)$.
We take the above isomorphism so that $\Psi^{-1}(M \times [0, \varepsilon) \times K\times J)$ is in $R$. 
If $\Pa_s R$ is empty, then we do not need such an isomorphism. 
In what follows, we will take vector fields with several conditions assigned to bold symbols.
Starred conditions are related with the isomorphism $\Psi$, which are not needed if $\Pa_s R$ is empty. 

We denote the union $\Pa_c R\cup \Psi^{-1}(M \times (-\varepsilon, \varepsilon) \times \Pa K\times J)$ by $\nu_f \Pa_c R$, which is a neighborhood of $\Pa_c R$ in $F^{-1}(\Pa K)$. (Note that $\nu_f\Pa_c R = \Pa_c R$ if $\Pa_s R$ is empty.)
By the assumptions the restriction $\pi: \nu_f \Pa_c R \to J$ is a fiber bundle. 
We take a vector field $X_1$ on $\nu_f \Pa_c R$ which satisfies the following conditions: 
\begin{enumerate}[\bf A{\color{white}${}^\ast$}]

\item\label{C:vectorfieldX1_1}
 the image $d\pi(X_1)$ is equal to $\frac{d}{dt} \in TJ$, 

\item[\bf B${}^\ast$]\label{C:vectorfieldX1_2}
 the image $d\Psi(X_1)$ is equal to $\frac{d}{dt} \in T(M \times (-\varepsilon, \varepsilon) \times K\times J)$ on $\nu \Pa_s R\cap \nu_f \Pa_c R$, where $t$ is the coordinate of the interval $J$. 

\end{enumerate}
Integrating $X_1$ we take the following bundle isomorphism: 
\begin{equation*}
\Theta :\nu_f\Pa_cR \to \nu_f\Pa_c R_{-1} \times J. 
\end{equation*}
We denote a neighborhood of $\nu_f \Pa_c R$ in $F^{-1}(\nu \Pa K)$ by $\nu\Pa_cR$. 
Since $f_t$ is transverse to $\Pa K\times \{s\}$ for any $s\in (-\delta, \delta)$, the following composition is a submersion: 
\[
\nu \Pa_cR \xrightarrow{F} \nu \Pa K \cong \Pa K \times (-\delta, \delta) \xrightarrow{\mathrm{proj}} (-\delta, \delta).  
\] 
We denote the above composition by $p: \nu \Pa_c R \to (-\delta, \delta)$. 
We take a vector field $X_2$ on $\nu\Pa_cR$ so that it satisfies the following conditions: 
\begin{enumerate}[\bf A{\color{white}${}^\ast$}]
\setcounter{enumi}{2}

\item\label{C:vectorfieldX2_1} 
the image $dp(X_2)$ is equal to the vector field $\frac{d}{ds}$, 

\item[\bf D${}^\ast$]\label{C:vectorfieldX2_2}
 the image $d\Psi(X_2)$ is $\frac{d}{ds} \in T(M \times (-\varepsilon, \varepsilon)\times K\times J)$, where $s$ is the coordinate for the second component of $\nu \Pa K \cong \Pa K \times (-\delta, \delta)$. 

\end{enumerate}
Using $X_2$ and $\Theta$ we obtain the following bundle isomorphism: 
\[
\hat{\Theta}:\nu \Pa_c R \to \nu_f \Pa_c R_{-1} \times (-\delta, \delta) \times J. 
\]
It is easy to verify that the composition 
\[
\tilde{F}\circ \hat{\Theta}^{-1}:\nu_f \Pa_c R_{-1} \times (-\delta, \delta) \times J \to \Pa K \times (-\delta, \delta) \times J
\]
preserves the latter two components, that is, there exists a two parameter family $g_{s,t} \in C^\infty (\nu_f \Pa_c R_{-1}, \Pa K)$ such that $\tilde{F}\circ \hat{\Theta}^{-1}(x,s,t)$ is equal to $(g_{s,t}(x),s,t)$ for any $(x,s,t)$. 
We denote the composition $\tilde{F}\circ \hat{\Theta}^{-1}$ by $\hat{F}$. 
By the assumption $g_{s,t}$ is locally infinitesimally stable for any $t$ and $s$. 
Thus by the same argument as that in the proofs of Propositions 4.3, 4.5 and Theorem 4.6 in \cite{GG}, we can obtain a vector field $\xi$ on $\nu \Pa_c R_{-1} \times (-\delta, \delta) \times J$ and a vector field $\eta$ on $\Pa K \times (-\delta, \delta) \times J$ which satisfy the following conditions:
\begin{enumerate}[\bf A{\color{white}${}^\ast$}]
\setcounter{enumi}{4}

\item\label{C:vectorfieldThomLevine_1}
 for any $(x,s,t) \in \nu_f \Pa_c R_{-1}\times (-\delta, \delta)\times J$ the vector $\xi(x,s,t)$ is contained in $T\nu_f \Pa_c R_{-1}$, 

\item\label{C:vectorfieldThomLevine_2}
 for any $(y,s,t) \in \Pa K \times (-\delta,\delta) \times J$ the vector $\eta(y,s,t)$ is contained in $T\Pa K$, 

\item\label{C:vectorfieldThomLevine_3}
 $t\hat{F}(\xi) + \omega \hat{F} (\eta) = t\hat{F}\left(d\hat{\Theta}(\hat{\Theta}^{-1})^\ast \left(\frac{\Pa}{\Pa t}\right)\right) - \hat{F}^\ast \left(\frac{\Pa}{\Pa t}\right)$, 

\item[\bf H${}^\ast$]\label{C:vectorfieldThomLevine_4}
 for any $x \in \hat{\Theta}^{-1}(\nu \Pa_s R \cap \nu \Pa_c R)$ the image $d(\Psi \circ \hat{\Theta}^{-1})(\xi(x))$ is contained in $T K$. 

\end{enumerate}

Let $\varrho_1:(-\delta, \delta) \to \mb{R}$ be a non-negative smooth function whose value is $1$ on $\left[-\frac{\delta}{3}, \frac{\delta}{3}\right]$ and $0$ outside $\left[-\frac{2\delta}{3}, \frac{2\delta}{3}\right]$. 
Since we identify $\nu \Pa K$ with $\Pa K \times (-\delta, \delta)$ we can regard $\varrho_1$ as a function defined on $\nu \Pa K$, and this function can be extended to that defined on  a neighborhood of $P\times J$ so that it vanishes on the complement of $\nu \Pa K\times J$. 
We also denote the extended function by $\varrho_1$, and let $\tilde{\eta} = \varrho_1 \eta$.   
Since $\tilde{\eta}$ is compactly supported the vector field $\tilde{\eta} + \frac{\Pa}{\Pa t}$ on $P \times J$ is complete. 
We denote the integral curve of it with the initial point $y\in P\times J$ by $c_{y}(t)$. 
We define a self-diffeomorphism $\phi$ of $P\times J$ as follows: 
\[
\phi(p,t) = c_{(p,-1)}(1+t). 
\]
We can deduce from the condition \ref{C:vectorfieldThomLevine_2} that the diffeomorphism $\phi$ preserves $\nu \Pa K$. 

We denote the set 
\[
\hat{\Theta}^{-1}\left(\left(\Psi^{-1}\left(M \times \left[-\frac{\varepsilon}{n},\frac{\varepsilon}{n} \right]\times \Pa K \times \{-1\}\right)\cup \nu_f \Pa_c R_{-1}\right) \times \left[-\frac{\delta}{n}, \frac{\delta}{n}\right]\times J\right)
\]
by $W_{1,n}$, and the set $\Psi^{-1}\left(M \times \left[-\frac{\varepsilon}{n}, \frac{\varepsilon}{n}\right] \times K\times J \right)$ by $W_{2,n}$. 
We take a vector field $\tilde{\xi}$ on $\nu \Pa_c R \cup \nu \Pa_s R$ which satisfies the following properties: 

\begin{enumerate}[\bf A{\color{white}${}^\ast$}]
\setcounter{enumi}{8}

\item\label{C:vectorfieldglobal_1} 
the vector field $\tilde{\xi}$ is equal to $d\hat{\Theta}^{-1}\hat{\Theta}^\ast(\xi)$ on $ W_{1,3}$, 

\item[\bf J${}^\ast$]\label{C:vectorfieldglobal_2}
for $x = \Psi^{-1}(z,r,q,t)\in W_{2,3}$ the vector $\tilde{\xi}(x)$ is equal to 
\[
\varrho_1(q)d\hat{\Theta}^{-1}(\xi(\hat{\Theta}(x))) \allowbreak + (1-\varrho_1(q))d\Psi^{-1}\left(d\tilde{F}_x \left(\frac{\Pa}{\Pa t}\right)_x +  d \phi_{\phi^{-1}(q,t)} \left(\frac{\Pa}{\Pa t}\right)_{\phi^{-1}(q,t)}\right),
\]
where the vectors $d\tilde{F}_x \left(\frac{\Pa}{\Pa t}\right)_x$ and $d \phi_{\phi^{-1}(q,t)} \left(\frac{\Pa}{\Pa t}\right)_{\phi^{-1}(q,t)}$ on $K\times J$ are regarded as those on $M \times (-\varepsilon, \varepsilon)\times K\times J$ in the obvious way.  

\end{enumerate}
We also take a function $\varrho_2 : N\times J \to \mathbb{R}$ so that $\varrho_2 \equiv 1$ on $W_{1,6} \cup W_{2,6}$ and $\varrho_2 \equiv 0$ outside $W_{1,4} \cup W_{2,4}$. 
Since the vector field $\varrho_2\tilde{\xi}$ is compactly supported, the vector field $-\varrho_2\tilde{\xi} + \frac{\Pa}{\Pa t}$ on $N\times \mb{R}$ is complete. 
We denote the integral curve of $-\varrho_2\tilde{\xi} + \frac{\Pa}{\Pa t}$ with the initial point $x\in \tilde{N}$ by $C_{x}(t)$. 
We define a self-diffeomorphism $\Phi$ of $N\times J$ as follows: 
\[
\Phi(x) = \begin{cases}
C_{\hat{\Theta}(y_{-1},s,t)}(1+t) & (x = \hat{\Theta}^{-1}(y,s,t)\in W_{1,3}) \\
C_{\Psi(z,r,q_{-1},t)}(1+t) & (x = \Psi^{-1}(z,r,q,t) \in W_{2,3}) \\
x & (\text{otherwise}), 
\end{cases}
\]
where $y_{-1}\in \psi^{-1}\left(\Pa(\Pa_cR)\times \left[-\frac{\varepsilon}{3}, \frac{\varepsilon}{3}\right]\right)$ is $\Theta^{-1}(p_1\circ \Theta(y),-1)$ and $q_{-1} \in K\times J$ is $(p_1(q), -1)$ ($p_1$ represents the projection onto the first component).  
Let $R_1$ and $R_2$ denote the unions $W_{1,6}\cup W_{2,6}$ and $W_{1,3}\cup W_{2,3}$, respectively.
The support of $\Phi$ is contained in $R_2$. 
Furthermore it is easy to see that $\Phi$ preserves the sets $W_{1,6}$ and $W_{2,6}$. 
The diffeomorphisms $\Phi$ and $\phi$ preserve the time components. 
We put $\Phi(x,t) = (\Phi_t(x), t)$ and $\phi(p,t) = (\phi_t(p),t)$. 

We denote the composition $\phi^{-1} \circ \tilde{F} \circ \Phi$ by $\Omega$. 
By the same argument as that in the proof of Theorem 3.3 in \cite{GG}, we can obtain the following equation: 
\[
d\Omega \left(\frac{\Pa}{\Pa t}\right) - \Omega^\ast \left(\frac{\Pa}{\Pa t}\right) = 0 \text{ on }R_1, 
\]
We can deduce from this equation that $\phi_t^{-1}\circ f_t \circ \Phi_t$ is equal to $f_{-1}$ on $R_1$ (see \cite{GG}*{Lemma~3.2}). 
\end{proof}

\begin{example}\label{E:making compactly supported}

Here we consider the cusp merge model $F = (\mu_s)$ defined as follows: 
\[
\mu_s:\mb{R}^4\to \mb{R}^2, \hspace{.5em} \mu_s(t,x,y,z) = (t, x^3 -3(t^2+s)x + y^2 -z^2). 
\]
For $\varepsilon >0$, we denote by $K_\varepsilon$ a closed disk neighborhood of the image of the joining curve obtained by smoothing the corners of the rectangle $[-1-\varepsilon, 1+\varepsilon]\times [-\varepsilon, \varepsilon]$. 
Assume that $K_\varepsilon$ is sufficiently close to the rectangle so that $\mu_s$ is transverse to the boundary $\Pa K_\varepsilon$. 
We take a tubular neighborhood $\nu \Pa K_\varepsilon \cong \Pa K_\varepsilon \times (-\delta, \delta)$ of $\Pa K_\varepsilon$ so that the map $\mu^{-1}(\nu \Pa K_\varepsilon)\xrightarrow{\mu_s} \nu \Pa K_\varepsilon \xrightarrow{p} \Pa K_\varepsilon$ has only Morse singularities and thus locally infinitesimally stable for any $s \in [-1,\varepsilon]$, where $p:\nu \Pa K_\varepsilon \to \Pa K_\varepsilon$ is the projection. 
Let $(\mathcal{H}_s)$ be a smooth one parameter family of distributions such that $\mathcal{H}_s$ is a connection of $\mu_s$ for each $s$. 
It is easy to verify that the fiber $\mu_s^{-1}(t,0)$ is a disk if $ t^2\leq -s$ and a once punctured torus if $t^2 >-s$. 
We take a two parameter family $(S_{s,t})$ of circles such that $S_{s,t}$ is contained in $\mu_s^{-1}(t,0)$ and $S_{s,t}$ bounds a compact subsurface containing $(t,0,0,0)$ in $\mu_s^{-1}(t,0)$. 
We can obtain a three dimensional submanifold $\Pa_s R\subset \tilde{F}^{-1}(K_\epsilon \times [-1, \varepsilon])$ by taking a flow from $(S_{s,t})$ with respect to $\mc{H}_s$.    
The complement of $\Pa_s R \cup \tilde{F}^{-1}(\Pa K_\varepsilon\times [-1,\varepsilon])$ has a compact component which contains the origin of $\mb{R}^4 \times \mb{R}$. 
We denote this component by $R$ and let $\Pa_e R = \Pa R \cap \mb{R}^4 \times \{-1,\varepsilon\}$ and $\Pa_c R = \overline{\Pa R \setminus (\Pa_s R\cup \Pa_e R)}$. 
It is easy to see that $K$ and $R$ satisfy the conditions in \cref{T:normalization lemma}. 
Thus we can obtain a compactly supported homotopy which coincides with $f_s$ in an open subset of its support. 

Note that the above construction of a compactly supported homotopy is also valid in much more general situation. 
Indeed, the same recipe can be applied to a homotopy $g_s$ if we have a disk $K$ such that the restriction of $g_s$ over $\Pa K$ is a Morse function and $\Crit (\tilde{G})$ is compact. 
In particular, using \cref{T:normalization lemma} we can obtain a compactly supported model for each versal unfolding appearing in \cref{T:1-parameter versal unfoldings}. 

\end{example}

%% file: appendix_symmetries.tex
\section{Symmetries of folds and cusps}
	\label{ch:symmetries}
	\newcommand{\InvD}[1]{\mathrm{Inv}^\sim(\mc D(#1))}
The purpose of this appendix is to explain how to deal with ambiguities that arise when working with local models mentioned.
Let~$f\colon N\ra P$ be a smooth map and suppose that the germ of~$f$ at~$x\in N$ has a local model~$\mf f\colon (\R^n,0)\ra(\R^p,0)$, that is, there are coordinates around~$x$ and~$f(x)$ in which~$f$ is represented by~$\mf f$.
If we choose a second set of coordinates with this property, then the germs of the two coordinate sets differ by germs of local diffeomorphisms $\rho$ of~$\R^n$ and~$\lambda$ of~$\R^p$, respectively, which preserve the germ of~$\mf f$ at the origin, that is, we have
	\begin{equation*}
	\mf f = \lambda\circ \mf f\circ \rho\inv.
	\end{equation*}
This leads us to study the group of pairs of such local diffeomorphisms.
\subsection*{Symmetries of map-germs}
More systematically, we denote by~$\mc E(n,p)$ the set of germs of smooth maps $(\R^n,0)\ra(\R^p,0)$ and let~$\mc R$ and~$\mc L$ be the groups of local diffeomorphism-germs defined near the origins of~$\R^n$ and~$\R^p$, respectively; that is, 
	\begin{equation*}
	\mc R=\Set{\rho\in\mc E(n,n)}{\det(d\rho_0)\neq0}
	\eqand
	\mc L=\Set{\lambda\in\mc E(p,p)}{\det(d\lambda_0)\neq0}.
	\end{equation*}
%
Writing $\mc A=\mc R\times \mc L$ for the product we get a split exact sequence
	\begin{equation}\label{eq:RAL sequence}
	1\lra \mc R  \overset{i}{\lra} \mc A \overset{p}{\lra} \mc L \lra 1
	\end{equation}
where~$i(\rho)=(\rho,\id)$ and~$p(\rho,\lambda)=\lambda$.
The group $\mc A$ naturally acts on~$\mc E(n,p)$ by the so called \emph{right-left action}
	$(\rho,\lambda)\cdot\mf f=\lambda\circ\mf f\circ\rho\inv$
which restricts to the \emph{right action} of~$\mc R$ and the \emph{left action} of~$\mc L$, whence the names.
The above considerations about local models lead us to study the \emph{right-left symmetries} or $\mc A$--symmetries of a given germ~$\mf f\in\mc E(n,p)$ which are defined as the stabilizer
	\begin{equation*}
	\mc A(\mf f)=\Set{(\rho,\lambda)\in\mc A}{\lambda\circ\mf f\circ\rho\inv=\mf f}.
	\end{equation*}
Of course, there are similar definitions involving~$\mc R$ and~$\mc L$.
As a last bit of notation, we denote groups of orientation preserving diffeomorphisms and symmetries by adding a superscript as in $\mc A^+$, $\mc A^+(\mf f)$, etc.
For $n\geq p$ the groups of $\mc A$--symmetries are typically quite large.
But it turns out that for sufficiently well-behaved germs 
-- more precisely for \emph{finitely $\mc A$--determined} germs which includes the stable ones -- 
they share a convenient property with finite dimensional Lie groups: 
according to Jänich~\cite{Jaenich} and Wall~\cite{Wall}, they have ``maximal compact Lie subgroups'' which are unique up to conjugacy and the quotient by such a group is ``contractible''.
However, as the quotation marks indicate, this statement should be taken with a grain of salt.
In fact, $\mc A$ and its subgroups are not topological groups in any natural way.%
	\footnote{This is due to the more general lack of a natural topologies on sets of smooth map-germs.}
In any case, these delicacies shall not concern us.
The interested reader is referred to the exposition of du~Plessis and Feragen in~\cite{duPlessis_Feragen} which contains the most recent and complete account on these matters.
In order to convince the reader that these are in fact the right notions we state the following immediate consequence of the definitions (that we did not give).
\begin{fact}\label{T:}
If~$G_{\mc A}(\mf f)\subset\mc A(\mf f)$ is a maximal compact subgroup, then any element of~$\mc A(\mf f)$ can be deformed to an element of~$G_{\mc A}(\mf f)$ by a smooth family in~$\mc A(\mf f)$ 
(where a family of germs~$\mf f_s\in\mc E(n,p)$, $s\in I$, is called smooth if there exists a smooth map $F\colon (I\times \R^n, I\times{0})\ra (\R^p,0)$ such that $F(s,\cdot)$ represents~$\mf f_s$).
\end{fact}
In order to determine maximal compact subgroups of $\mc A$--symmetries for the fold and cusp singularities, we will use results of du~Plessis--Wilson~\cite{duPlessis_Wilson} and du~Plessis--Feragen~\cite{duPlessis_Feragen}.
These provide an analogue of the exact sequence~\eqref{eq:RAL sequence} for maximal compact subgroups and, in addition, completely determine the topology of maximal compact subgroups of~$\mc R(\mf f)$.
We begin with the latter.
\begin{theorem}[du~Plessis--Wilson~\cite{duPlessis_Wilson}]\label{T:right symmetries}
Let~$\mf f\in \mc E(n,p)$ be a stable germ.
Then the group~$\mc R(\mf f)$ of right symmetries of~$\mf f$ has a maximal compact subgroup~$G_{\mc R}(\mf f)$ which is unique up to conjugation.
Moreover, if~$n\geq p$, then 
	\begin{equation}\label{eq:right symmetries}
	G_{\mc R}(\mf f)\cong
	\begin{cases}
	\mathrm{O}(n-p)							& \text{if $\rk(d\mf{f}|_0)=p$} \\
	\mathrm{O}(i)\times\mathrm{O}(r-i)	& \text{if $\rk(d\mf{f}|_0)=p-1$} \\
	\{1\}											& \text{if $\rk(d\mf{f}|_0)\leq p-2$} 
	\end{cases}
	\end{equation}
where in the second case~$r$ and~$i$ denote the rank and index of the intrinsic second derivative~$\delta^2|_0\mf{f}\colon \ker(d\mf f|_0)\times\ker(d\mf f|_0)\ra\coker(d\mf f|_0)$.
\end{theorem}
Before we can state the exact sequence involving maximal compact subgroups of $\mc R(\mf f)$ and~$\mc A(\mf f)$, we have to understand the image of~$\mc A(\mf f)$ in~$\mc L$ for stable~$\mf f$.
It is easy to see that if~$(\rho,\lambda)\in\mc A(\mf f)$, then $\lambda\in\mc L$ satisfies the following two conditions.
\begin{itemize}
	\item \textit{The discriminant condition:}
		$\lambda$~preserves the (germ of) the discriminant $\mc D(\mf f)=f(\mc C(f))$.
	\item \textit{The cokernel condition:}
		If $\rk d\mf f_0=p-1$ and $r\neq 2i$, then $d\lambda_0$ preserves the orientation of~$\coker(d\mf f_0)$.
\end{itemize}
These conditions define a subgroup of~$\mc L$ which is denoted by~$\InvD{\mf f}$ in~\cite{duPlessis_Wilson}.
According to~\cite{duPlessis_Wilson}*{Lemma~4.1}, $\mc A(\mf f)$ is mapped onto~$\InvD{\mf f}$ if $\mf f$~is stable.
%
%
%
\begin{theorem}[du~Plessis--Feragen~\cite{duPlessis_Feragen}]\label{T:du Plessis-Feragen}
Let~$\mf f\in\mc E(n,p)$ be a stable germ such that any linear compact subgroup of~$\InvD{\mf f}$ is contained in the image of a compact subgroup of~$\mc A(\mf f)$.
Then there are maximal compact subgroups 
	$G_{\mc R}(\mf f)\subset\mc R(\mf f)$, 
	$G_{\mc A}(\mf f)\subset\mc A(\mf f)$, and 
	$G_{\mc D}(\mf f)\subset \InvD{\mf f}\cap\GL(p)$
such that \eqref{eq:RAL sequence} restricts to a split exact sequence
	\begin{equation}\label{eq:RAL sequence_compact}
	1\lra G_{\mc R}(\mf f)  \overset{i}{\lra} G_{\mc A}(\mf f) \overset{p}{\lra} G_{\mc D}(\mf f) \lra 1.
	\end{equation}
\end{theorem}
\begin{remark}\label{R:du Plessis-Wilson error}
A stronger form of this theorem was originally stated by du~Plessis and Wilson in~\cite{duPlessis_Wilson}*{(1.5)} but an error in their proof was pointed out in~\cite{duPlessis_Feragen}*{Ch.~4.2}.
\end{remark}

\subsection*{Symmetries of folds and cusps}
Now we turn to the computation of maximal compact subgroups of $\mc A$--symmetries of the fold and cusp singularities for~$n=4$ and~$p=2$.
Recall that these were defined in terms of the local models
	\begin{align*}
	\mf F_\pm(t,x,y,z)&=(t, x^2+y^2\pm z^2) \\
	\mf C_\pm(t,x,y,z)&=(t, x^3+3tx+y^2\pm z^2).
	\end{align*}
From the formulas some immediate linear symmetries are apparent.
In fact, the linear $\mc R$--symmetries are easily determined as
	\begin{equation*}
	\begin{array}{cc}
	\mc R(\mf F_+)\cap\GL(4) \cong \O(3) 
	\hspace{10pt}&
	\mc R(\mf F_-)\cap\GL(4) \cong \O(2,1)
	\\[5pt]
	\mc R(\mf C_+)\cap\GL(4) \cong \O(2)
	\hspace{10pt}&
	\mc R(\mf C_-)\cap\GL(4) \cong \O(1,1)
	\end{array}
	\end{equation*}
where $\O(3)$ and~$\O(2,1)$ act on the $(x,y,z)$--coordinates while $\O(2)$ and~$\O(1,1)$ act on the $(y,z)$--plane.
(Note that $\O(r,s)$ is not compact for $r,s\geq1$, a maximal compact subgroup is given by $\O(r)\times\O(s)=\O(r,s)\cap\O(r+s)$.)
Moreover, there are some additional linear $\mc A$--symmetries, one for each $\mf F_\pm$, given by
	\begin{equation}\label{eq:special fold symmetry}
	(t,x,y,z)\mapsto(-t,x,y,z)
	\eqand
	(u,v)\mapsto(-u, v),
	\end{equation}
and another one for~$\mc C_-$, namely
	\begin{equation}\label{eq:special cusp symmetry}
	(t,x,y,z)\mapsto(t,-x,z,y)
	\eqand
	(u,v)\mapsto(u, -v).
	\end{equation}
Obviously, these have order two and therefore contribute copies of $\Z_2$ in~$\mc A(\mf F_\pm)$ and~$\mc A(\mf C_-)$.
We will show that the maximal compact subgroups of linear symmetries are also maximal in the full groups of $\mc A$--symmetries.
\begin{proposition}\label{T:symmetries of folds and cusps}
The $\mc A$--symmetries of the fold and cusp singularities have maximal compact subgroups of the form
	\begin{equation*}
	\begin{array}{cc}
	G_{\mc A}(\mf F_+) \cong \O(3) \ltimes \Z_2
	\hspace{10pt}&
	G_{\mc A}(\mf F_-) \cong \big(\O(2)\times\O(1)\big) \ltimes \Z_2
	\\[5pt]
	G_{\mc A}(\mf C_+) \cong \O(2)
	\hspace{10pt}&
	G_{\mc A}(\mf C_-) \cong \big(\O(1)\times\O(1)\big) \ltimes \Z_2
	\end{array}
	\end{equation*}
generated by the linear symmetries exhibited above.
In particular, for orientation preserving $\mc A$--symmetries we have maximal compact subgroups
	\begin{equation*}
	\begin{array}{cc}
	G^+_{\mc A}(\mf F_+) \cong \SO(3)
	\hspace{10pt}&
	G^+_{\mc A}(\mf F_-) \cong \O(2)
	\\[5pt]
	G^+_{\mc A}(\mf C_+) \cong \SO(2)
	\hspace{10pt}&
	G^+_{\mc A}(\mf C_-) \cong \Z_2
	\end{array}
	\end{equation*}
where $\O(2)$ is embedded into $\O(2)\times\O(1)$ via $A\mapsto(A,\det A)$ and the $\Z_2$ inside~$\O(1)\times\O(1)$ is generated by~$(t,x,y,z)\mapsto(t,x,-y,-z)$.
\end{proposition}
\begin{proof}
Let $\mf f$~be either $\mf F_\pm$ or~$\mf C_\pm$.
One readily checks that 
	\begin{equation*}
	\ker(d\mf f|_0)=\scp{\del_x,\del_y,\del_z}
	\eqand
	\im(d\mf f|_0)=\scp{\del_u=d\mf f|_0(\del_t)}
	\end{equation*}
so that $\coker(d\mf f|_0)$ can be identified with the span of~$\del_v$.
Furthermore, it is simply a matter of writing out the definitions to see that in these identifications the intrinsic second derivatives $\delta^2_0\mf F_\pm$ and~$\delta^2_0\mf C_\pm$ correspond to the Hessians of the functions~$x^2+y^2\pm z^2$ and~$x^3+y^2\pm z ^2$, respectively.
So according to \cref{T:right symmetries}, the maximal compact subgroups of $\mc R$--symmetries satisfy
	\begin{equation*}
	\begin{array}{cc}
	G_{\mc R}(\mf F_+) \cong \O(3)
	\hspace{10pt}&
	G_{\mc R}(\mf F_-) \cong \O(2)\times\O(1)
	\\[5pt]
	G_{\mc R}(\mf C_+) \cong \O(2)
	\hspace{10pt}&
	G_{\mc R}(\mf C_-) \cong \O(1)\times\O(1)
	\end{array}
	\end{equation*}
and these groups are realized by the linear symmetries discussed above.
In order to apply \cref{T:du Plessis-Feragen} we have to study the compact subgroups of $\InvD{\mf f}\cap\GL(2)$ which can be determined completely in the cases at hand.
Indeed, it is easy to see that the only non-trivial linear diffeomorphism that preserves the discriminants of the cusps
	\begin{equation*}
	\mc D(\mf C_\pm)=\Set{(-\tau^2,-2\tau^3)}{\tau\in\R}
	\end{equation*}
is given by~$(u,v)\mapsto(u,-v)$ which violates to cokernel condition for~$\mf C_+$ whereas the cokernel is empty for $\mf C_-$.
Therefore, we have
	\begin{equation*}
	\InvD{\mf C_+}\cap\GL(2)=1
	\eqand
	\InvD{\mf C_-}\cap\GL(2)=\Z_2.
	\end{equation*}
and similar reasoning yields
	\begin{equation*}
	\InvD{\mf F_\pm}\cap\GL(2)=\Set{\left(\begin{smallmatrix}
	a&b\\0&c
	\end{smallmatrix}\right)}{ac\neq0,\;c>0}
	\end{equation*}
which has a unique compact subgroup, a copy of~$\Z_2$ generated by $(u,v)\mapsto(-u,v)$.
Finally, according to~\eqref{eq:special cusp symmetry} and~\eqref{eq:special fold symmetry} any compact subgroup of $\InvD{\mf f}\cap\GL(2)$ can be lifted to~$\mc A(\mf f)$ and the claims about the structure of~$\mc G_{\mc A}(\mf f)$ now follow from \cref{T:du Plessis-Feragen}.
\end{proof}